\documentclass[12pt]{elsarticle}
\usepackage[utf8]{inputenc}

\usepackage[margin=1in]{geometry}
\usepackage{amsmath}  
\usepackage{amsthm}
\usepackage{amsfonts}
\usepackage{amssymb}
\usepackage{graphicx}  
\usepackage{float}
\usepackage{natbib}
\usepackage{hyperref}
\usepackage{bm}

\usepackage{caption}
\usepackage{subcaption}

\newtheorem{proposition}{Proposition}

\theoremstyle{definition}

\newtheorem{remark}{Remark}

\newcommand{\projS}[1]{\bm{\Pi}^S_{k,E} #1}

\newcommand{\projLtwo}[2][\ell]{\bm{\Pi}^{0}_{#1,E} \bm{\varepsilon} (#2)}

\newcommand{\Pk}[1][k]{\mathbb{P}_{#1}}
\newcommand{\grad}[1]{\nabla #1}
\newcommand{\diverge}[1]{\nabla \cdot #1}
\newcommand{\polysym}{{\mathbb{P}}_{\ell}(E)^{2\times 2}_{\text{sym}}}
\newcommand{\VEM}{\bm{V}_{k,\ell}^E}
\newcommand{\symgrad}[1]{\nabla_s #1}
\newcommand{\pder}[1]{\frac{\partial }{\partial #1}}
\newcommand{\largeset}{\mathcal{EN}_{k,\ell}^E}

\newcommand{\projLtwodisp}[2][k]{\bm{\Pi}^{0}_{#1,E}#2}
\newcommand{\Nmatrix}{\bm{N}^{\partial E}}

\newcommand{\vm}[1]{\bm{#1}}
\newcommand{\vx}{\bm{x}}

\newcommand{\ac}[1]{{#1}}

\hypersetup{
    colorlinks,
    citecolor=black,
    filecolor=black,
    linkcolor=black,
    urlcolor=black
}

\begin{document}
\begin{frontmatter}

\title{Stabilization-free serendipity virtual element method for plane elasticity}

\cortext[cor1]{Corresponding authors}

\author[inst1]{Alvin Chen\corref{cor1}}
\ead{avnchen@ucdavis.edu}

\affiliation[inst1]{organization={Department of Mathematics},
            addressline={University of California}, 
            city={Davis},
            postcode={95616}, 
            state={CA},
            country={USA}}

\author[inst2]{N. Sukumar\corref{cor1}}
\ead{nsukumar@ucdavis.edu}
\affiliation[inst2]{organization={Department of Civil and Environmental Engineering},
            addressline={University of California}, 
            city={Davis},
            postcode={95616}, 
            state={CA},
            country={USA}}

\begin{abstract}
We present a higher order stabilization-free virtual element method applied to plane elasticity problems. We utilize a serendipity approach to reduce the total number of degrees of freedom from the corresponding high-order approximations. The well-posedness of the problem is numerically studied via an eigenanalysis. The method is then applied to several benchmark problems from linear elasticity and we show that the method delivers optimal convergence rates in $L^2$ \ac{norm} and energy seminorm that match theoretical 
estimates as well as the convergence rates from higher order virtual element methods. 
\end{abstract}


\end{frontmatter}

\section{Introduction} The \ac{Virtual Element Method} (VEM) is an extension of the classical finite element method (FEM) to arbitrary polygonal meshes. In early works~\cite{AHMAD2013376,Artoli:2017:cmech,basicprinciple,elasticdaveiga,Veiga2014TheHG,daveiga2014virtual,Gain:2014:VEM}, the method was developed for scalar elliptic boundary-value problems as well as extended to linear elasticity. These methods when applied to high-order VEM required additional internal degrees of freedom; however, from the study of serendipity finite elements~\cite{Arnold:2011:fcm}, it was shown that the number of internal degrees of freedom can be greatly reduced. \ac{This property is also important for the extension to three-dimensional VEM. For polyhedral elements, the face degrees of freedoms cannot be reduced using standard static condensation techniques, but the serendipity approach can still be applied to reduce these face degrees of freedom~\cite{veiga:2018:cham}.} In~\cite{Veiga:2016:caf,veiga:2018:cham} and~\cite{debellis:2019:cas}, the serendipity virtual element method was developed for scalar problems and nonlinear elasticity, respectively. Unlike serendipity FEM, the serendipity VEM is found to be robust for general polygonal meshes, including on meshes with distorted elements. In all prior studies with the VEM, a stabilization term is required in order to ensure 
that the element stiffness matrix has the correct rank. Recently in~\cite{berrone2021lowest}, a stabilization-free VEM was introduced for the Poisson equation, and extended to linear plane elasticity in~\cite{chen:2022:arxiv}. The main idea in this approach is to modify the approximation space to be able to compute a higher order polynomial $L^2$ projection of the strain. A secondary projection operator is used to fill in the additional degrees of freedom introduced by the higher order polynomials. A similar approach was explored in~\cite{enhanced:VEM}; however, instead of using additional projections, they used static condensation to eliminate the extra degrees of freedom. In the standard VEM, the stabilization term is problem dependent and there is no general method of constructing it; therefore, devising a stabilization-free virtual element method is desirable. 

In this paper, we combine the serendipity elements~\cite{Veiga:2016:caf} to extend the stabilization-free techniques from~\cite{berrone2021lowest,chen:2022:arxiv} to higher order methods for two-dimensional linear elasticity. The resulting method will not have a stabilization term and in many cases will not require any additional internal degrees of freedom. In Section~\ref{sec:elastostatics}, we introduce the model problem of plane elasticity. In Sections~\ref{sec:mathprelim} and~\ref{sec:projections} we set up the necessary polynomial spaces and recall some properties of the serendipity VEM. In Section~\ref{sec:vemspace}, we define the higher order VEM spaces from~\cite{berrone2021lowest,chen:2022:arxiv} and in Section~\ref{sec:numerical_implementation} we discuss the construction and implementation of the projections, element stiffness and forcing terms. Section~\ref{sec:choice_ell} contains a numerical study of an upper bound for the order of polynomial enhancement for second- and third-order methods. In Section~\ref{sec:numerical_results}, we apply the second- and third-order methods to the patch test, two manufactured problems, a beam under a sinusoidal load and an infinite plate with circular hole under uniaxial tension. The rate of convergence for the manufactured problems and the beam problem agree with theoretical estimates; however, consistent with expectations~\cite{Artioli:2020:cmame,daveiga:2019:esiam}, the presence of the curved boundary in the infinite plate with circular hole problem results in suboptimal convergence rates.  We close 
with some concluding remarks in Section~\ref{sec:conclusions}.

\section{Elastostatic Model Problem and Weak Form}\label{sec:elastostatics}
We consider an elastic body that occupies the region $\Omega \subset \mathbb{R}^2$ with
boundary $\partial \Omega$. Assume that the boundary $\partial \Omega$ can be written as the disjoint union of two parts $\Gamma_D$ and $\Gamma_N$ with prescribed Dirichlet and Neumann conditions on $\Gamma_D$ and $\Gamma_N$, respectively. 
The strong form for the elastostatic problem is:
\begin{subequations}\label{strongproblem}
\begin{align}
&\diverge{\bm{\sigma}} + \bm{f} = \bm{0}  \ \ \textrm{in } \Omega, 
\quad \bm{\sigma} = \bm{\sigma}^T \ \ \textrm{in } \Omega, \\
&\bm{\varepsilon}(\bm{u}) =\symgrad{\bm{u}} = \frac{1}{2}(\grad{\bm{u}} + \grad{\bm{u}}^T),\\
&\bm{\sigma}(\bm{u}) = \mathbb{C} : \bm{\varepsilon}(\bm{u}), \\ 
&\bm{u} = \ac{{\bm{u}_0}} \quad \text{on } \Gamma_D ,\\
&\bm{\sigma} \cdot \mathbf{n} = \ac{{\bm{t}_0}} \quad \text{on } \Gamma_N,
\end{align}
\end{subequations}
where $\bm{f} \in [L^2(\Omega)]^2$  is the body force per unit volume, $\bm{\sigma}$ is the Cauchy stress tensor,
$\bm{\varepsilon}$ is the small-strain tensor with $\symgrad{(\cdot)}$ being the symmetric gradient operator, 
$\mathbf{u}$ is the displacement field, 
$\ac{\bm{u}_0}$ and $\ac{\bm{t}_0}$ are the imposed essential boundary and traction boundary data, and 
$\bm{n}$ is the unit outward normal on the boundary. Linear elastic constitutive material relation ($\mathbb{C}$ is the material moduli tensor) and small-strain kinematics are assumed. 

The associated weak form of the boundary-value problem posed in~\eqref{strongproblem} is to find the displacement field $\bm{u} \in \bm{U}$, where $\bm{U} := \{\bm{u} : \bm{u} \in [H^{1}(\Omega)]^2, \ 
\bm{u} = \ac{\bm{u}_0} \ \textrm{on } \Gamma_D \}$, such that
\begin{subequations}\label{weakbilinearform}
\begin{align}
    a(\bm{u},\bm{v}) &= b(\bm{v}) \quad \forall \bm{v} \in \bm{U}_0,
\intertext{where $\bm{U}_0 = [H^1_0(\Omega)]^2$ and}
    a(\bm{u},\bm{v}) &= \int_{\Omega} \bm{\sigma}(\bm{u}): \bm{\varepsilon}(\bm{v}) \, d \bm{x}, \\ 
    b(\bm{v}) &= \int_{\Omega} \bm{f} \cdot \bm{v} \, d \bm{x} + \int_{\Gamma_N} \ac{\bm{t}_0}\cdot \bm{v} \, ds .
\end{align}
\end{subequations}
In~\eqref{weakbilinearform}, $H^1(\Omega)$ is the Hilbert space that consists of square-integrable functions up to order $1$ and $H_0^1(\Omega)$ is the subspace of
$H^1(\Omega)$ that contains functions that vanish on $\ac{\Gamma_D}$.

\section{Mathematical preliminaries}\label{sec:mathprelim}
Let $\mathcal{T}^h$ be the decomposition of the region $\Omega$ into nonoverlapping polygons with standard mesh assumptions~\cite{basicprinciple}. For each polygon $E \in \mathcal{T}^h$, we denote its diameter by $h_E$ and its centroid by $\bm{x}_E$. Each
polygon $E$ consists of $N_E$ vertices (nodes) with $N_E$ edges. 
Let the coordinate of each vertex be $\bm{x}_i := (x_i,y_i)$.
We denote the
$i$-th edge by $e_i$ and let the set $\mathcal{E}_E$ represent the collection of all edges of $E$. 

\subsection{Polynomial basis}
Over each element $E$, we define $[{\mathbb{P}_k}(E)]^2$ 
as the space of of two-dimensional vector-valued polynomials of degree less than or equal to $k$. On each $E$, we choose the scaled monomial vectorial basis set
as: 
\begin{subequations}\label{eq:Mhat}
\begin{align}
    \bm{\widehat{M}}(E) = \begin{bmatrix}  
    \begin{Bmatrix}
    1 \\ 0 
    \end{Bmatrix}, 
    \begin{Bmatrix}
    0 \\ 1
    \end{Bmatrix}, 
    \begin{Bmatrix}
    -\eta \\ \xi
    \end{Bmatrix}, 
    \begin{Bmatrix}
    \eta \\ \xi
    \end{Bmatrix}, 
    \begin{Bmatrix}
    \xi \\ 0
    \end{Bmatrix}, 
    \begin{Bmatrix}
    0 \\ \xi
    \end{Bmatrix},
    \dots ,
    \begin{Bmatrix}
        \eta^k \\ 0
    \end{Bmatrix},
    \begin{Bmatrix}
        0 \\ \eta^k
    \end{Bmatrix}
    \end{bmatrix},\label{poly_2_basis}
\end{align}
where 
\begin{align}
    \xi = \frac{x-x_E}{h_E} , \quad \eta = \frac{y-y_E}{h_E} .
\end{align}
The $\alpha$-th element of the set $\bm{\widehat{M}}(E)$ is denoted by $\bm{m}_\alpha$, and \ac{we} define the matrix $\tilde{\vm{N}}^p$ that contains the basis elements as 
\begin{align}
    \tilde{\vm{N}}^p := \begin{bmatrix}
    1 & 0 & -\eta & \eta & \xi & 0 & \dots & \eta^k & 0 \\
    0 & 1 & \xi & -\xi & 0 & \xi & \dots & 0 & \eta^k 
    \end{bmatrix}.
\end{align}
\end{subequations}

\smallskip
We also define the space 
$\polysym$ that represents $2\times 2$ symmetric matrix polynomials of degree less than or equal to $\ell$. 
\ac{We adopt Voigt notation to represent symmetric
$2\times 2$ matrices as an equivalent $3 \times 1$ vector. Let $\vm{A}$ be a symmetric $2\times 2$ matrix whose Voigt representation is $\overline{\vm{A}}$: 
\begin{equation*}
    \vm{A} = \begin{bmatrix}
    a_{11} & a_{12} \\
    a_{12} & a_{22}
    \end{bmatrix}, \quad
    \overline{\vm{A}} = \begin{Bmatrix}
        a_{11} \\a_{22} \\ a_{12}
    \end{Bmatrix}. 
\end{equation*}
}
On using Voigt notation for $\polysym$, the basis set over element $E$ is written as: 
\begin{subequations}\label{eq:M2x2hat}
\begin{align}
    \widehat{\bm{M}}^{2\times 2}(E)=
    \begin{bmatrix}
    \begin{Bmatrix}
        1 \\ 0 \\ 0 
    \end{Bmatrix},
    \begin{Bmatrix}
        0 \\ 1 \\ 0
    \end{Bmatrix},
    \begin{Bmatrix}
        0 \\ 0 \\ 1
    \end{Bmatrix},
    \begin{Bmatrix}
        \xi \\ 0 \\ 0
    \end{Bmatrix},
    \begin{Bmatrix}
        0 \\ \xi \\ 0
    \end{Bmatrix},
    \begin{Bmatrix}
        0 \\ 0 \\ \xi
    \end{Bmatrix},\dots,
    \begin{Bmatrix}
        \eta^\ell \\ 0 \\ 0
    \end{Bmatrix},
    \begin{Bmatrix}
        0 \\ \eta^\ell \\0 
    \end{Bmatrix},
    \begin{Bmatrix}
        0 \\ 0 \\ \eta^\ell
    \end{Bmatrix}
    \end{bmatrix}.
\end{align}
We denote the $\alpha$-th vector in this set as $\widehat{\bm{m}}_\alpha$ and define the matrix $\bm{N}^p$ that
contains these basis elements as 
\begin{align}
    \bm{N}^p := \begin{bmatrix}
        1 & 0 & 0 & \xi & 0 & 0 &\dots & \eta^\ell & 0 & 0 \\
        0 & 1 & 0 & 0 &\xi & 0 &\dots & 0 & \eta^\ell & 0 \\
        0 & 0 & 1 & 0 & 0 & \xi & \dots& 0 & 0 & \eta^\ell
    \end{bmatrix}.\label{matrix_N_p}
\end{align}
\end{subequations}

\ac{
On using Voigt notation, the stress and strain tensor are represented as:
\begin{equation*}
\overline{\bm{\sigma}} = \begin{Bmatrix}
        \sigma_{11} \\ \sigma_{22} \\ \sigma_{12}. 
    \end{Bmatrix}, \quad
    \overline{\bm{\varepsilon}} = \begin{Bmatrix}
       \varepsilon_{11} \\ \varepsilon_{22} \\ 2\varepsilon_{12}
  \end{Bmatrix}.
  \end{equation*}
Now with the vector representation of stress and strain, we can write the strain-displacement relation and the constitutive law in matrix form as:
\begin{subequations}
\begin{align}
    \overline{\bm{\sigma}} &= \bm{C}\overline{\bm{\varepsilon}}, \quad
    \overline{\bm{\varepsilon}} = \bm{S}\bm{u}, \\
\intertext{where $\bm{S}$ is a matrix differential operator that is given by}
    \bm{S} &= \begin{bmatrix}
        \pder{x} & 0 \\ 0 & \pder{y} \\ \pder{y} & \pder{x}
    \end{bmatrix},
\end{align}
\end{subequations}
and $\bm{C}$ is the associated matrix representation of the material tensor that is given by
\begin{align*}
    \bm{C} &= \frac{E_Y}{(1-\nu^2)} \begin{bmatrix}
    1 & \nu & 0 \\ \nu & 1 & 0 \\ 0 & 0 & \frac{1-\nu}{2} 
    \end{bmatrix} \qquad \qquad \qquad \qquad  (\textrm{plane stress}),\\
    \bm{C} &= \frac{E_Y}{(1+\nu)(1-2\nu)} \begin{bmatrix}
    1-\nu & \nu & 0 \\ \nu & 1-\nu & 0 \\ 0 & 0 & \frac{1-2\nu}{2}
    \end{bmatrix} \quad (\textrm{plane strain}),
\end{align*}
where $E_Y$ is the 
Young's modulus and $\nu$ is the Poisson's ratio of the material. 
}

\subsection{Properties of serendipity virtual elements}
We recall some results on serendipity virtual element methods for scalar problems from~\cite{Veiga:2016:caf}.
Let $E$ be a polygon with $N_E$ edges and let $\eta_E$ be the minimum number of unique lines to cover $\ac{\partial E}$. For a $k$-th order method there are a total of $kN_E$ boundary degrees of freedom and $\frac{k(k-1)}{2}$ internal degrees of freedom. The idea of serendipity VEM is that we are able to fully define a computable projection operator by only 
retaining a subset of all the degrees of freedom. To do this, we introduce two propositions as proven in~\cite{Veiga:2016:caf}.
\begin{proposition}\label{proposition:1}
For $k<\eta_E$, if the set of $S$ degrees of freedom $\{\delta_1,\delta_2,\dots \delta_S\}$ contains all of the $kN_E$ boundary degrees of freedom, then the following property holds true:
\begin{align}
    \delta_1(p_k)=\delta_2(p_k) = \dots = \delta_S(p_k)=0 \implies p_k \equiv 0 \quad \forall p_k \in \Pk(E), \label{eq:seren_key_prop}
\end{align}
where $\delta_i(\cdot)$ is the $i$-th degree of freedom of its argument.
\end{proposition}
\begin{proposition}\label{proposition:2}
For $k\geq \eta_E$, if the set of $S$ degrees of freedom
$\{\delta_1,\delta_2,\dots, \delta_S\}$ contains all $kN_E$ boundary degrees of freedom and contain all internal moments of order $\leq k-\eta_E$, then the set satisfies 
\begin{align}
    \delta_1(p_k)=\delta_2(p_k) = \dots = \delta_S(p_k)=0 \implies p_k \equiv 0 \quad \forall p_k \in \Pk(E).
\end{align}
\end{proposition}
 Once these degrees of freedom are chosen, we can construct a serendipity projection operator $\projS$ such that it satisfies the properties:
 \begin{subequations}
 \begin{align}
     &\projS \text{ can be fully computed using } \delta_1, \delta_2, \dots, \delta_S, \\ 
     \intertext{and}
     &\projS{p_k} = p_k \quad \forall p_k\in \Pk(E).
 \end{align}
 \end{subequations}
 This operator is used to define a serendipity virtual element space for a vector field and to construct two $L^2$ operators. 
 
 \begin{remark}
For $k=2$, on any polygonal element $E$ it is sufficient to take $\{\delta_1,\delta_2,\dots, \delta_S\}$ to be the vertex and edge degrees of freedom. For $k=3$, if $E$ is at least a quadrilateral with four distinct sides then it is also sufficient to take $\{\delta_1,\delta_2,\dots, \delta_S\}$ as the vertex and edge degrees of freedom.
\end{remark}
\section{Projection operators}\label{sec:projections}
We first define the serendipity projection of the displacement field as proposed in~\cite{Veiga:2016:caf}.
We then present the derivation of two $L^2$ projection operators, the $L^2$ projection of the displacement and the $L^2$ projection of the strain. 
\subsection{Serendipity projection}
\ac{For any element $E$, denote $\bm{H}^1(E):= [H^1(E)]^2$ and $\bm{C}^0(\bar{E}):= [C^0(\bar{E})]^2$}. Let $S$ be the number of sufficient degrees of freedom for a scalar function as defined in Proposition~\ref{proposition:1} and~\ref{proposition:2}, and then define the operator $D: \vm{H}^1(E) \ac{\cap \bm{C}^0(\bar{E})} \to \mathbb{R}^{2S}$ by 
\begin{align}
    D(\vm{v}) = \bigl( \delta_1(\vm{v}),\delta_2(\vm{v}), \dots, \delta_{2S}(\vm{v}) \bigr),
\end{align}
where $\delta_i(\vm{v})$ is the $i$-th degree of freedom of the vector field $\vm{v}$. We define the serendipity projection operator $\projS : \vm{H}^1(E)\ac{\cap \bm{C}^0(\bar{E})} \to [\mathbb{P}_{k}(E)]^2$ as the unique function that satisfies the orthogonality condition: 
\begin{subequations}\label{eq:seren_proj}
\begin{align}
    \left( D(\projS{\vm{v}}-\vm{v}) , D(\vm{m}_\alpha) \right)_{\mathbb{R}^{2S}} = 0 \quad \forall \vm{m}_\alpha\in [\Pk(E)]^2.\label{eq:seren_proj_a} \\
\intertext{On writing out the expressions, we get the equivalent system: }
    \sum_{j=1}^{2S}{\delta_{j}(\projS{\vm{v}})\delta_j({\vm{m}_\alpha}}) = \sum_{j=1}^{2S}{\delta_{j}({\vm{v}})\delta_{j}({\vm{m}_\alpha})} \quad \forall \vm{m}_\alpha\in [\Pk(E)]^2.\label{eq:seren_proj_b}
\end{align}
\end{subequations}

\subsection{\texorpdfstring{L\textsuperscript{2}}{L2} projection of the displacement field}
We define the $L^2$ projection operator $\projLtwodisp : \vm{H}^1(E) \to [\Pk(E)]^2$ of the displacement field by the function that satisfies the $L^2$ orthogonality relation: 
\begin{subequations}
\begin{align}
        (\vm{p},\vm{v}-\projLtwodisp\vm{v})_{E} &= 0 \quad \forall \vm{p} \in [\Pk(E)]^2,\label{eq:L2disp_a}\\
\intertext{where we use the standard $L^2$ inner product for vector fields:}
    (\vm{p},\vm{v})_E &=\int_{E}{\vm{p} \cdot \vm{v} \,d\vm{x}}.\label{eq:L2disp_b}\\
\intertext{Expanding~\eqref{eq:L2disp_a} and rewriting in matrix-vector operations, we have}
    \int_{E}{\vm{p}^T \projLtwodisp{\vm{v}}\,d\vm{x}} &= \int_{E}{\vm{p}^T \vm{v}\, d\vm{x}} \quad \forall \vm{p} \in [\Pk(E)]^2.\label{eq:L2disp_c}
\end{align}
\end{subequations}
\subsection{\texorpdfstring{L\textsuperscript{2}}{L2} projection of the strain field}
We define the associated $L^2$ projection operator $\projLtwo{.} :\bm{H}^1(E) \to \polysym$ of the strain tensor by the unique operator that satisfies 
\begin{subequations}\label{defL2proj}
\begin{align}
    \label{defL2proj-a}
    (\bm{\varepsilon}^p,\bm{\varepsilon}(\bm{v})- \projLtwo{\bm{v}})_E &= 0 \quad \forall \bm{\varepsilon}^p \in \polysym , \\
    \intertext{where we use the $L^2$ inner product (double contraction) for rank-2 tensor fields:}
    \label{defL2proj-b}
    (\bm{\varepsilon}^p,\bm{\varepsilon})_E &= \int_{E}{\bm{\varepsilon}^p:\bm{\varepsilon} \,d\bm{x}}.
\end{align}
\end{subequations}
After using Voigt notation and simplifying, we obtain the system for the $L^2$ strain projection as 
\begin{subequations}
\begin{align}
    \int_{E}{\left(\overline{\projLtwo{\bm{v}}}\right)^T \overline{\bm{\varepsilon}^p}\,d\bm{x}} &= \int_{\partial E}{\bm{v}^T \ac{\Nmatrix}\overline{\bm{\varepsilon}^p} \,ds} + \int_{E}{\bm{v}^T\bm{\partial} \overline{\bm{\varepsilon}^p}\,d\bm{x}} \quad \forall  {\bm{\varepsilon}^p} \in \polysym, \label{eq:L2_matrix_eq} \\ 
\intertext{where}
    &\qquad  \ac{\Nmatrix} := \begin{bmatrix}
        n_1 & 0 & n_2 \\
        0 & n_2 & n_1 
    \end{bmatrix}\label{normal_matrix},\\
    &\qquad  \bm{\partial} :=\begin{bmatrix}
    \pder{x} & 0 &\pder{y} \\
    0 & \pder{y} & \pder{x}
    \end{bmatrix},
\end{align}
\end{subequations}
and $\overline{\bm{\varepsilon}^p}$, $\overline{\projLtwo{\bm{v}}}$ are the Voigt representations of $\bm{\varepsilon}^p$ and $\projLtwo{\bm{v}}$, respectively.
\ac{\begin{remark}
The integral defined on the right-hand side of~\eqref{eq:L2disp_c} and the last integral on the right-hand side of~\eqref{eq:L2_matrix_eq} are not computable from the degrees of freedom of a standard $k$-th order virtual element space. In Section~\ref{sec:vemspace}, we construct an enlarged VEM space where it is possible to compute these integrals. 
\end{remark}}

\ac{\begin{remark}
For the implementation of the projection operators, 
we use the monomial basis given in~\eqref{eq:Mhat} and~\eqref{eq:M2x2hat}. 
\end{remark}}

\section{Enlarged Enhanced Serendipity Virtual Element Space}\label{sec:vemspace}
For any element $E$, fix a $\ell = \ell(E)$ and define the set $\largeset$ as 
\begin{equation}
\begin{split}
   \largeset := \left\{ \bm{v} \in \vm{H}^1(E)\ac{\cap \bm{C}^0(\bar{E})} :  
    \int_{E}{\bm{v}\cdot \bm{p}\,d\bm{x}} = \int_{E}{\projS{\bm{v}}\cdot \bm{p} \,d\bm{x} }\right.  \\ 
   \left. \forall \bm{p} \in \ac{[\mathbb{P}_{\ell-1}(E)]^2\slash [\mathbb{P}_{k-\eta_E} (E)]^2 }\right\}, \label{enrichment_condition}
\end{split}
\end{equation}
\ac{where $[\mathbb{P}_{\ell-1}(E)]^2\slash [\mathbb{P}_{k-\eta_E} (E)]^2$ denotes the set of vector polynomials in $[\mathbb{P}_{\ell-1}(E)]^2$ that are orthogonal to $[\mathbb{P}_{k-\eta_E}(E)]^2$ with respect to the $L^2$ inner product on $E$.} We then define the local enlarged virtual element space as:
\begin{align}
    \bm{V}_{k,\ell}^{E} := \left\{ \bm{v}_h \in \largeset : \Delta \bm{v}_h \in [\mathbb{P}_{\ell-1}(E)]^2 , \ \gamma^{e}(\bm{v}_h) \in [\mathbb{P}_{k}(e)]^2  \ \forall e \in \mathcal{E}_E , \ \bm{v}_h \in [C^{0}(\partial E)]^2\right\} \!, 
    \label{local_vem} 
\end{align}
where $\gamma^{e_i}(\cdot)$ is the trace of a function (its argument) on an edge $e_i$.
In the above space we require functions to be $k$-th order vector polynomials on the edges, and by the serendipity condition we take the degrees of freedom to be the values of the function at the vertices and edges of the polygon $E$ and possibly all the internal moments up to order $k-\eta_E$. In general, there are a total of 
\[ 2S = \max{\{2kN_E,2kN_E+(k-\eta_E+1)(k-\eta_E+2)\}} \]
degrees of freedom.
With the local space defined, we define the global enlarged virtual element space as
\begin{align}
    \bm{V}_{k,\bm{\ell}} := \{\bm{v}_h \in [H^1(\Omega)]^2 : \bm{v}_h\rvert_E \in \bm{V}_{k,\ell}^{E} \ \ \text{for } \ell=\ell(E)\}.
\end{align}
For each $E$, we assign a suitable basis to the local virtual element space $\VEM$. Let $\{\phi_{i}\}$ be the set of \ac{canonical basis functions~\cite{basicprinciple,Veiga2014TheHG}} that satisfy
$\delta_j(\phi_i) = \delta_{ij}$, 
\ac{where $\delta_{ij}$ is the Kronecker-delta.}
We express the components of  any  $\bm{v}_h \in \VEM$ as the sum of these basis functions:
\begin{subequations}\label{vem_basis_matrix}
\begin{align}
\bm{v}_h &= \begin{Bmatrix}
    v_h^1 \\ v_h^2
\end{Bmatrix}   
=\begin{bmatrix}
    \phi_1 & \phi_2 & \dots &\phi_{S} & 0 &0 & \dots & 0 \\
    0 & 0 & \dots & 0 & \phi_1 & \phi_2 &\dots & \phi_{S}
\end{bmatrix}\begin{Bmatrix}
    v^1_1 \\ v^1_2 \\ \vdots \\ v^2_{S}
\end{Bmatrix}
:= \bm{N}^v \tilde{\bm{v}}_h, \\
\intertext{where we define $\bm{N}^v$ as the matrix of vectorial basis functions:}
    \bm{N}^v &= \begin{bmatrix}
        \phi_1 & \phi_2 & \dots &\phi_{S} & 0 &0 & \dots & 0 \\
    0 & 0 & \dots & 0 & \phi_1 & \phi_2 &\dots & \phi_{S}
    \end{bmatrix} 
    := \begin{bmatrix}
        \bm{\varphi}_1
        & \dots & \bm{\varphi}_{S} 
        & \dots & \bm{\varphi}_{2S} 
    \end{bmatrix}.
\end{align}
\end{subequations}

We now define the weak form of the virtual element method on this space.  On defining a discrete bilinear operator $a^E_h : \VEM \times \VEM \to \mathbb{R} $ and a discrete linear functional $b^{E}_h: \VEM \to \mathbb{R}$, we seek the solution to the problem: find $\bm{u}_h \in \VEM$ such that
\begin{align}\label{vem_weak_problem}
a_h^E(\bm{u}_h,\bm{v}_h) = b^{E}_h(\bm{v}_h) \quad \forall \bm{v}_h \in \VEM .
\end{align}
Following~\cite{berrone2021lowest}, we introduce the local discrete bilinear form in matrix-vector form:
\begin{align}
    a_h^E(\bm{u}_h,\bm{v}_h) := \int_{E}{\left(\overline{\projLtwo{\bm{v_h}}}\right)^T \bm{C} \, \overline{\projLtwo{\bm{u}_h}} \, d\bm{x}},\label{local_bilinearform}
\end{align}
with the associated global operator defined as 
\begin{align}\label{global_bilinearform}
    a_h(\bm{u}_h,\bm{v}_h) := \sum_{E}{a_h^E(\bm{u}_h,\bm{v}_h)}.
\end{align}
We also define a local linear functional by 
\begin{align}
    b^E_h(\bm{v}_h)= \int_{E}{\bm{v}_h^T\bm{f}_h \, d\bm{x}} + \int_{\Gamma_N \cap \partial E}{\bm{v}_h^T\ac{\bm{t}_0} \, ds},\label{local_force}
\end{align}
with the associated global functional 
\begin{align}
    b_h(\bm{v}_h) = \sum_{E}{b^E_h(\bm{v}_h)},
\end{align}
where $\vm{f}_h$ is some approximation to $\vm{f}$. For a $k$-th order method we use $\vm{f}_h = \projLtwodisp{\vm{f}}$, but from~\cite{Artoli:2017:cmech}, it is sufficient to take $\vm{f}_h = \projLtwodisp[k-2]{\vm{f}}$.
\section{Numerical implementation}\label{sec:numerical_implementation}
\ac{For simplicity of implementation, we only consider the case $k<\eta_E$
(meshes do not contain triangles for $k = 3$).  This removes the need for internal moment degrees of freedom in the construction of the serendipity projection and simplifies the space $\largeset$.}
\subsection{Implementation of serendipity projector}
We start with the implementation of the serendipity projector. From~\eqref{eq:seren_proj_b}, we have for $\alpha=1,2, \dots , N_k$, where
$N_k = \textrm{dim} ([\mathbb{P}_k(E)]^2) = (k+1)(k+2)$, the condition
\begin{align*}
    \sum_{j=1}^{2S}{\delta_{j}(\projS{\vm{v}_h})\delta_j({\vm{m}_\alpha}}) = \sum_{j=1}^{2S}{\delta_{j}({\vm{v}_h})\delta_{j}({\vm{m}_\alpha})}  .
\end{align*}
We choose $\vm{v}_h=\vm{\varphi}_i$, the basis functions in $\VEM$, and expand $\projS{\vm{\varphi}_i}$ in terms of the scaled monomial basis functions:
\begin{align}
    \projS{\vm{\varphi}_i} = \sum_{\beta=1}^{N_k}{s^i_\beta \vm{m}_\beta} ,
\end{align}
where $\vm{m}_\alpha$ is an element of $\widehat{\vm{M}}(E)$ in~\eqref{eq:Mhat}.
Expanding the left-hand side of~\eqref{eq:seren_proj_b}, we have 
\begin{align}
     \sum_{j=1}^{2S}{\delta_{j}(\projS{\vm{v}_h})\delta_j({\vm{m}_\alpha}}) =  \sum_{\beta=1}^{N_k}s^i_\beta\sum_{j=1}^{2S}{{\delta_{j}(\vm{m}_\beta)}\delta_j({\vm{m}_\alpha}}).
\end{align}
Define the matrix $\vm{\widehat{G}}$ ($\alpha, \beta =1,2,\dots,N_k)$ by 
\begin{subequations}
\begin{align}
    \widehat{\vm{G}}_{\alpha \beta} &=  \sum_{j=1}^{2S}{{\delta_{j}(\vm{m}_\beta)}\delta_j({\vm{m}_\alpha}}).\label{eq:G_seren}
\intertext{Similarly, we define the matrix $\widehat{\vm{B}}$ representing the right-hand side of~\eqref{eq:seren_proj_b} by }
    \widehat{\vm{B}}_{\alpha i} &= \sum_{j=1}^{2S}{\delta_j({\vm{\varphi}_i}})\delta_j({\vm{m}_\alpha}).\label{eq:B_seren}
\intertext{Now combining these linear equations we can determine the coefficients $\{s^i_\beta\}$ for the serendipity projection by solving the linear system:}
    &\vm{\Pi}^S =\widehat{\vm{G}}^{-1} \widehat{\vm{B}}, 
\end{align}
\end{subequations}
where $(\vm{\Pi}^S)_{\beta i} = s^i_{\beta}$ is the matrix representation of the serendipity projection operator in the scaled monomial vectorial basis set. 

\ac{\begin{remark}
To compute the matrix $\widehat{\vm{G}}$,  it is convenient to use 
\begin{align*}
    \widehat{\vm{G}} = \vm{D}^T\vm{D},
\end{align*}
where $\bm{D}$ is the $2S\times N_k$  matrix that is defined by
\begin{align*}
    \vm{D}_{j \alpha}:= \delta_{j}(\bm{m}_\alpha) \quad ( j=1,2,\dots 2S, \ \alpha=1,2,\dots N_k).
\end{align*}.  
\end{remark}}

\subsection{Implementation of the \texorpdfstring{L\textsuperscript{2}}{L2} displacement projector}
With the serendipity projection matrix on hand, we now construct the remaining projection matrices. We start with the construction of the $L^2$ projection operator of the displacement field. From~\eqref{eq:L2disp_c}, we have the relation 
\begin{align}\label{eq:L2dispproj}
    \int_E{\vm{p}^T\projLtwodisp{\vm{v}_h} \, d\vm{x}} = \int_E{\vm{p}^T\vm{v}_h \,d\vm{x}}.
\end{align}
Expanding $\vm{v}_h$ in terms of the basis in $\VEM$, we have $\vm{v}_h = \vm{N}^v\tilde{\vm{v}}_h$. Similarly we expand $\vm{p}$ and $\projLtwodisp{\vm{v}_h}$ in terms of the polynomial basis in $[\Pk(E)]^2$. In particular we obtain $\vm{p} = \tilde{\vm{N}}^p\tilde{\vm{p}}$ and $\projLtwodisp{\vm{v}_h} = \tilde{\vm{N}}^p\tilde{\vm{\Pi}}^0\tilde{\vm{v}}_h$. On substituting into~\eqref{eq:L2disp_c} and simplifying, we obtain
\begin{align}\label{eq:L2dispproj-2nd}
    \tilde{\vm{p}}^T\left(\int_E{(\tilde{\vm{N}}^p)^T\tilde{\vm{N}}^p\,d\vm{x}}\right)\tilde{\vm{\Pi}}^0\tilde{\vm{v}}_h = \tilde{\vm{p}}^T\left(\int_E{(\tilde{\vm{N}}^p)^T\vm{N}^v\, d\vm{x}}\right)\tilde{\vm{v}}_h.
\end{align}
Define the matrix 
\begin{subequations}
\begin{align}
    \tilde{\vm{G}} &:= \int_E{(\tilde{\vm{N}}^p)^T\tilde{\vm{N}}^p\,d\vm{x}}.\\
    \intertext{\ac{The integral on the right-hand side 
     of~\eqref{eq:L2dispproj-2nd} is not computable directly; 
     however by applying the property of the space~\eqref{enrichment_condition}, we can realize an equivalent computable matrix:} }
    \tilde{\vm{B}} &:= \int_E{(\tilde{\vm{N}}^p)^T\projS{\vm{N}^v}\, d\vm{x}}.\label{eq:Btilde_matrix}
\end{align}
\end{subequations}
Then, we solve for the projection matrix $\tilde{\vm{\Pi}}^0$ in terms of the matrices $\tilde{\vm{G}}$ and $\tilde{\vm{B}}$:
\begin{align*}
    \tilde{\vm{\Pi}}^0 = \tilde{\vm{G}}^{-1}\tilde{\vm{B}}.
\end{align*}
This projection matrix is used to compute the element force integral which appears in~\eqref{eq:loc_force}.
\subsection{Implementation of the \texorpdfstring{L\textsuperscript{2}}{L2} strain projector}
To compute the $L^2$ projection of the strain we follow the construction in~\cite{chen:2022:arxiv}. Expand $\vm{v}_h = \vm{N}^v\tilde{\vm{v}}_h$ , $\overline{\vm{\varepsilon}^p} =\vm{N}^p\tilde{\vm{\varepsilon}}^p $ and $\overline{\projLtwo{\vm{v}_h}} = \vm{N}^p\vm{\Pi}\tilde{\vm{v}}_h$. Substituting into~\eqref{eq:L2_matrix_eq} and simplifying, we get the expression:
\begin{align*}
            (\tilde{\bm{\varepsilon}}^p)^T\left(\int_{E}{\left(\bm{N}^p\right)^T\bm{N}^p d\bm{x} }\right)\vm{\Pi}\tilde{\bm{v}}_h = (\tilde{\bm{\varepsilon}}^p)^T\left(\int_{\bm{\partial} E}{\left(\ac{\Nmatrix}\bm{N}^p\right)^T\bm{N}^v ds } - \int_{E}{\left(\bm{\partial} \bm{N}^p\right)^T\bm{N}^v d\bm{x}}\right)\tilde{\bm{v}}_h. 
\end{align*}
Define the matrix 
\begin{subequations}
\begin{align}
    \bm{G} &:= \int_{E} {\left(\bm{N}^p\right)^T\bm{N}^p \, d\bm{x} }. \label{G_matrix}\\
\intertext{\ac{Similar to~\eqref{eq:Btilde_matrix}, the last integral in~\eqref{eq:L2_matrix_eq} is not computable, so we again use the property in~\eqref{enrichment_condition} to construct an equivalent computable matrix:}}
    \bm{B} &:=\int_{\partial E} 
    {\left(\ac{\Nmatrix}\bm{N}^p\right)^T\bm{N}^v  \, ds}  - \int_{E}{\left(\bm{\partial} \bm{N}^p\right)^T\projS{\bm{N}^v} \, d\bm{x}}\label{B_matrix}.
\intertext{We now solve for the strain projection matrix $\vm{\Pi}$ in terms of $\vm{G}$ and $\vm{B}$ by}
    &\qquad \qquad \qquad \vm{\Pi} = \vm{G}^{-1}\vm{B}\label{def:Pi_matrix}.
\end{align}
\end{subequations}

\subsection{Implementation of element stiffness}
To construct the element stiffness, we first rewrite~\eqref{local_bilinearform} in terms of the matrices that have been constructed:
\begin{align*}
    a_h^E(\bm{u}_h,\bm{v}_h) &:= \int_{E}{\left(\overline{\projLtwo{\bm{v}_h}}\right)^T \bm{C}\, \overline{\projLtwo{\bm{u}_h}} \, d\bm{x}} \\ 
    &= \int_{E}{\left(\bm{N}^p\vm{\Pi}\tilde{\bm{v}}_h\right)^T\bm{C} \left(\bm{N}^p\vm{\Pi}\tilde{\bm{u}}_h\right) \, d\bm{x}} \\ 
    &= \tilde{\bm{v}}_h^T \left(\vm{\Pi}\right)^T
    \left(\int_{E}{\left(\bm{N}^p\right)^T}\bm{C}\bm{N}^p \, d\bm{x}\right)\vm{\Pi} \tilde{\bm{u}}_h .
\end{align*}
Then, define the element stiffness matrix $\bm{K}_E$ by
\begin{align}\label{eq:loc_stiff}
    \bm{K}_E :=\left(\vm{\Pi}\right)^T\left(\int_{E}{\left(\bm{N}^p\right)^T}\bm{C}\bm{N}^p \, d\bm{x}\right)\vm{\Pi} ,
\end{align}
where $\vm{\Pi}$ is given in~\eqref{def:Pi_matrix}.

\subsection{Implementation of element force vector}
We now construct the element forcing term given in~\eqref{local_force}
as
\begin{align*}
    b_h^E(\vm{v}_h) = \int_E{\vm{v}_h^T\vm{f_h} \,d\vm{x}} + \int_{\Gamma_N \cap \partial E}{\vm{v}_h^T\,\ac{\vm{t}_0}\,ds}.
\end{align*}
After making the approximations $\vm{f}_h = \projLtwodisp{\bm{f}}$ and simplifying, we can rewrite the expression in the form
\begin{subequations}
\begin{align}
    b_h^E(\vm{v}_h) &= \int_{E}{\vm{v}_h^T \projLtwodisp{\vm{f}}\,d\vm{x}} + \int_{\Gamma_N \cap \partial E}{\vm{v}_h^T\,{\ac{\vm{t}_0}}\,ds} \nonumber \\&=\int_{E}{\left(\projLtwodisp{\vm{v}_h}\right)^T \vm{f}\,d\vm{x}} + \int_{\Gamma_N \cap \partial E}{\left({\vm{N}^v\tilde{\vm{v}}_h}\right)^T\,\ac{\vm{t}_0}\,ds} \nonumber \\
    &=\tilde{\vm{v}}_h^T\left[\int_{E}{\left(\tilde{\vm{N}}^p\tilde{\vm{\Pi}}^0\right)^T \vm{f}\,d\vm{x}} + \int_{\Gamma_N \cap \partial E}{\left(\vm{N}^v\right)^T\,\ac{\vm{t}_0}\,ds}\right].
\intertext{Then, define the local forcing vector by}
    \bm{b}_h^E &:= \left[\int_{E}{\left(\tilde{\vm{N}}^p\tilde{\vm{\Pi}}^0\right)^T \vm{f}\,d\vm{x}} + \int_{\Gamma_N \cap \partial E}{\left(\vm{N}^v\right)^T\,\ac{\vm{t}_0}\,ds}\right]\label{eq:loc_force}.
\end{align}
\end{subequations}
All integrals that are required to form the element stiffness matrix in~\eqref{eq:loc_stiff} and the element force vector in~\eqref{eq:loc_force} are computed with the scaled boundary cubature (SBC) scheme~\cite{chin:2021:cmame}.

\ac{In general, the computational cost of the stabilization-free VEM is higher than that of standard VEM. For example, on a quadrilateral element with $k=2$, to construct the strain projection matrix $\vm{\Pi}$ requires the inversion of 
a $ 18 \times 18 $ matrix, while in the standard VEM~\cite{Artoli:2017:cmech} the strain projection only 
involves the inversion of a  $9 \times 9$ matrix. This cost significantly increases for polygonal meshes with elements that have many edges. However, hexagons tend to be dominant 
in polygonal (Voronoi) discretizations of solid geometries, 
so the increase in 
costs will not be substantial on such meshes. 
For nonlinear continua or problems such as 
acoustics or other wave phenomena that require
higher accuracy, finding a suitable stabilization term is  more involved and is not settled, and therefore the benefits that accrue 
on using a quadratic or cubic
stabilization-free virtual element formulation 
might outweigh the increase in computational costs. }

\section{Choice of \texorpdfstring{$\ell$}{l}}\label{sec:choice_ell}
In the previous sections we have left the choice of $\ell := \ell(E)$ open. We now numerically establish a choice of $\ell$ that results in a well-posed, stable discrete problem. We examine specifically the case of second- and third-order methods. \ac{Discussion of the choice of $\ell$ for the first-order method can be found in~\cite{berrone2021lowest,chen:2022:arxiv}}.

\subsection{Eigenanalysis for regular polygons}
\begin{figure}[H]
     \centering
     \begin{subfigure}{0.32\textwidth}
         \centering
         \includegraphics[width=\textwidth]{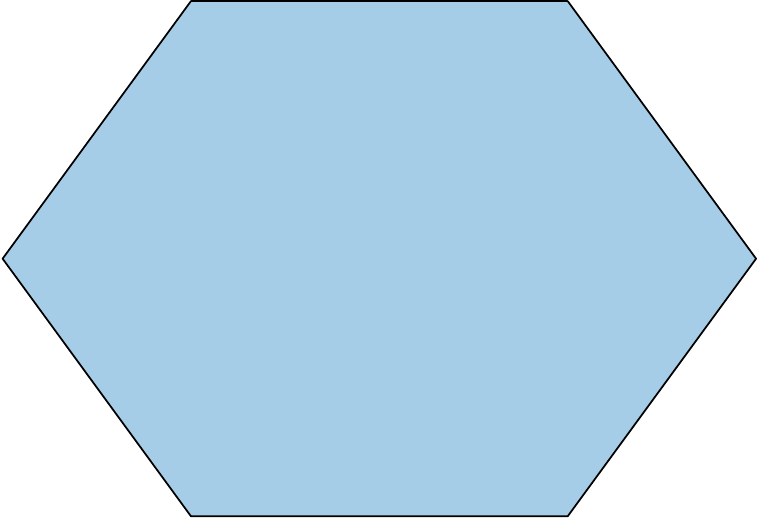}
         \caption{}
         \label{fig:hex}
     \end{subfigure}
     \hfill
     \begin{subfigure}{0.32\textwidth}
         \centering
         \includegraphics[width=\textwidth]{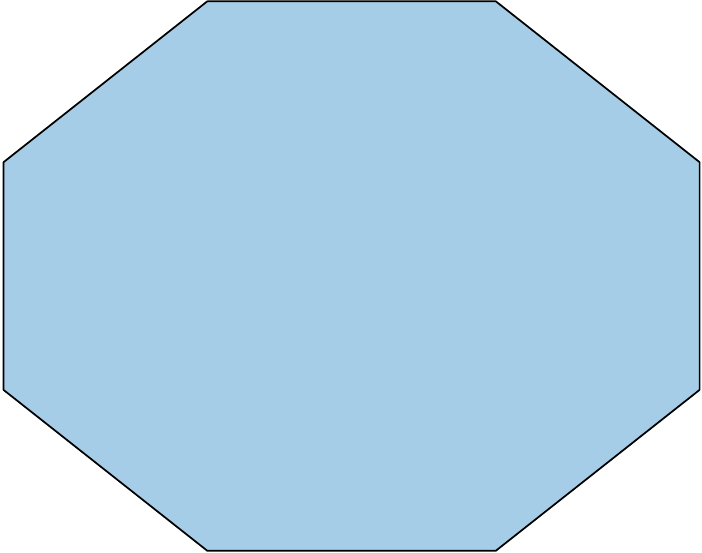}
         \caption{}
         \label{fig:octa}
     \end{subfigure}
     \hfill
     \begin{subfigure}{0.32\textwidth}
         \centering
         \includegraphics[width=\textwidth]{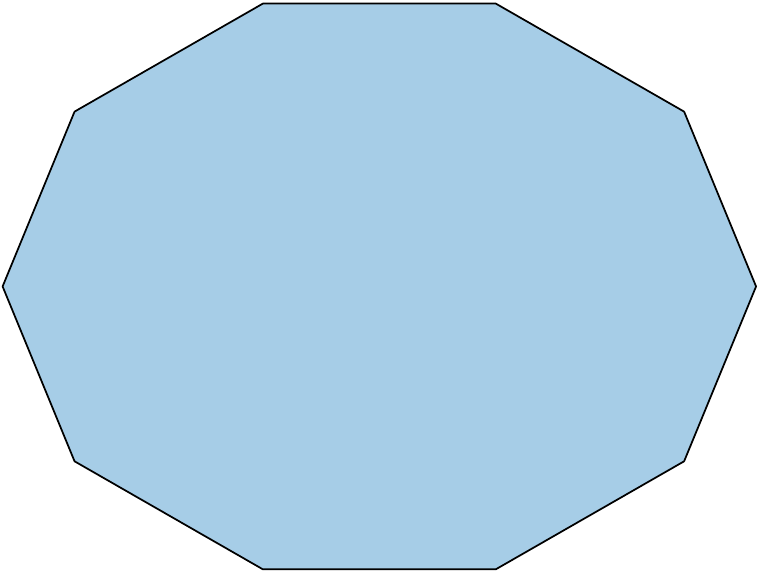}
         \caption{}
         \label{fig:deca}
     \end{subfigure}
        \caption{Regular polygons that are used in the eigenanalysis.}
        \label{fig:three_polygons}
\end{figure}

\begin{figure}[H]
     \centering
     \begin{subfigure}{0.32\textwidth}
         \centering
         \includegraphics[width=\textwidth]{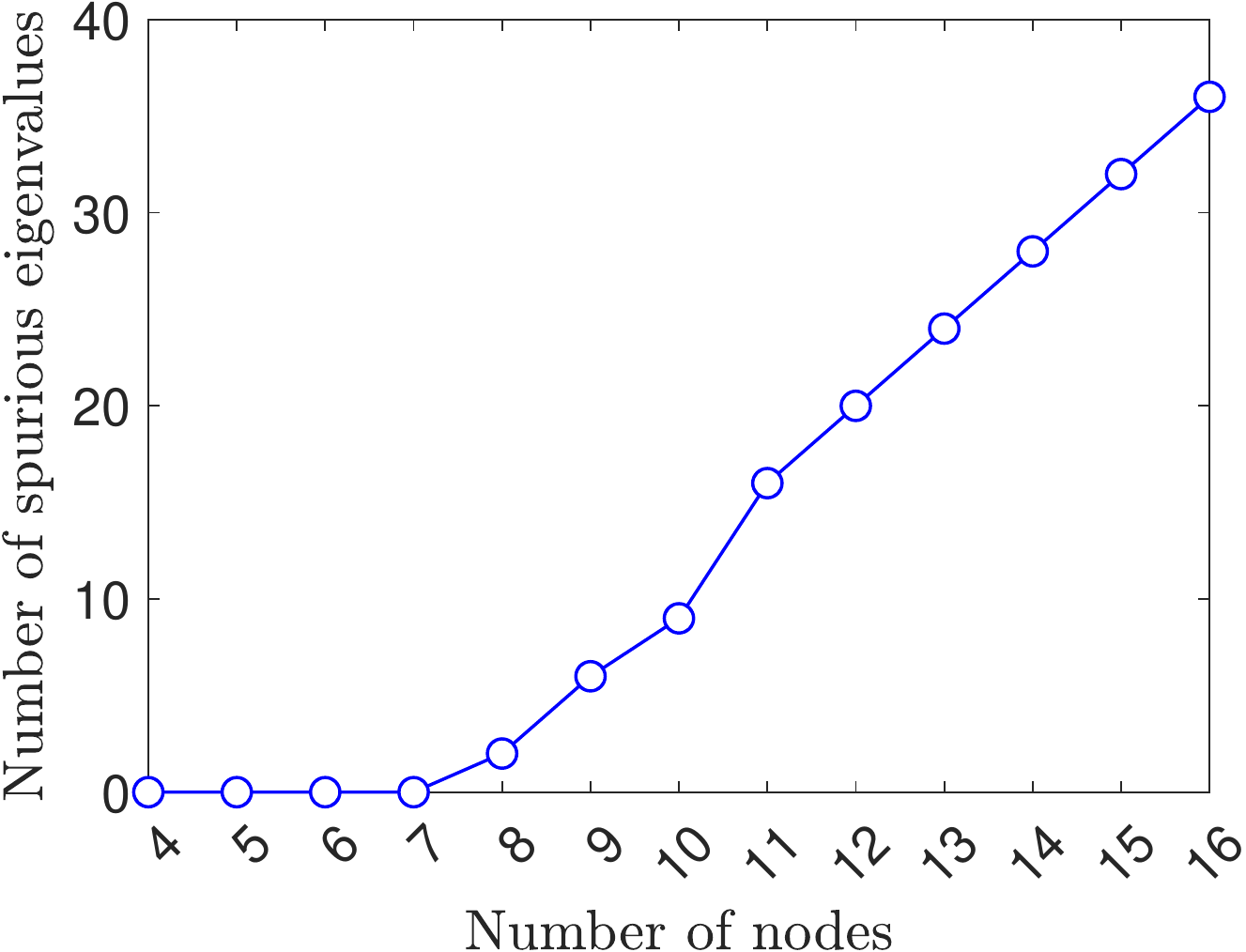}
         \caption{}
         \label{}
     \end{subfigure}
     \hfill
     \begin{subfigure}{0.32\textwidth}
         \centering
         \includegraphics[width=\textwidth]{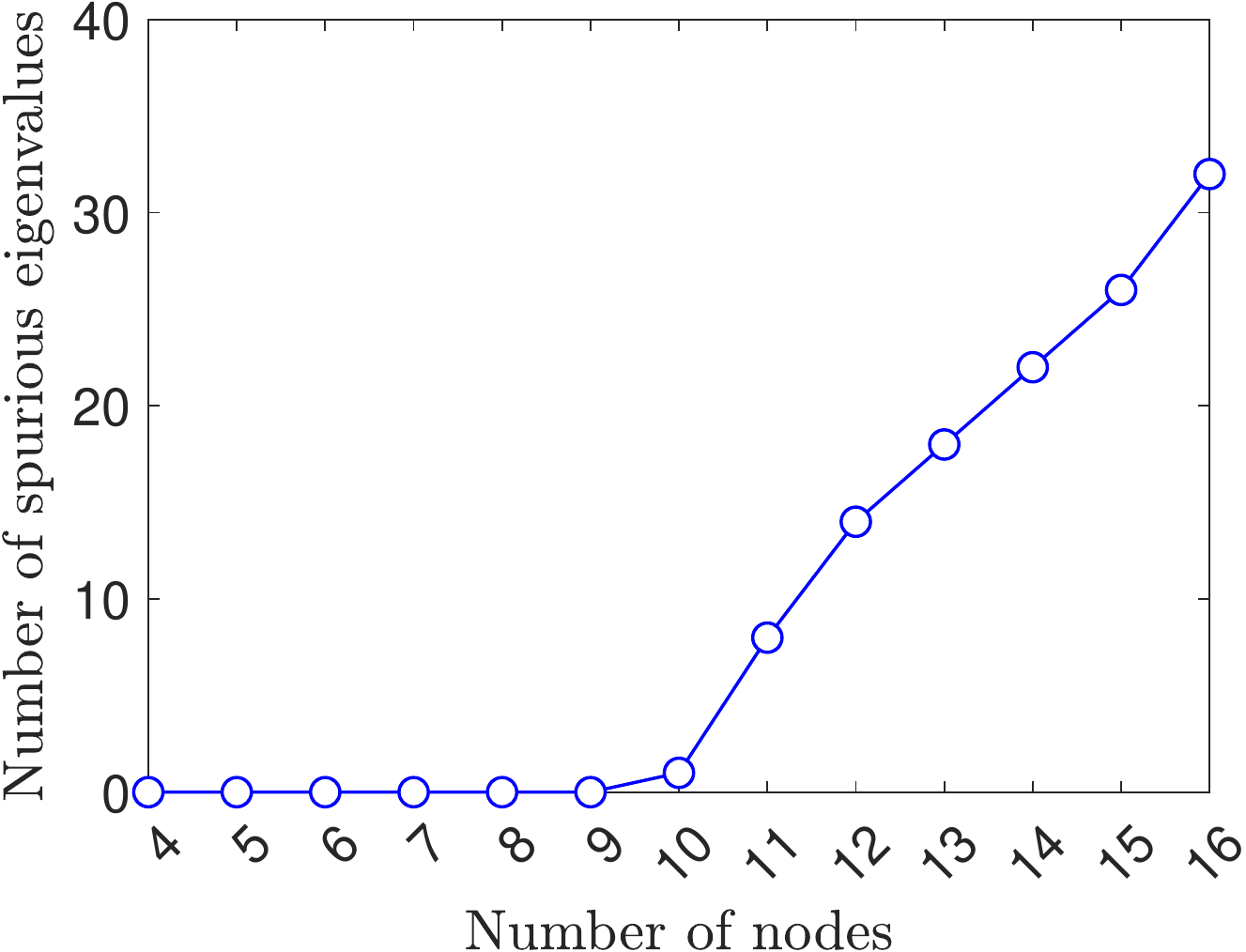}
         \caption{}
         \label{}
     \end{subfigure}
     \hfill
     \begin{subfigure}{0.32\textwidth}
         \centering
         \includegraphics[width=\textwidth]{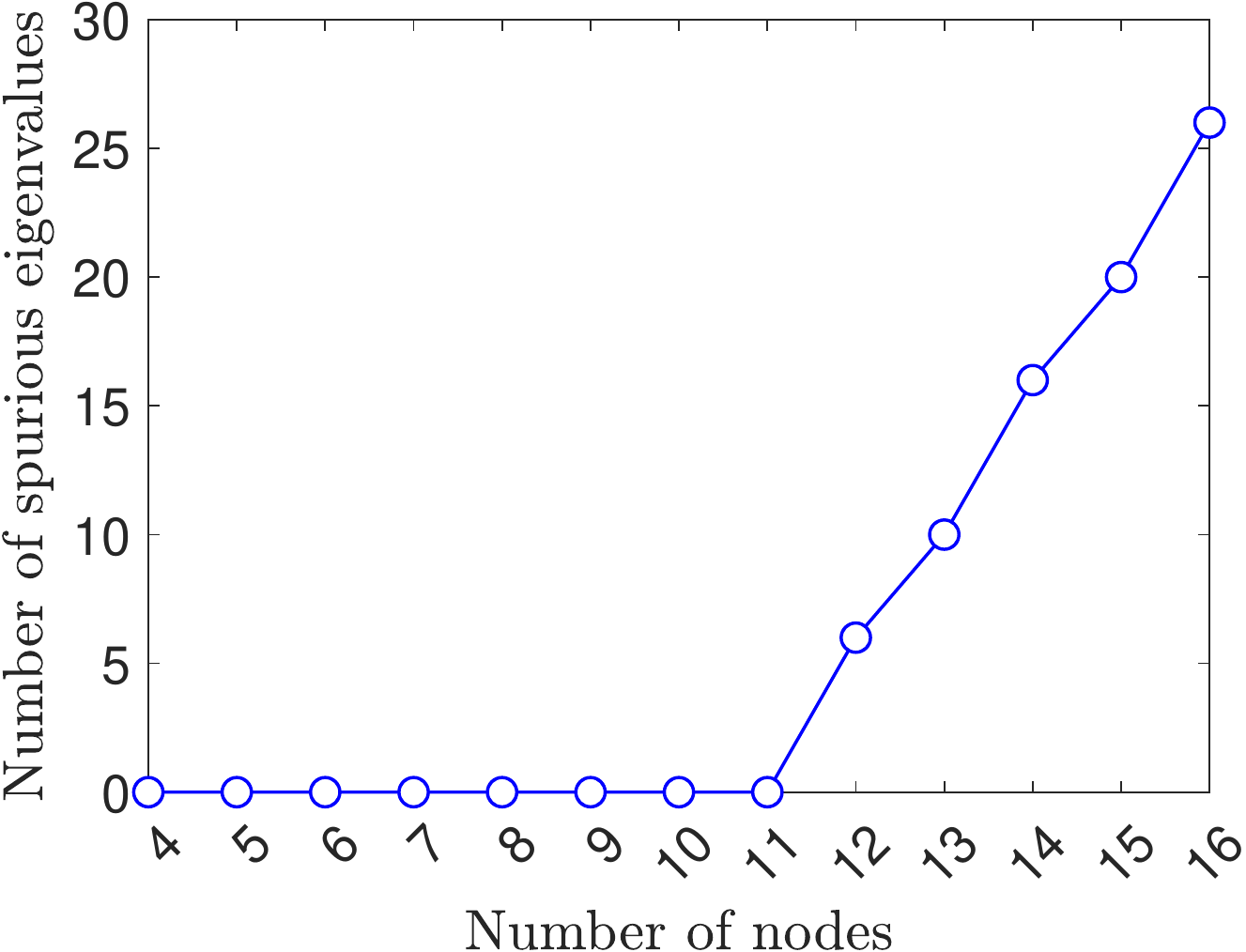}
         \caption{}
         \label{}
     \end{subfigure}
        \caption{Eigenvalue analysis on regular polygons with the second-order
        method. (a) $\ell =3$, (b) $\ell=4$ and (c) $\ell=5$. }
        \label{fig:eig_plots_k2_poly}
\end{figure}

\begin{figure}[H]
     \centering
     \begin{subfigure}{0.32\textwidth}
         \centering
         \includegraphics[width=\textwidth]{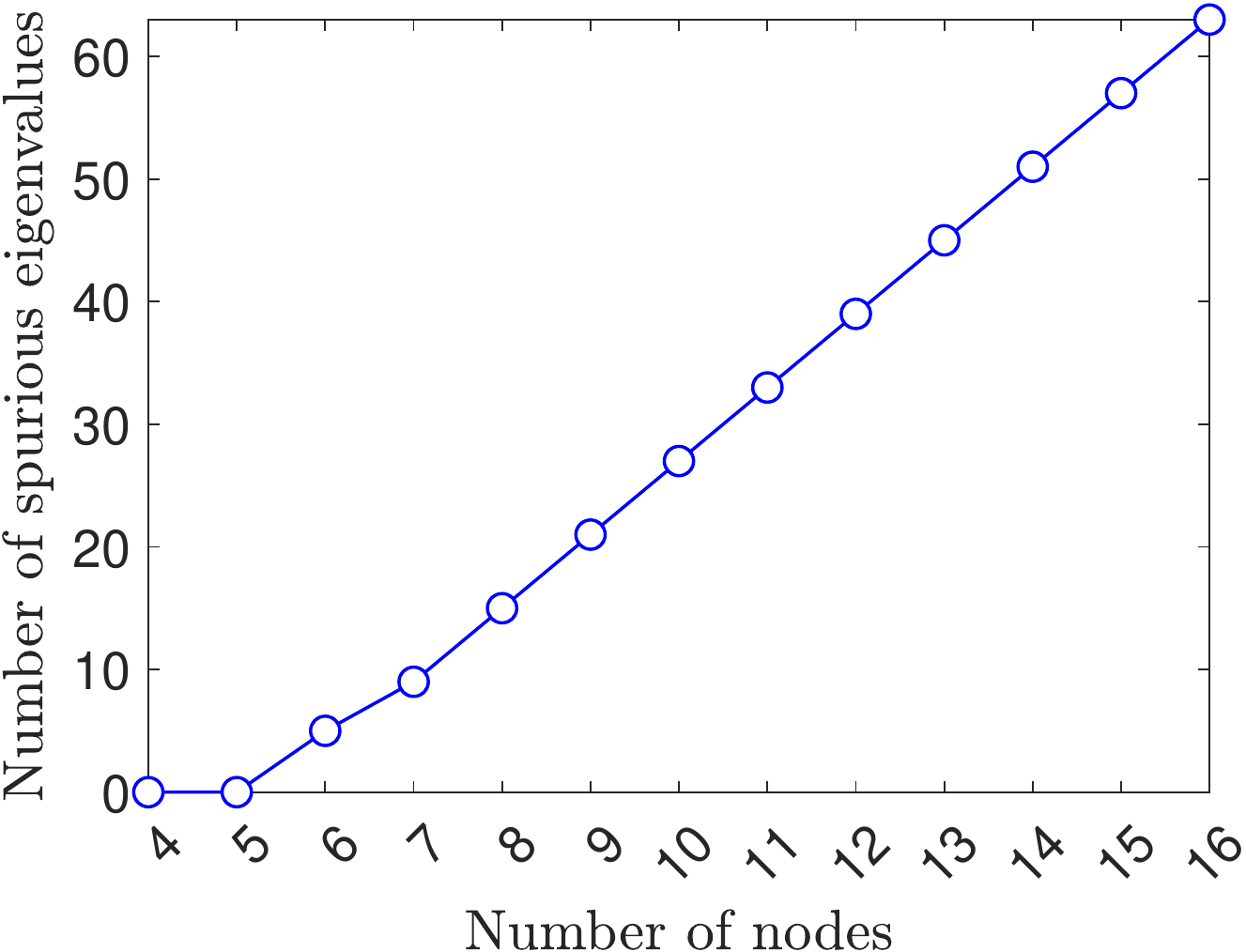}
         \caption{}
         \label{}
     \end{subfigure}
     \hfill
     \begin{subfigure}{0.32\textwidth}
         \centering
         \includegraphics[width=\textwidth]{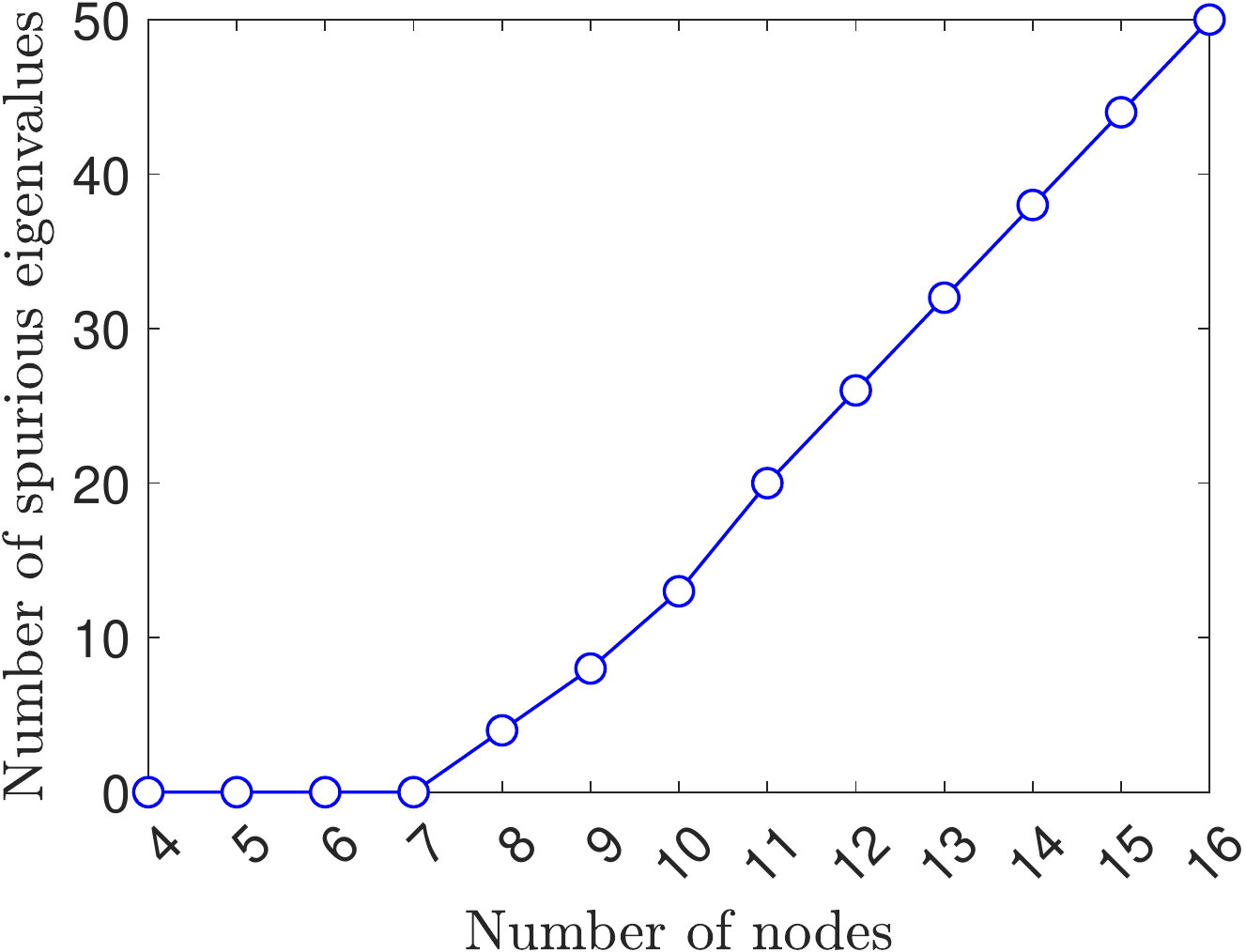}
         \caption{}
         \label{}
     \end{subfigure}
     \hfill
     \begin{subfigure}{0.32\textwidth}
         \centering
         \includegraphics[width=\textwidth]{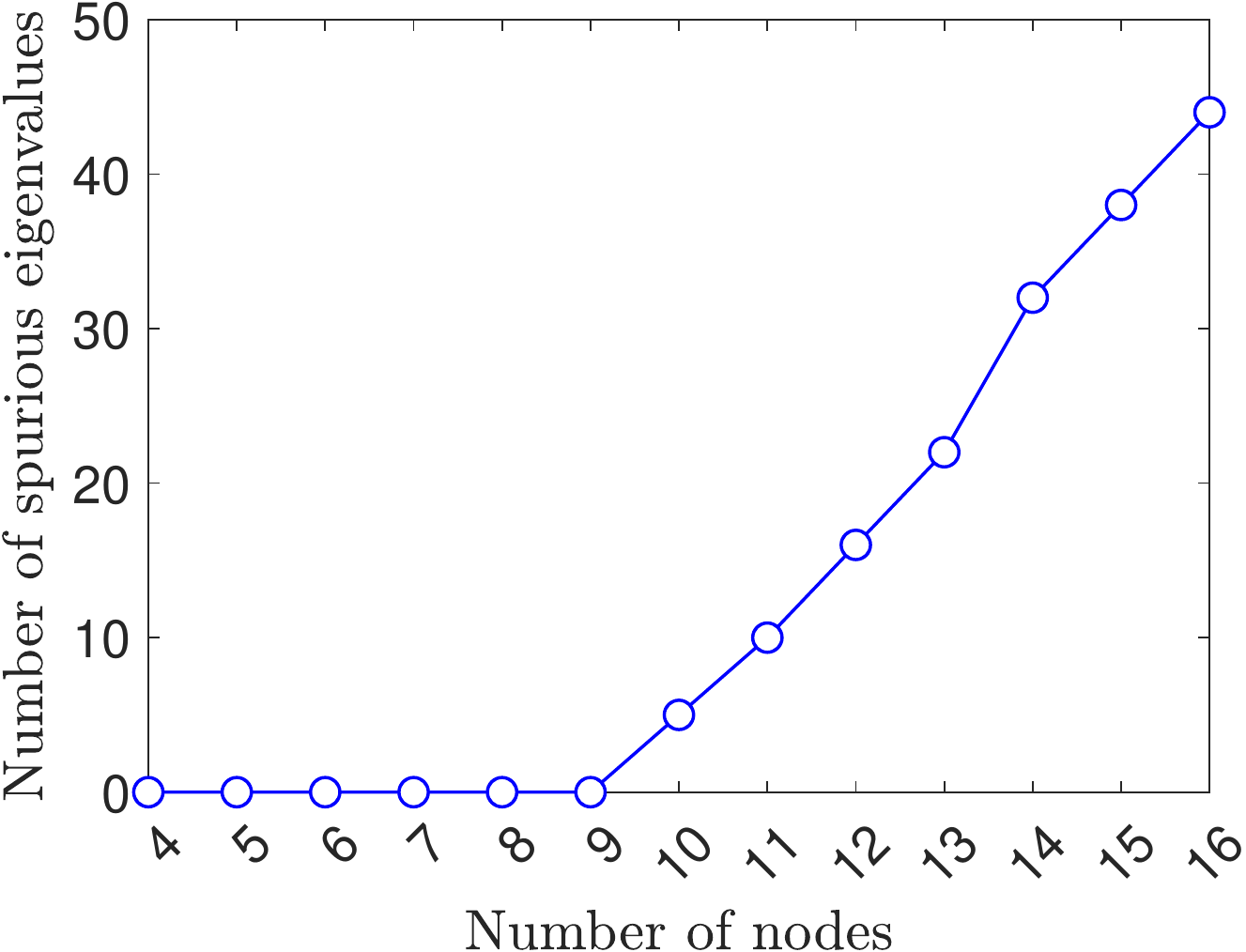}
         \caption{}
         \label{}
     \end{subfigure}
        \caption{Eigenvalue analysis on regular polygons with the third-order
        method. (a) $\ell =3$, (b) $\ell=4$ and (c) $\ell=5$. }
        \label{fig:eig_plots_k3_poly}
\end{figure}

 We first study the stability on regular polygons by considering the element eigenvalue problem $\bm{K}_E \bm{d}_E = \lambda \bm{d}_E$. For plane elasticity, the element stiffness should have three zero eigenvalues that correspond to the three rigid-body modes, with any additional zero eigenvalue being a non-physical (spurious) mode.
We measure the number of spurious eigenvalues of the local stiffness matrix over the set of regular $n$-gons.  We fix $\ell=3,4,5$ and measure 
the number of spurious eigenvalues on a given regular polygon. In Figure~\ref{fig:three_polygons} a few sample polygons are shown, and in Figures~\ref{fig:eig_plots_k2_poly} and~\ref{fig:eig_plots_k3_poly} we plot the number of spurious eigenvalues as a function of the vertices of the corresponding polygon. The analyses over regular polygons reveals that the element stiffness matrix is stable with the correct rank if the inequalities $N_E \leq 2\ell+1$ and $N_E \leq 2\ell-1$ hold for $k=2,3$ respectively. In~\cite{chen:2022:arxiv}, it was shown that for $k=1$ the inequality is given by $N_E \leq 2\ell +3$ for regular polygons. We conjecture that this pattern holds, and for a general $k$-th order method, a sufficient inequality is given by  $N_E \leq 2\ell -2k+5$.

We are also interested in robustness of the inequality when the vertices of an element are perturbed. In particular, for $k=2,3$ we first fix $\ell =3,4,5$, then take the respective regular hexagon, octagon, decagon and perturb one component of a vertex by $\delta$. We measure the number of spurious modes as a function of $\delta$. For $k=2$ the three elements will satisfy the inequality $N_E \leq 2\ell+1$, so we expect no spurious eigenvalues to appear, but for $k=3$ the inequality $N_E\leq 2\ell-1$ is not satisfied so we expect to see some additional spurious eigenvalues.  
\begin{figure}[H]
     \centering
     \begin{subfigure}{0.32\textwidth}
         \centering
         \includegraphics[width=\textwidth]{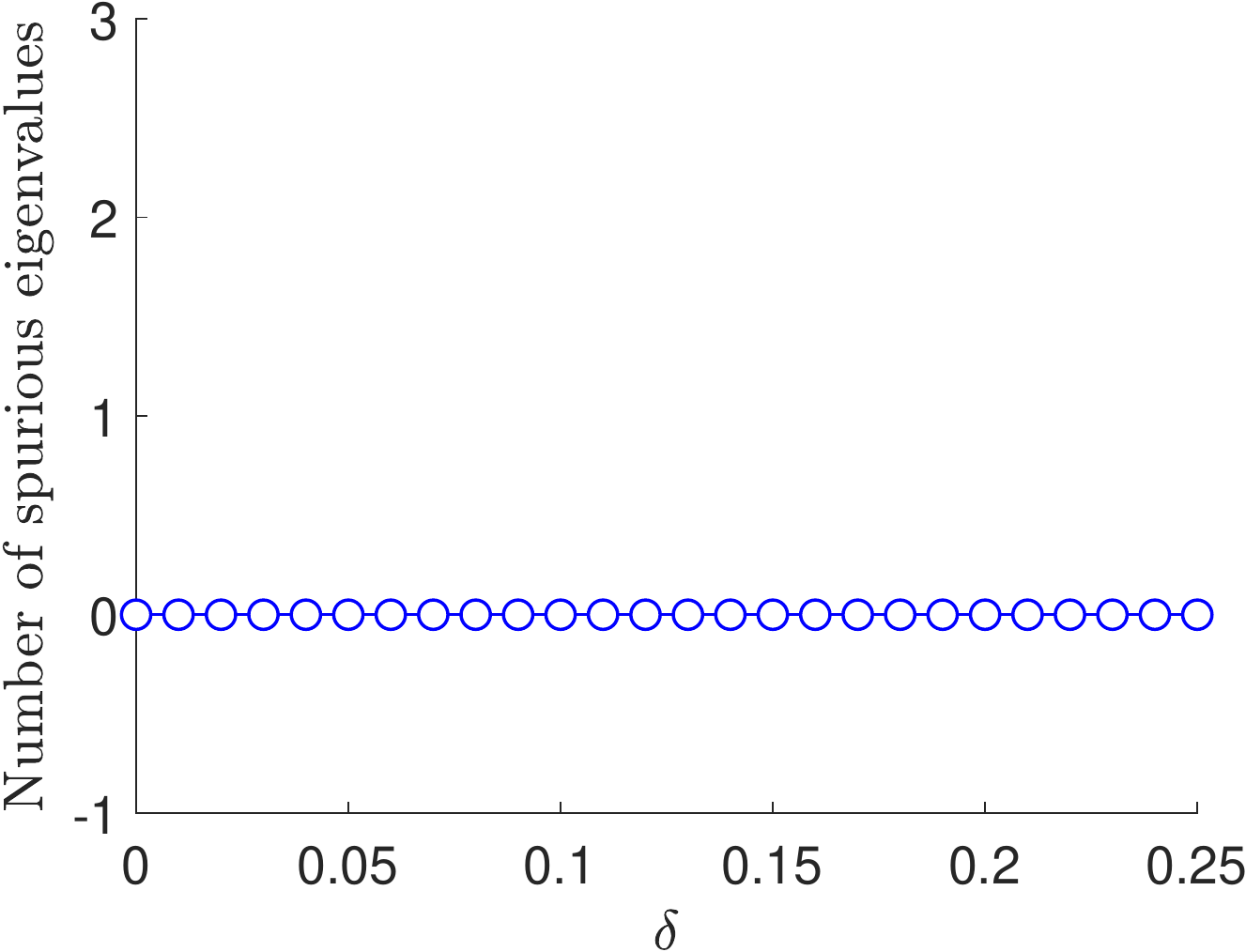}
         \caption{}
         \label{}
     \end{subfigure}
     \hfill
     \begin{subfigure}{0.32\textwidth}
         \centering
         \includegraphics[width=\textwidth]{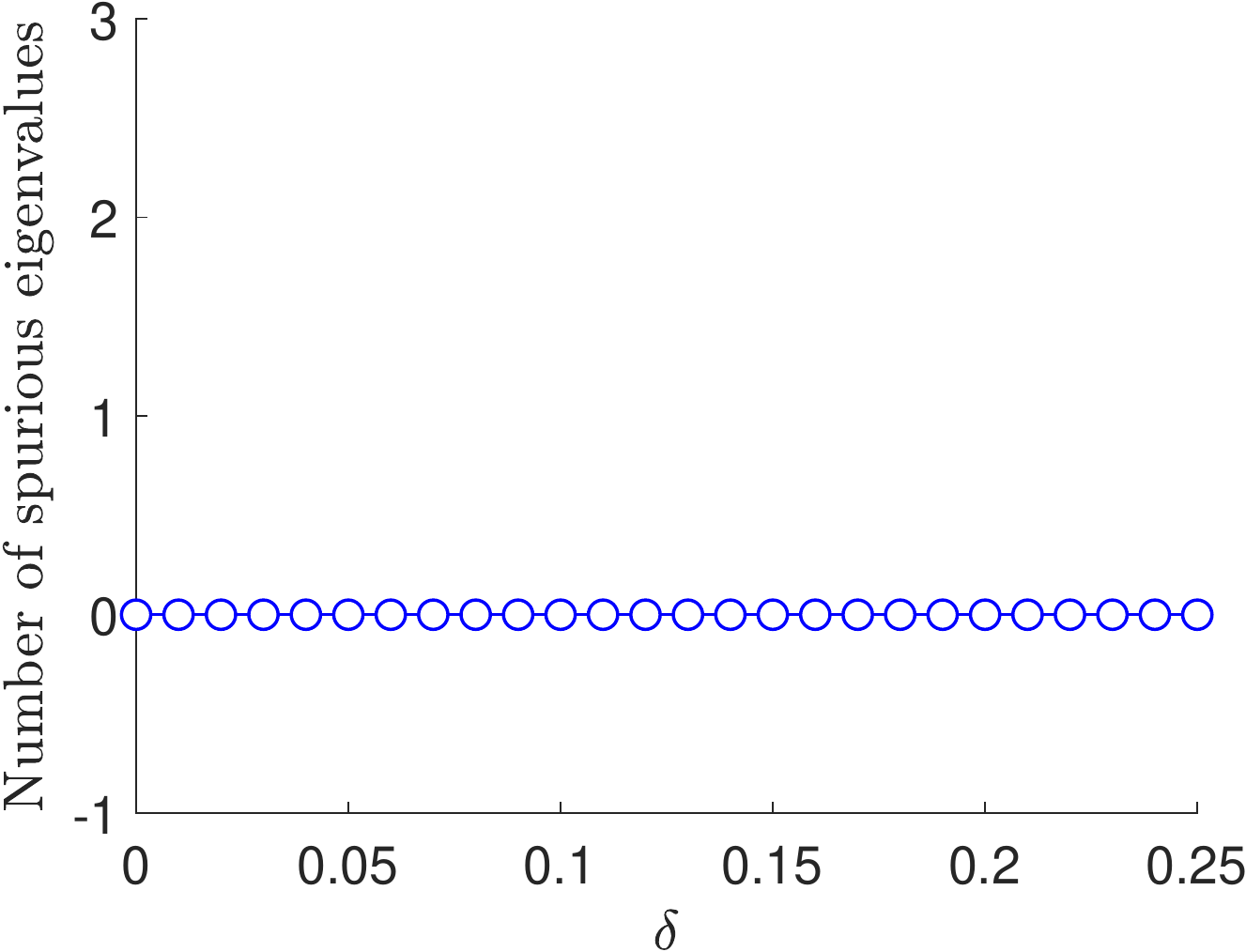}
         \caption{}
         \label{}
     \end{subfigure}
     \hfill
     \begin{subfigure}{0.32\textwidth}
         \centering
         \includegraphics[width=\textwidth]{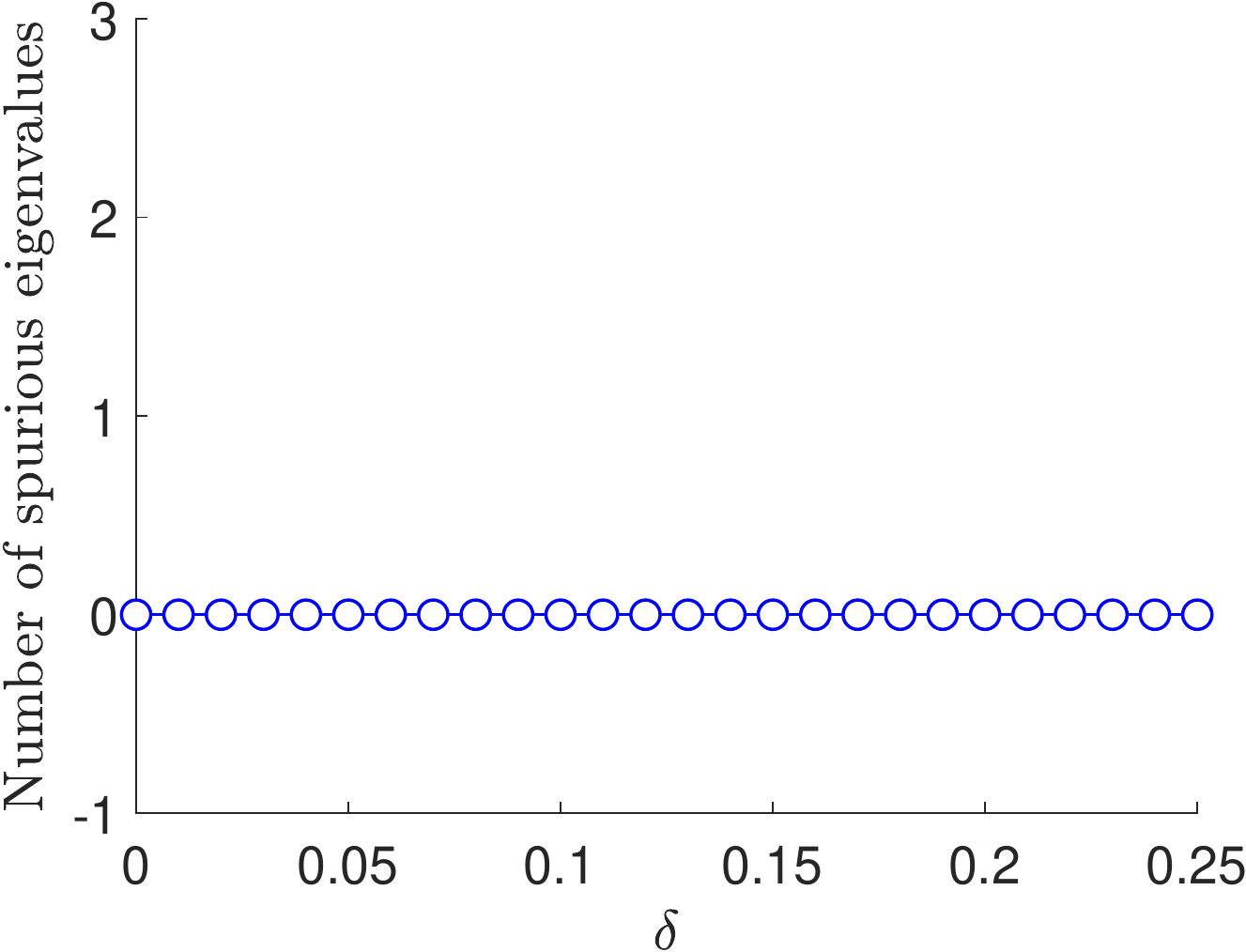}
         \caption{}
         \label{}
     \end{subfigure}
        \caption{Eigenvalue analysis on the perturbed regular polygons with the second-order method. 
        (a) $\ell =3$ on hexagon, (b) $\ell=4$ on octagon and (c) $\ell=5$ on decagon. }
        \label{fig:eig_perturbed_k2}
\end{figure}
\begin{figure}[H]
     \centering
     \begin{subfigure}{0.32\textwidth}
         \centering
         \includegraphics[width=\textwidth]{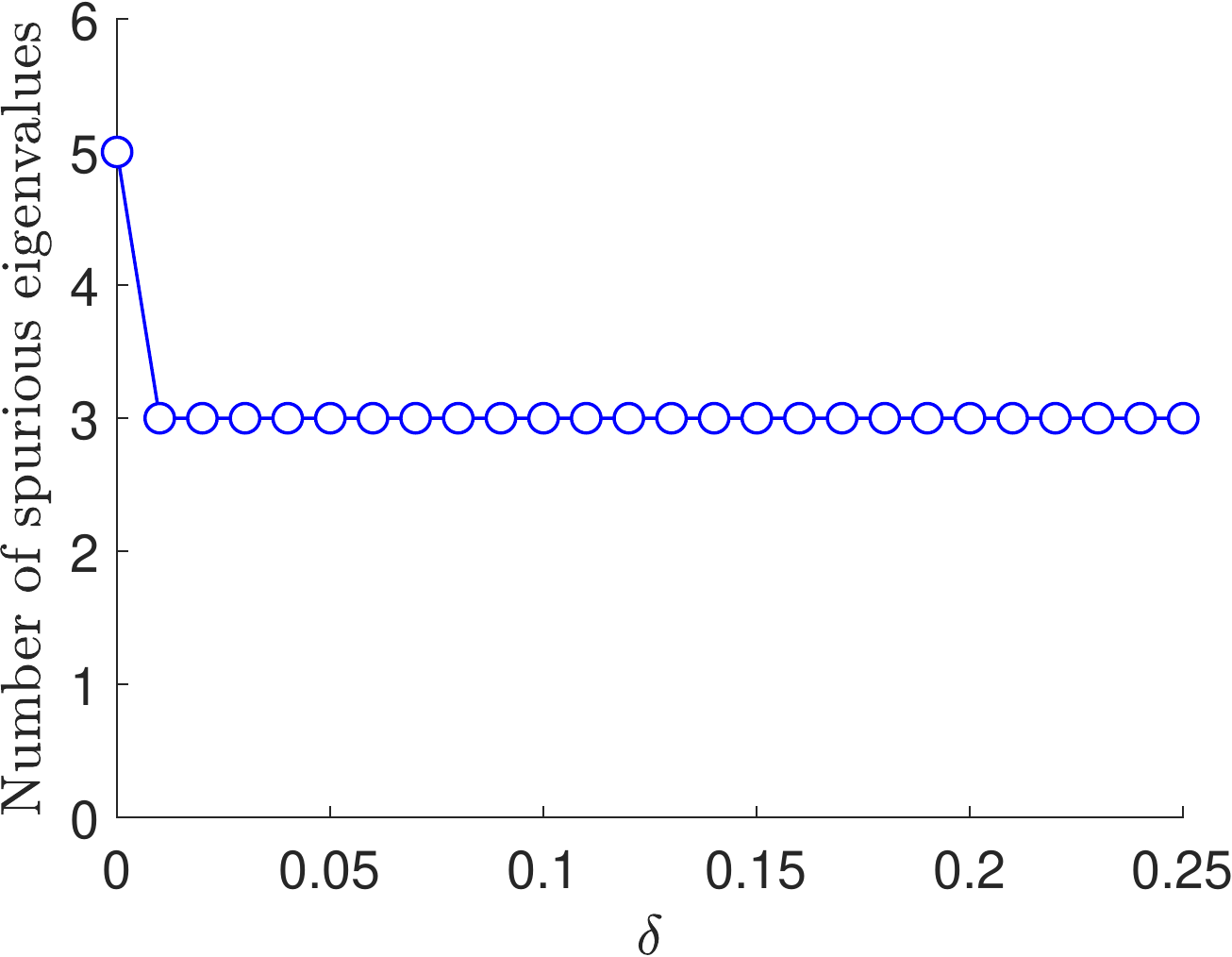}
         \caption{}
         \label{}
     \end{subfigure}
     \hfill
     \begin{subfigure}{0.32\textwidth}
         \centering
         \includegraphics[width=\textwidth]{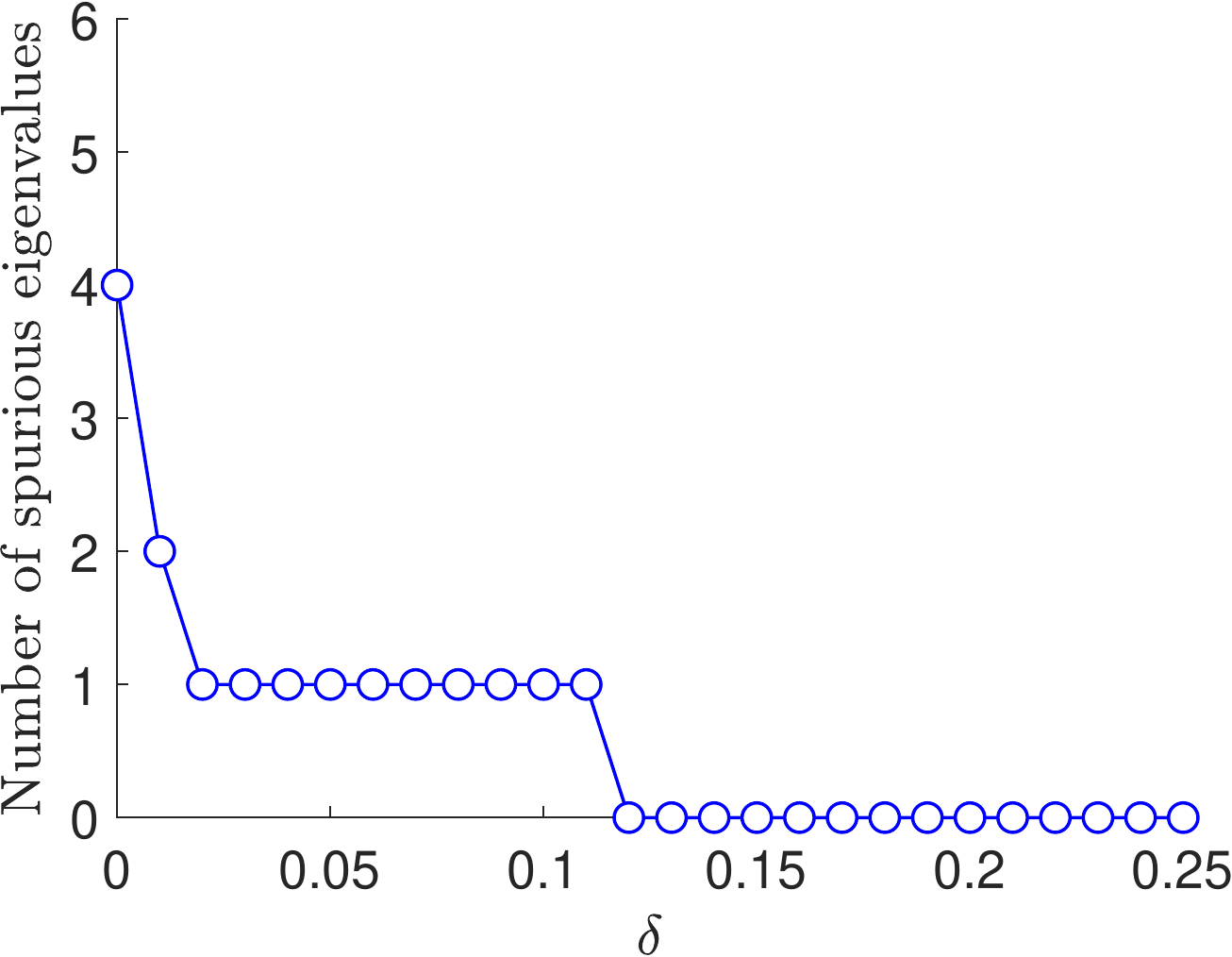}
         \caption{}
         \label{}
     \end{subfigure}
     \hfill
     \begin{subfigure}{0.32\textwidth}
         \centering
         \includegraphics[width=\textwidth]{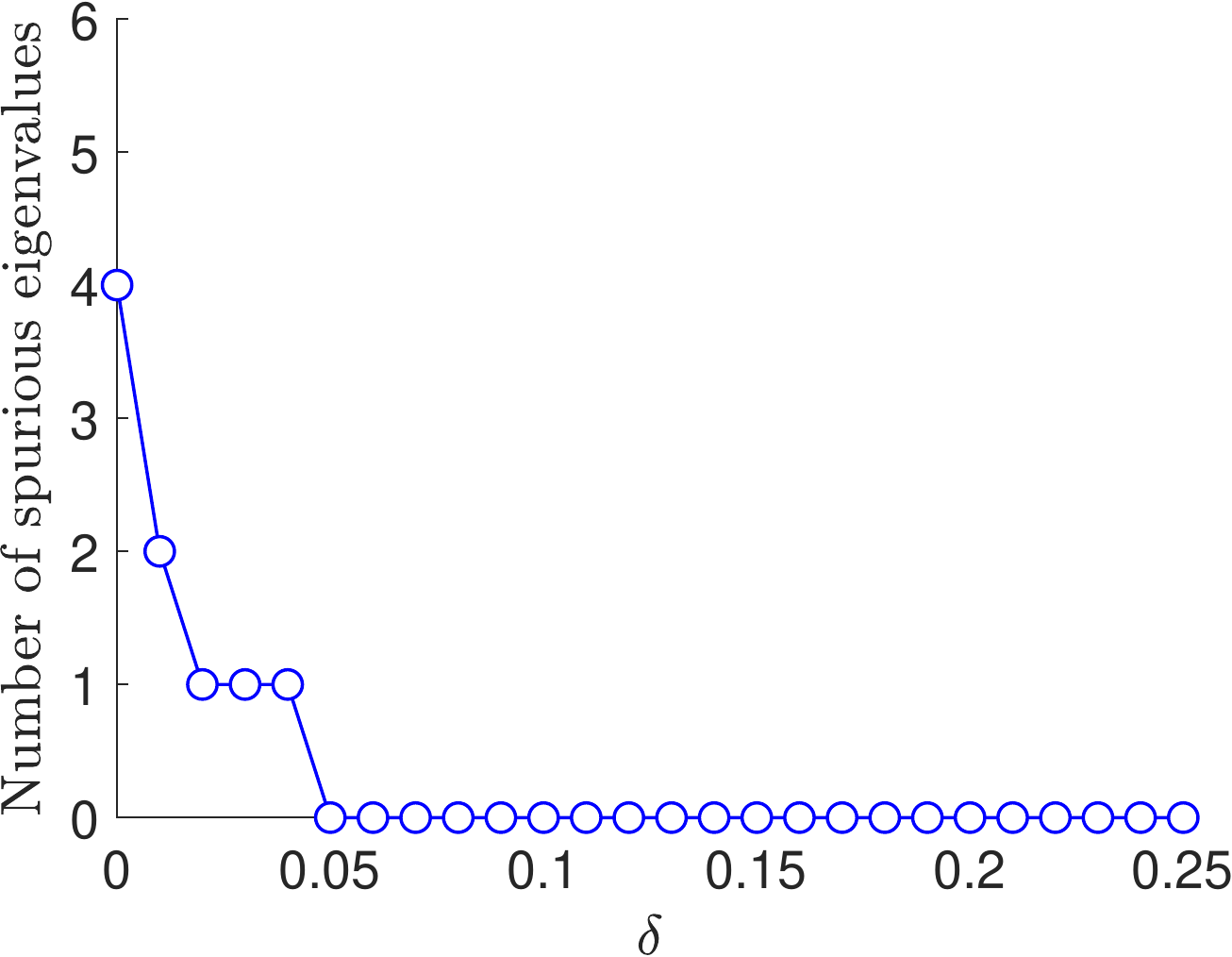}
         \caption{}
         \label{}
     \end{subfigure}
        \caption{Eigenvalue analysis on the perturbed regular polygons with the third-order method. 
        (a) $\ell =3$ on hexagon, (b) $\ell=4$ on octagon and (c) $\ell=5$ on decagon. }
        \label{fig:eig_perturbed_k3}
\end{figure}
From Figure~\ref{fig:eig_perturbed_k2}, we observe that for small perturbations of the hexagon, octagon, and decagon that no spurious eigenvalues arise. In Figure~\ref{fig:eig_perturbed_k3}, we see that by perturbing the octagon and decagon, we are able to reduce the number of spurious eigenvalues to zero when using $\ell=4,5$ respectively. 
 
\subsection{Eigenanalysis for general polygons}
\begin{figure}[H]
     \centering
     \begin{subfigure}{0.32\textwidth}
         \centering
         \includegraphics[width=\textwidth]{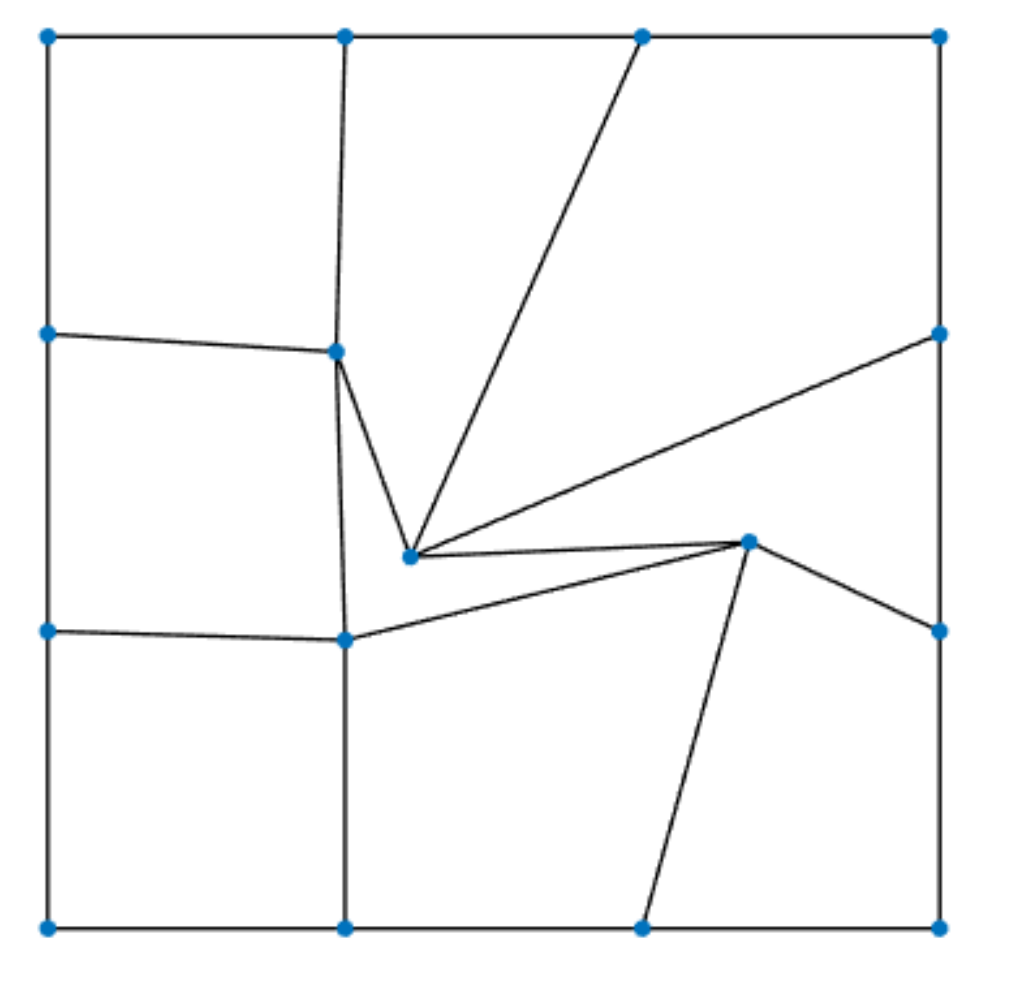}
         \caption{}
         \label{}
     \end{subfigure}
     \hfill
     \begin{subfigure}{0.32\textwidth}
         \centering
         \includegraphics[width=\textwidth]{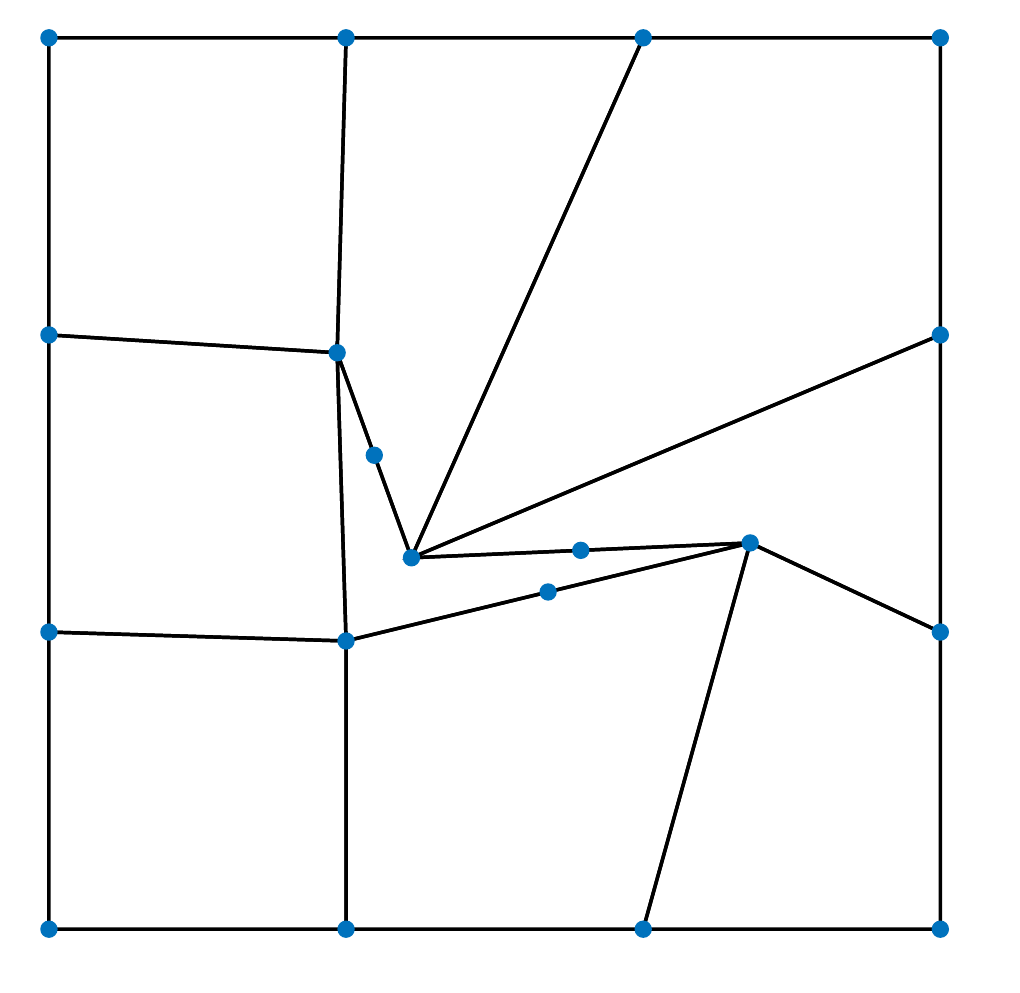}
         \caption{}
         \label{}
     \end{subfigure}
     \hfill
     \begin{subfigure}{0.32\textwidth}
         \centering
         \includegraphics[width=\textwidth]{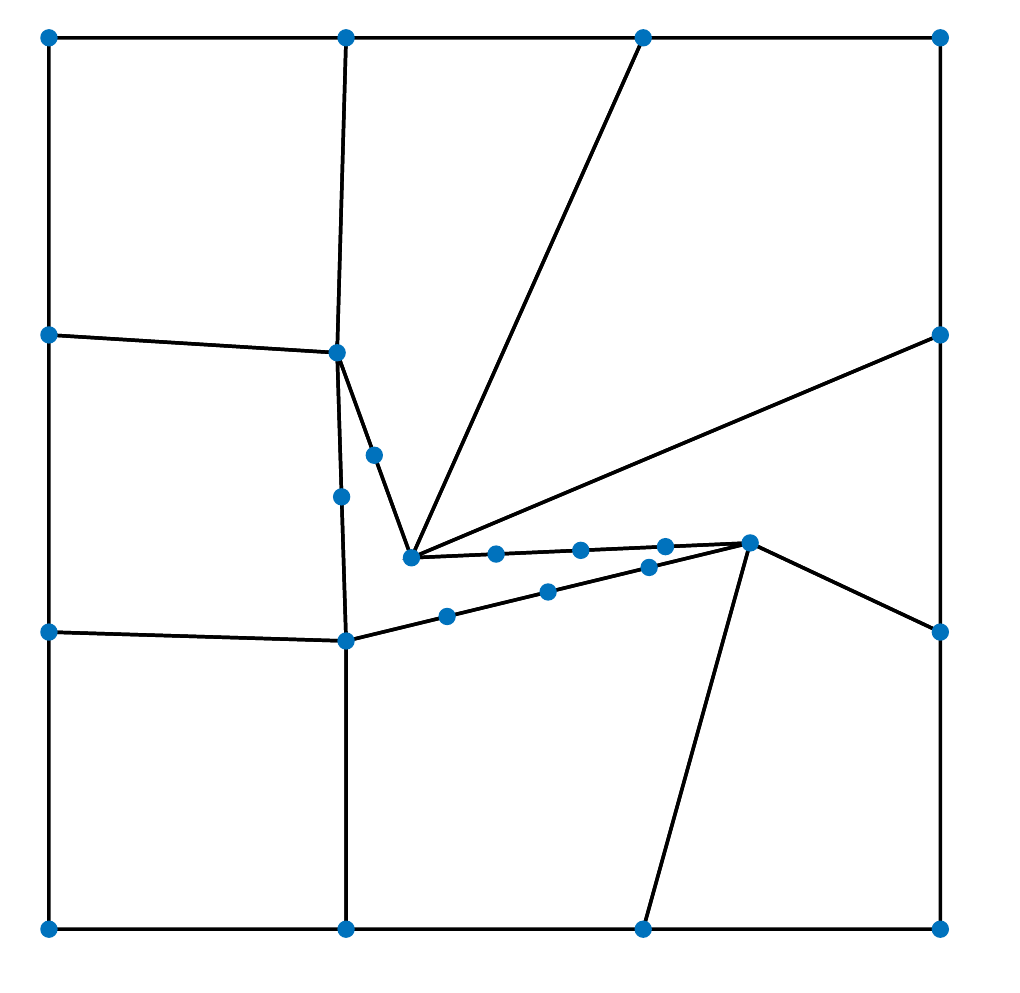}
         \caption{}
         \label{}
     \end{subfigure}
        \caption{Sample meshes used in the element-eigenvalue analysis
        for $\ell = 3,4,5$. The central quadrilateral element has (a) 4 nodes, (b) 7 nodes, and (c) 12 nodes. }
        \label{fig:eig_mesh}
\end{figure}

\begin{figure}[H]
     \centering
     \begin{subfigure}{0.32\textwidth}
         \centering
         \includegraphics[width=\textwidth]{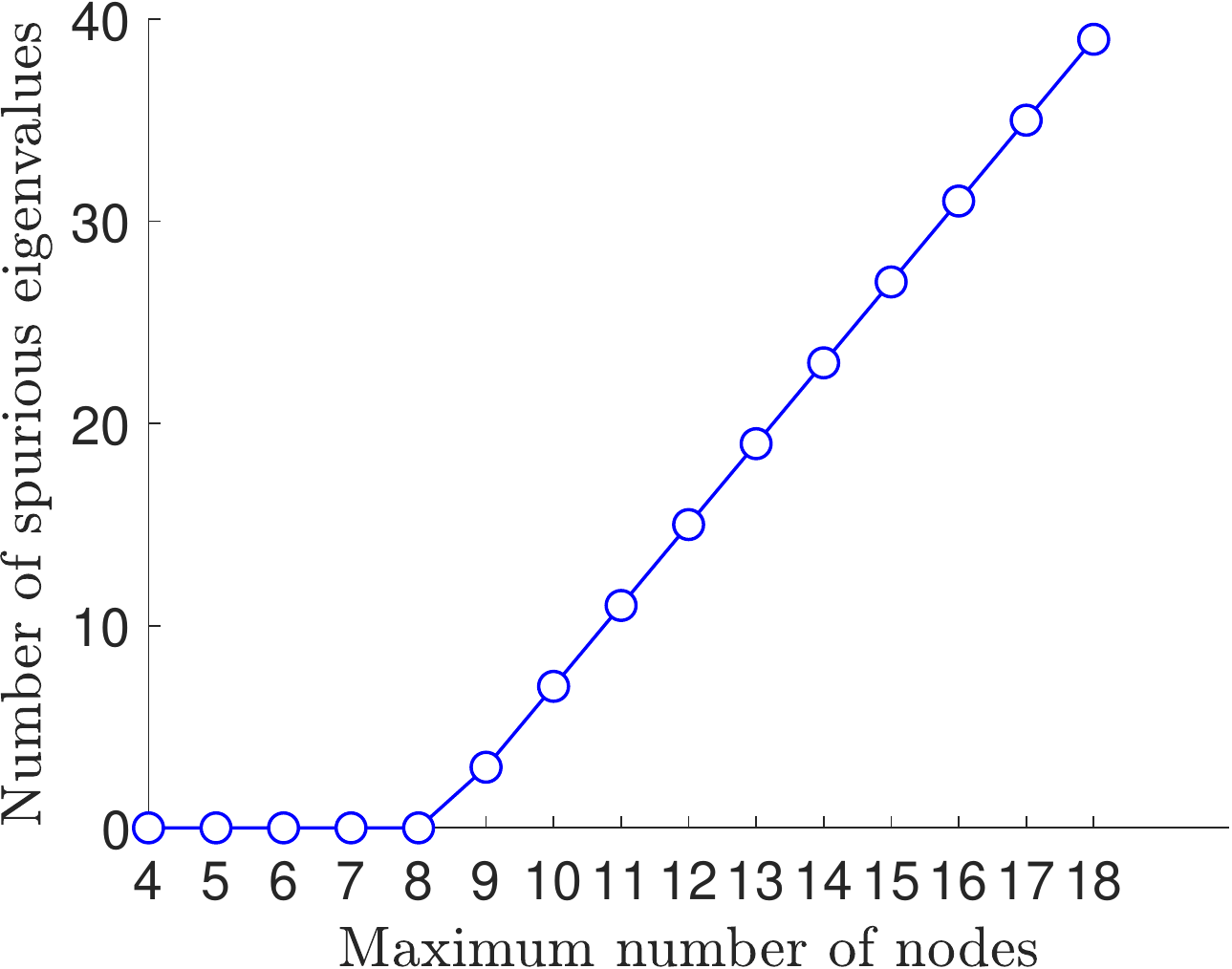}
         \caption{}
         \label{}
     \end{subfigure}
     \hfill
     \begin{subfigure}{0.32\textwidth}
         \centering
         \includegraphics[width=\textwidth]{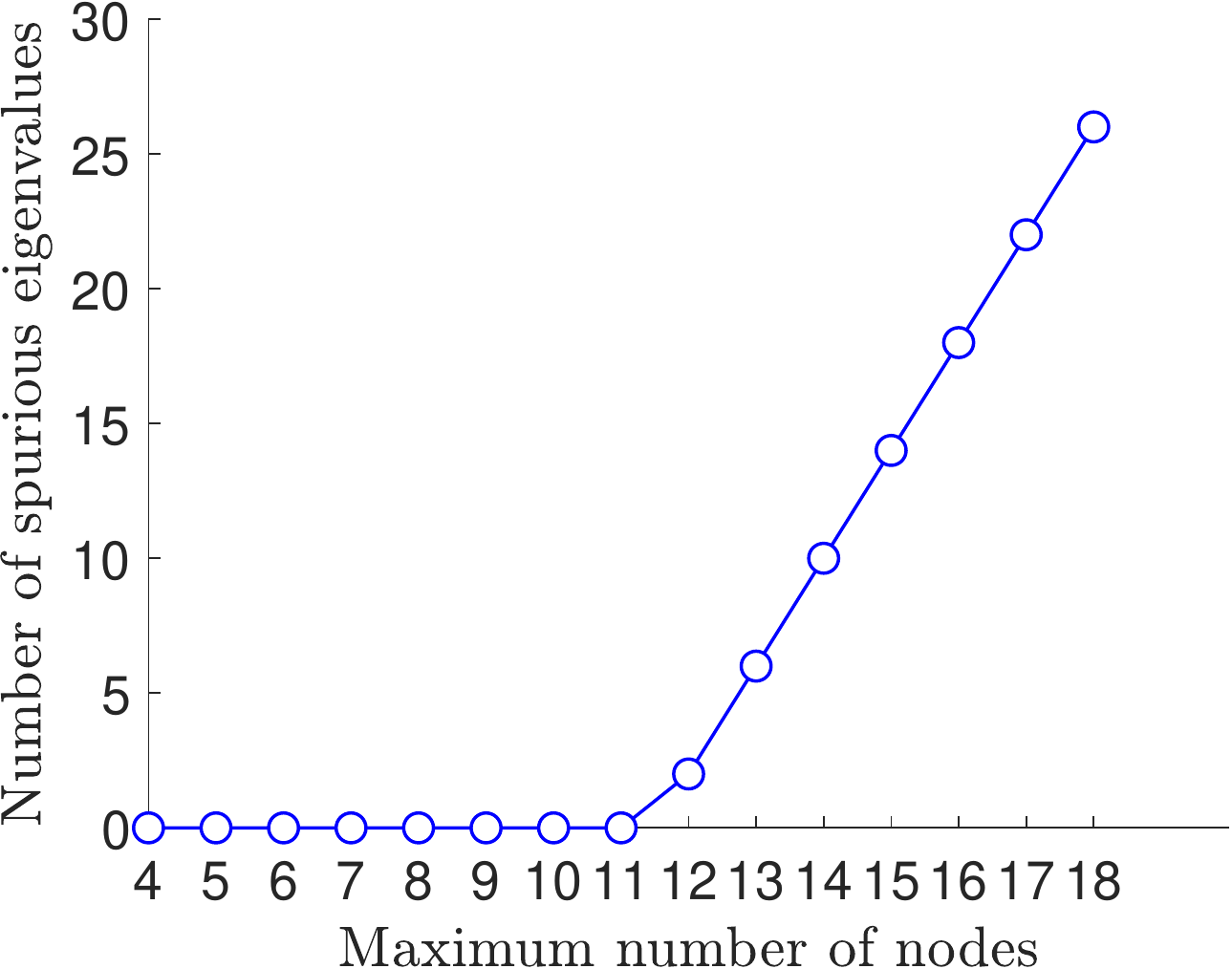}
         \caption{}
         \label{}
     \end{subfigure}
     \hfill
     \begin{subfigure}{0.32\textwidth}
         \centering
         \includegraphics[width=\textwidth]{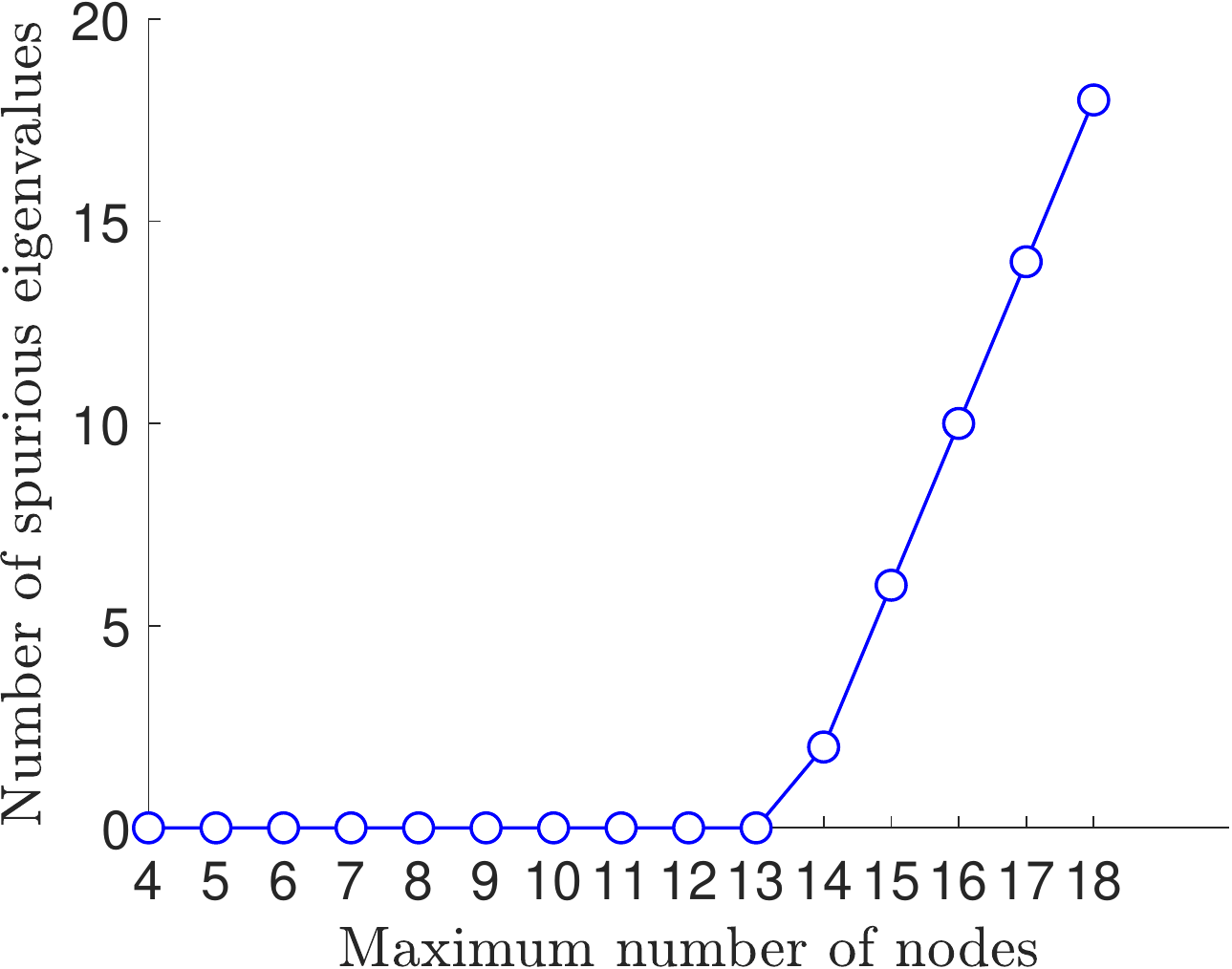}
         \caption{}
         \label{}
     \end{subfigure}
        \caption{Eigenvalue analysis on the meshes shown 
        in Figure~\protect\ref{fig:eig_mesh} with the second-order method. 
        (a) $\ell =3$, (b) $\ell=4$ and (c) $\ell=5$. }
        \label{fig:eig_plots_k2}
\end{figure}

\begin{figure}[H]
     \centering
     \begin{subfigure}{0.32\textwidth}
         \centering
         \includegraphics[width=\textwidth]{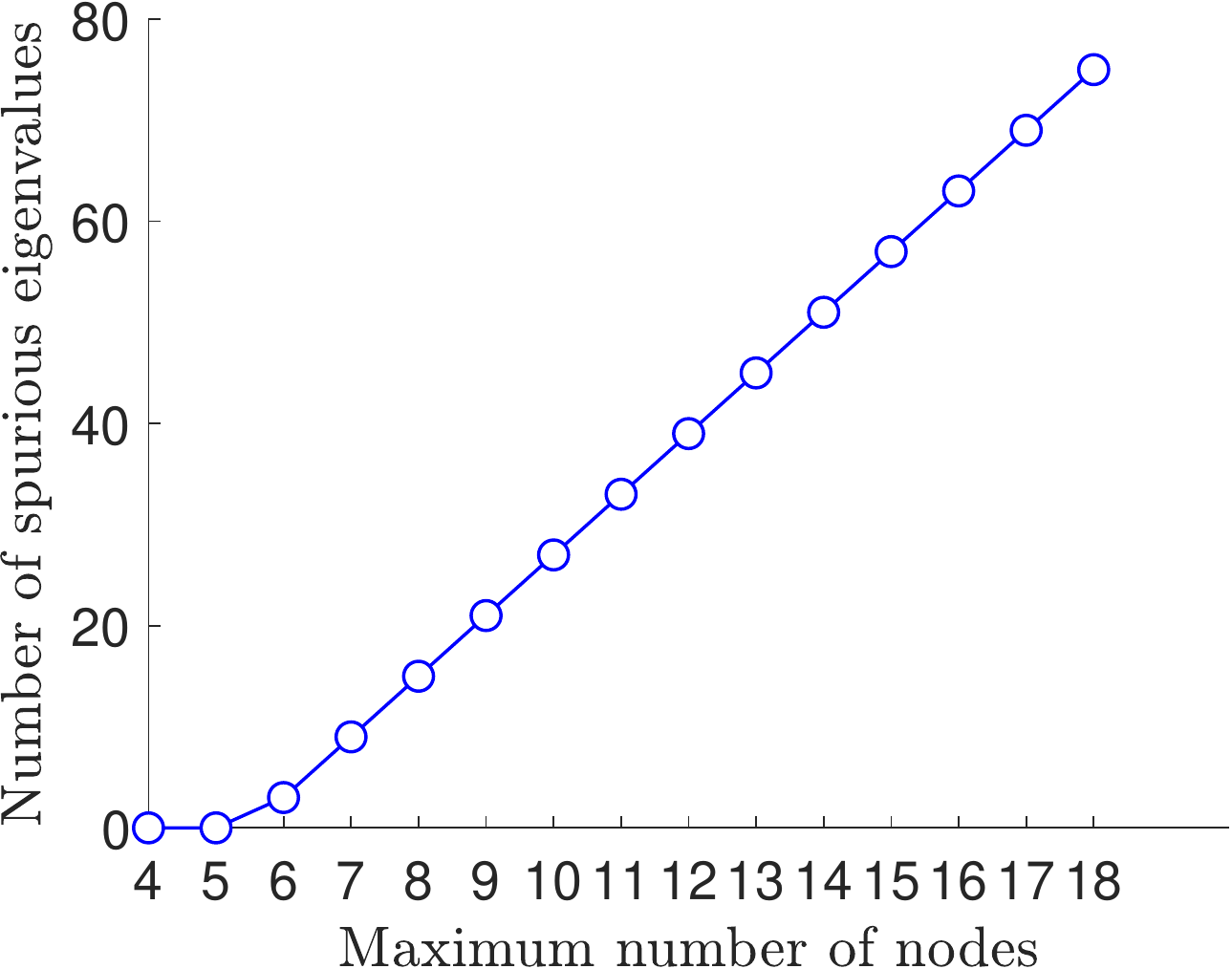}
                  \caption{}
         \label{}
     \end{subfigure}
     \hfill
     \begin{subfigure}{0.32\textwidth}
         \centering
         \includegraphics[width=\textwidth]{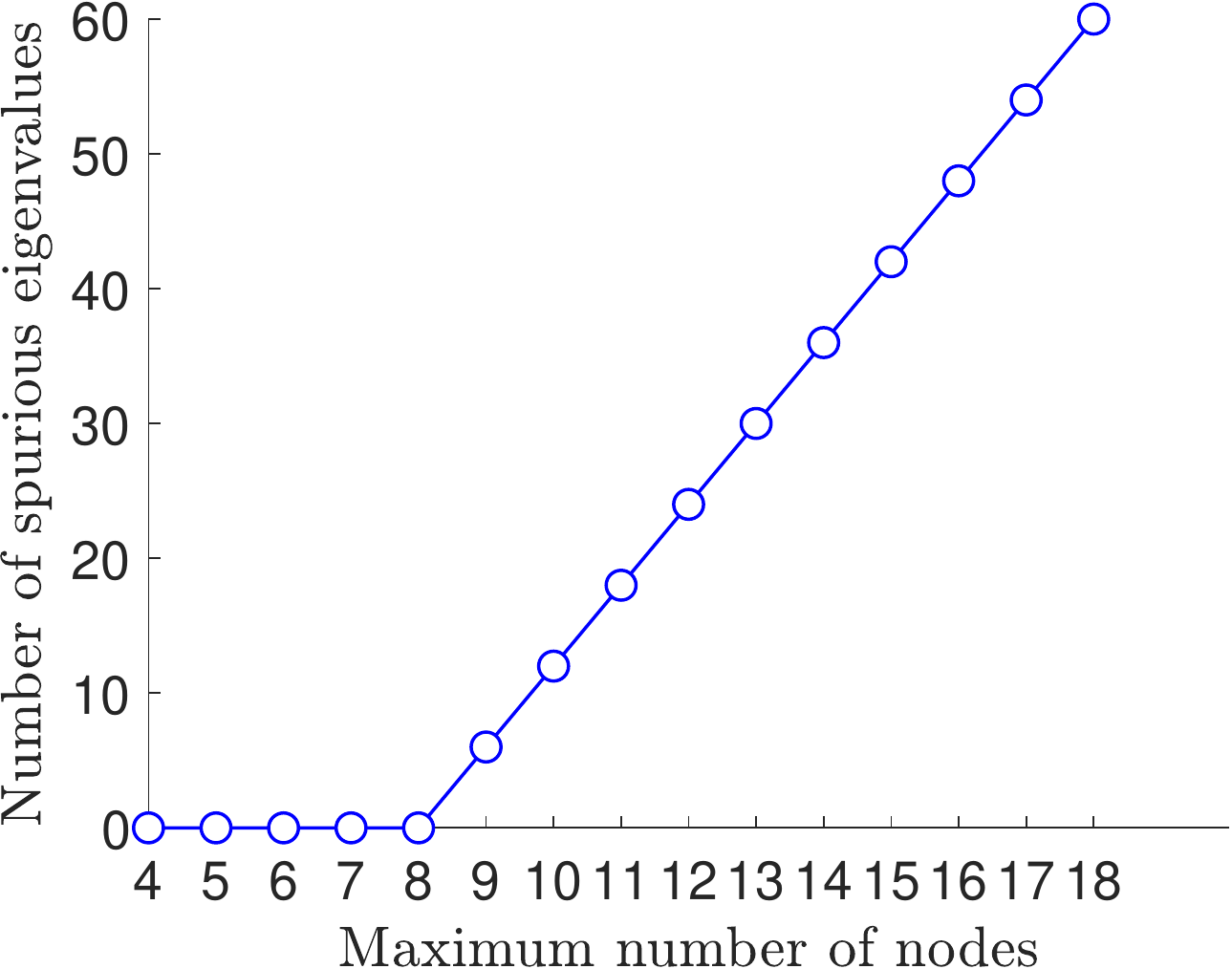}
         \caption{}
         \label{}
     \end{subfigure}
     \hfill
     \begin{subfigure}{0.32\textwidth}
         \centering
         \includegraphics[width=\textwidth]{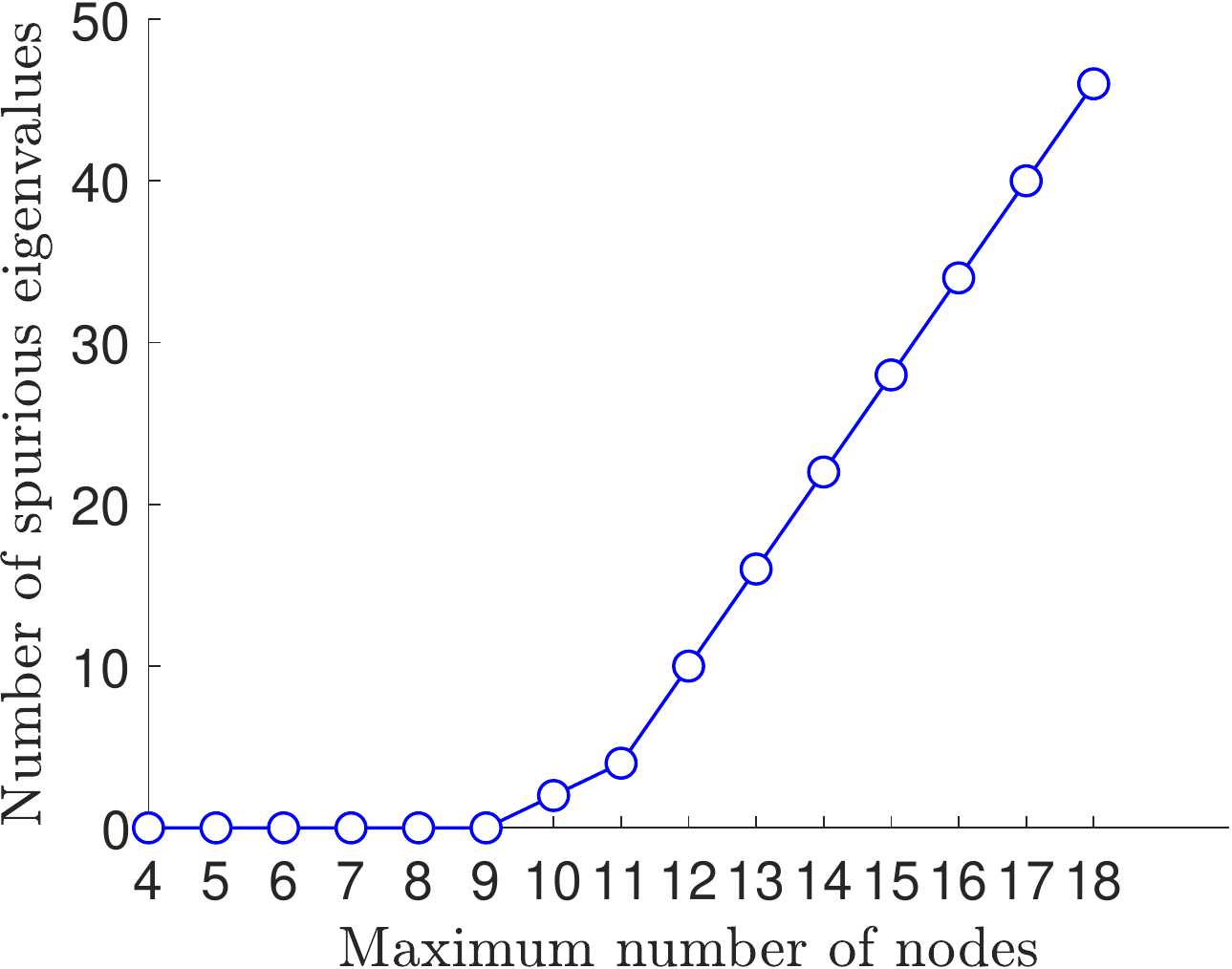}
         \caption{}
         \label{}
     \end{subfigure}
        \caption{Eigenvalue analysis on the meshes shown 
        in Figure~\protect\ref{fig:eig_mesh} with the third-order method.
        (a) $\ell =3$, (b) $\ell=4$ and (c) $\ell=5$. }
        \label{fig:eig_plots_k3}
\end{figure}

We now consider a more general polygonal mesh. Consider the 
unit square, \ac{${\Omega} = (0,1)^2$}, which is
discretized using nine quadrilateral elements. We again solve the
element-eigenvalue problem, $\bm{K}_E \bm{d}_E = \lambda \bm{d}_E$. We choose $\ell =3 ,4 ,5$ and measure the maximum number of spurious eigenvalues of the element
stiffness matrix as we artificially increase the number of 
nodes of the central element. We show a few sample meshes in Figure~\ref{fig:eig_mesh}. In Figure~\ref{fig:eig_plots_k2}, the resulting number of spurious eigenvalues as a function of the number of nodes of an element from the second-order method are plotted for $\ell = 3,4,5$ and similarly the results of the third-order method is plotted in Figure~\ref{fig:eig_plots_k3}. We see that for $k=2$, the spurious modes seem to appear later than in the regular polygons, while for $k=3$ the results are closer to the regular polygonal case. This suggests that the inequalities $N_E \leq 2\ell+1$ and $N_E\leq 2\ell-1$ provide an upper bound for the choice of $\ell$ for $k=2,3$, respectively.  

\section{Numerical results}\label{sec:numerical_results}
We present a series of numerical examples in plane elasticity for second- and third-order serendipity methods. \ac{For these tests we use the inequalities $N_E\leq 2\ell+1$ and $N_E\leq 2\ell-1$ for $k=2$ and $k=3$, respectively.} We examine the errors using the $L^\infty$ and $L^2$ norms, as well as the energy seminorm, and compare the convergence rates of the method with the theoretical estimates of standard VEM. In particular, we use the following discrete measures:
\begin{subequations}\label{eq:error_norms}
\begin{align}
    \|\bm{u} -\bm{u}_h\|_{\bm{L}^\infty(\Omega)} &= \max_{\ac{\bm{x}_i} \in \Omega}{|\bm{u}(\ac{\bm{x}_i}) - \bm{u}_h(\ac{\bm{x}_i})|}, \\
    \|\bm{u} -\bm{u}_h \|_{\bm{L}^2(\Omega)} &= \sqrt{\sum_{E}{\int_{E}{|\bm{u}-\ac{\projS{\bm{u}_h}}}|^2 \,d\bm{x}}},  \\
    \|\bm{u} -\bm{u}_h\|_{a} &= \sqrt{\sum_{E}{\int_{E}{(\overline{\bm{\varepsilon}}-\overline{\projLtwo{\bm{u}_h}})^T \bm{C}{(\overline{\bm{\varepsilon}}-\overline{\projLtwo{\bm{u}_h}}) }\,d\bm{x}}}}.
\end{align}
\end{subequations}
In order to compute the integrals in~\eqref{eq:error_norms}, we adopt the scaled boundary cubature scheme~\cite{chin:2021:cmame}. In the SBC method, an integral over a general polygonal element $E$ is written as the sum of integrals over triangles that are mapped onto the unit square. Let $f$ be any scalar function and
$\vm{c}_i(t)$ be the parametric representation of the edge $e_i$. Define $\ell_i$ to be the signed distance from a fixed point $\vx_0$ to the line containing $e_i$ and  $|e_i|$ denote the length of the $i$-th edge. On using the scaled-boundary parametrization, $\vx = \vm{\varphi}_i (\xi,t) = \vx_0 + \xi(\vm{c}_i(t) - \vx_0)$,
the integral of $f$ over $E$ can be expressed as~\cite{chin:2021:cmame}:
\begin{equation}\label{eq:SBC}
        \int_{E}{f \,d\bm{x}} = \sum_{i=1}^{N_E}\ell_i |e_i| \int_{0}^{1}\int_{0}^1 \xi  f \bigl(\vm{\varphi}_i(\xi,t) \bigr)\, d\xi dt ,
\end{equation}
where in the computations we set $\vx_0$ to be a vertex of the polygon.
To compute the integral over the unit square in~\eqref{eq:SBC}, we use a tensor-product Gauss quadrature rule. 
\subsection{Patch tests}
To test the second- and third-order methods, we first consider the quadratic and cubic displacement patch test. Let $\Omega = (0,1)^2$, $\ac{E_Y}=1$ psi and $\nu=0.3$ be the material properties. For the quadratic patch test, we impose a quadratic displacement field on the boundary and an associated load vector:
\begin{align*}
    &u(\bm{x}) = x^2+3xy+7y^2+5x+2y+8,  \\
    &v(\bm{x}) = 6x^2+3xy+y^2+4x+9y+1 \quad \textrm{on } \partial \Omega, \\ 
    & f(\vm{x}) = \begin{Bmatrix}
    \frac{-E}{1-\nu^2}\left(2+3\nu+\frac{17}{2}(1-\nu)\right) \\ \frac{-E}{1-\nu^2}\left(2+3\nu+\frac{15}{2}(1-\nu)\right)
    \end{Bmatrix}.
\end{align*}

For the cubic patch test we impose a cubic displacement field and load vector:
\begin{align*}
    &u(\bm{x}) = 3x^3+6x^2y+7xy^2+8y^3 + x^2+3xy+y^2+5x+2y+4, \\
    &v(\bm{x}) = 4x^3+7x^2y+8xy^2+11y^3+2x^2+xy+4y^2+8x+9y+11 \quad \textrm{on } \partial \Omega, \\
    &f(\vm{x}) = \begin{Bmatrix}
        \frac{-E}{1-\nu^2}\left(18x+12y+2+\nu(14x+16y+1)+\frac{1-\nu}{2}(36x+28y+7)\right) \\ \frac{-E}{1-\nu^2}\left(\frac{1-\nu}{2}(36x+28y+7)+\nu(14y+3+12x)+16x+66y+8\right)
    \end{Bmatrix}.
\end{align*}

The exact solutions is the extension of the boundary data onto the entire domain $\Omega$. We test the numerical solution for the two methods for four different meshes with 16 elements in each case. First we have a uniform square mesh, second we use a random Voronoi mesh, next we use a Voronoi mesh after applying three Lloyd iterations and finally we use a non-convex mesh. The results for the quadratic test are listed in Table~\ref{tab:quadpatch}, and the cubic test in Table~\ref{tab:cubicpatch}. They show that the errors are near machine precision, which indicate that the second- and third-order method passes the quadratic and cubic patch tests respectively.    
\begin{figure}[H]
     \centering
     \begin{subfigure}{0.22\textwidth}
         \centering
         \includegraphics[width=\textwidth]{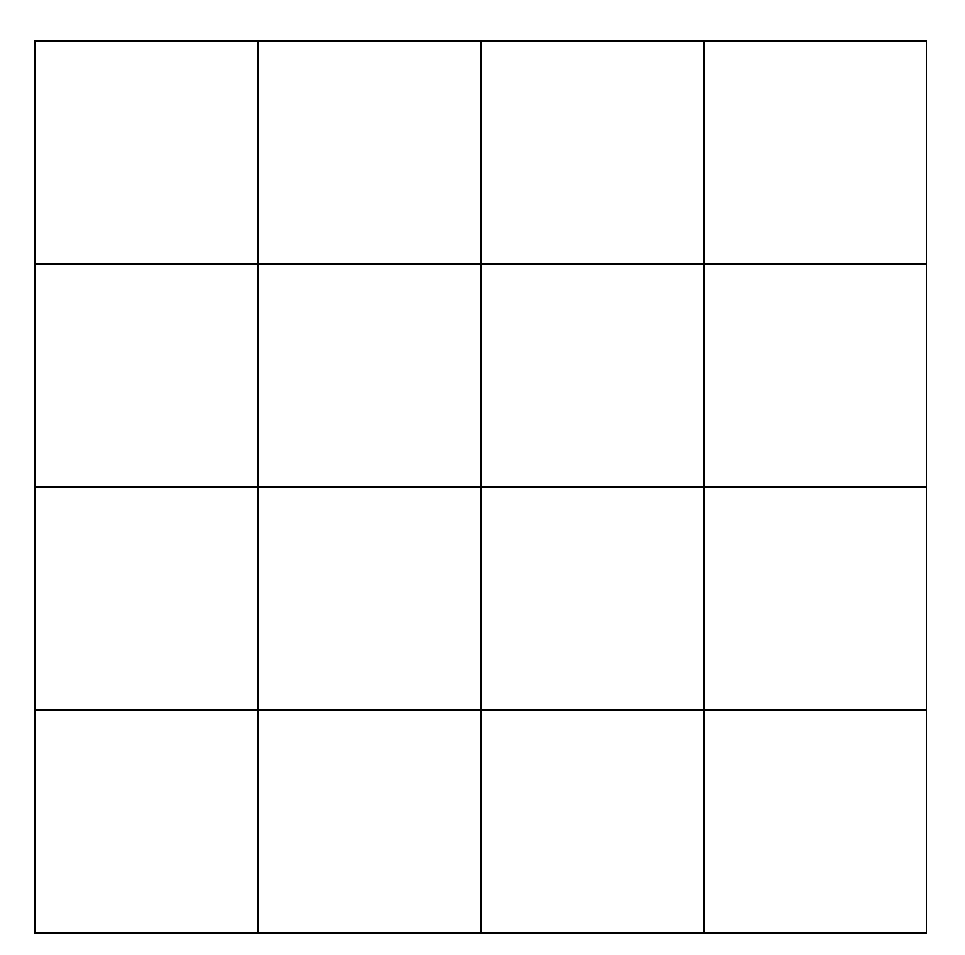}
         \caption{Uniform}
         \label{}
     \end{subfigure}
     \hfill
     \begin{subfigure}{0.22\textwidth}
         \centering
         \includegraphics[width=\textwidth]{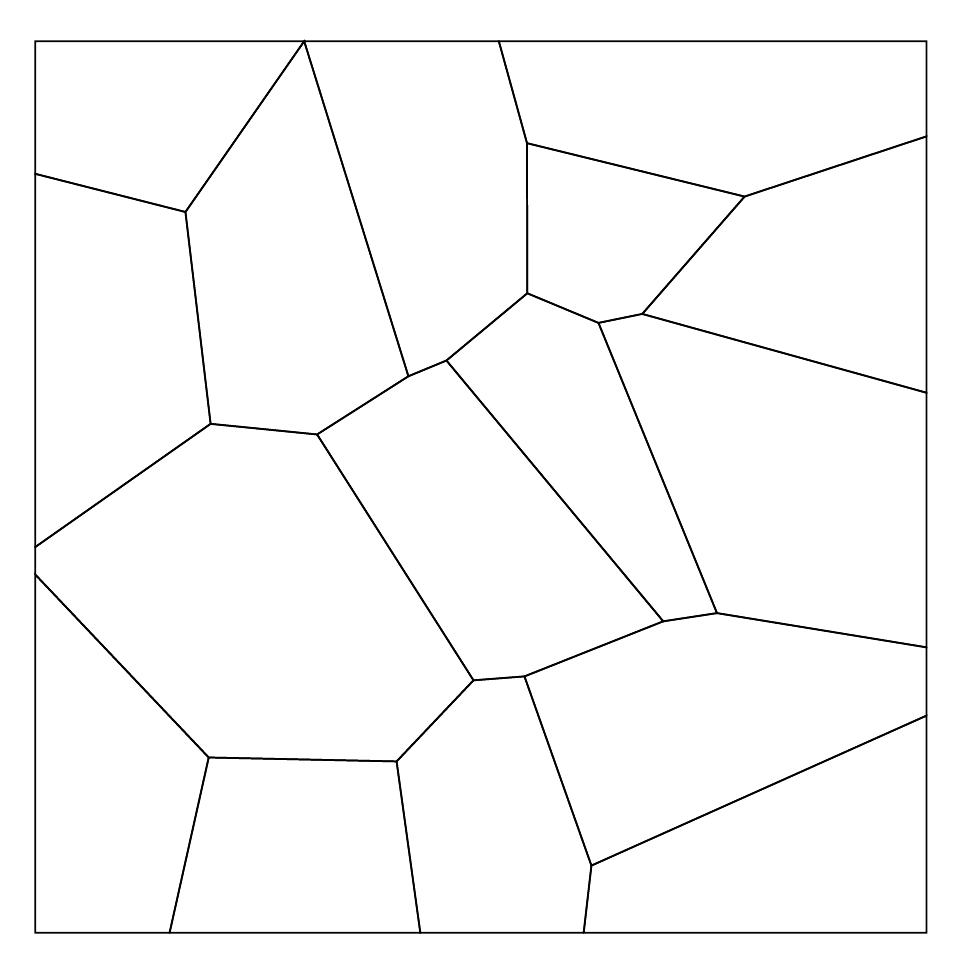}
         \caption{Random}
         \label{}
     \end{subfigure}
     \hfill
     \begin{subfigure}{0.22\textwidth}
         \centering
         \includegraphics[width=\textwidth]{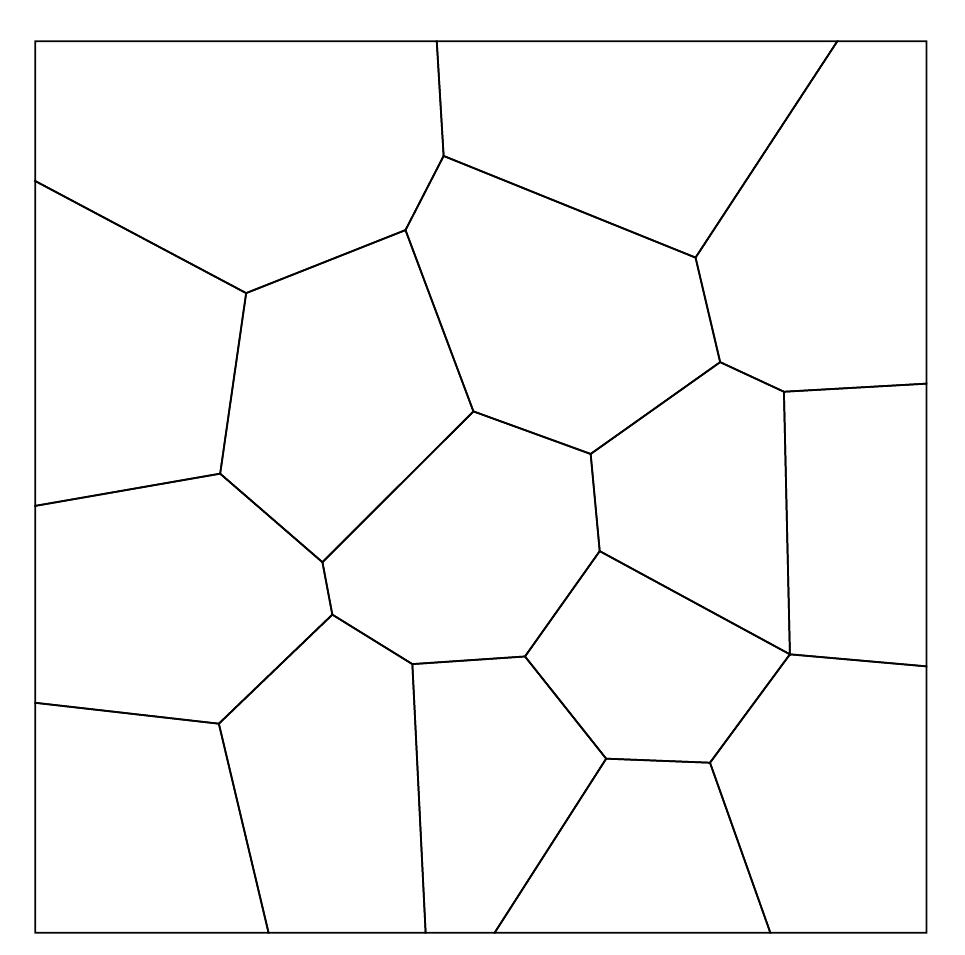}
         \caption{Lloyd iterated}
         \label{}
     \end{subfigure}
     \hfill
     \begin{subfigure}{0.22\textwidth}
         \centering
         \includegraphics[width=\textwidth]{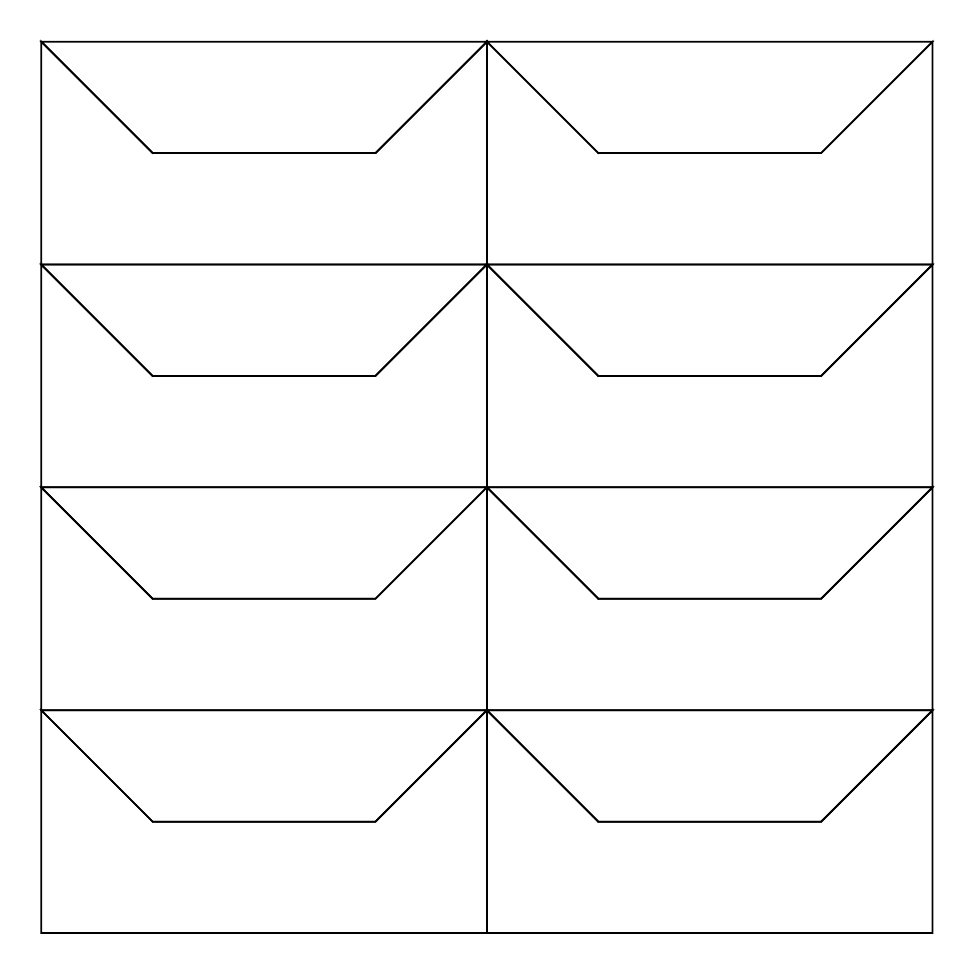}
         \caption{Nonconvex elements}
         \label{}
     \end{subfigure}
        \caption{Sample meshes used for the displacement patch test.}
        \label{fig:mesh_patch}
\end{figure}

\begin{table}[H]
    \centering
    \begin{tabular}{ |c| c | c |c| }
    \hline
     Mesh type & $L^\infty$ error 
     & $L^2$ error & Energy error \\ [0.5ex] 
    \hline
        Uniform & $1 \times10^{-15}$  & $1 \times10^{-15}$ &  $8 \times10^{-14}$
\\
        Random &$6 \times10^{-15}$&  $1 \times10^{-15}$ & $1 \times10^{-14}$ \\
        Lloyd iterated & $7 \times10^{-15}$ & $1 \times10^{-15}$ & $2 \times10^{-14}$\\
        Nonconvex & $4\times 10^{-15}$  & $3 \times 10^{-15}$  & $3 \times 10^{-14}$ \\  
    \hline
    \end{tabular}
    \caption{Errors for the quadratic displacement patch test on different types of meshes.}
    \label{tab:quadpatch}
\end{table}

\begin{table}[H]
    \centering
    \begin{tabular}{ |c| c | c |c| }
    \hline
     Mesh type & $L^\infty$ error 
     & $L^2$ error & Energy error \\ [0.5ex] 
    \hline
        Uniform & $7 \times10^{-15}$  & $2 \times10^{-15}$ &  $2 \times10^{-14}$
\\
        Random &$6 \times10^{-14}$&  $4 \times10^{-14}$ & $3 \times10^{-13}$ \\
        Lloyd iterated & $1 \times10^{-15}$ & $6 \times10^{-15}$ & $4 \times10^{-14}$\\
        Nonconvex & $3\times 10^{-15}$  & $2 \times 10^{-15}$  & $2 \times 10^{-14}$ \\  
    \hline
    \end{tabular}
    \caption{Errors for the cubic displacement patch test on different types of meshes.}
    \label{tab:cubicpatch}
\end{table}

We are also interested in the patch test when Neumann boundary conditions are imposed. Let $\Omega=(0,8) \times (-0.5 , 0.5)$ be a long slender bar with material properties $\ac{E_Y}=1$ psi and $\nu=0.3$. For $k=2$, we construct the following exact solution:
\begin{align*}
    &u(\bm{x}) = xy \ \ \textrm{and} \  \
    v(\bm{x}) = x \quad \textrm{in} \ \Omega, \\ 
    &f(\vm{x}) = \begin{Bmatrix}
    0 \\ -\dfrac{E}{2(1-\nu)}
    \end{Bmatrix},
\end{align*}
where the Dirichlet boundary is imposed along $x=0$, and the remaining boundary conditions on other edges are set to the exact tractions. For $k=3$, we use the cantilever beam under shear end load~\cite{timoshenko1951theory}. We obtain the numerical solutions over a set of three meshes with $16$ elements in each. The results for the quadratic and cubic cases are listed in Tables~\ref{tab:quadpatch_equilibrium} and~\ref{tab:cubicpatch_equilibrium}, respectively. The results show that both the second- and third-order method pass this patch test with errors at worst of
${\cal O}(10^{-11})$.

\begin{table}[H]
    \centering
    \begin{tabular}{ |c| c | c |c| }
    \hline
     Mesh type & $L^\infty$ error 
     & $L^2$ error & Energy error \\ [0.5ex] 
    \hline
        Uniform & $2 \times10^{-11}$  & $1 \times10^{-13}$ &  $1 \times10^{-13}$
\\
        Random &$3 \times10^{-13}$&  $3 \times10^{-13}$ & $4 \times10^{-14}$ \\
        Lloyd iterated & $2 \times10^{-12}$ & $2 \times10^{-12}$ & $4 \times10^{-14}$\\
    \hline
    \end{tabular}
    \caption{Errors for the quadratic equilibrium patch test on different types of meshes.}
    \label{tab:quadpatch_equilibrium}
\end{table}

\begin{table}[H]
    \centering
    \begin{tabular}{ |c| c | c |c| }
    \hline
     Mesh type & $L^\infty$ error 
     & $L^2$ error & Energy error \\ [0.5ex] 
    \hline
        Uniform & $2 \times10^{-11}$  & $2 \times10^{-11}$ &  $1 \times10^{-11}$
\\
        Random &$6 \times10^{-12}$&  $6 \times10^{-12}$ & $6 \times10^{-12}$ \\
        Lloyd iterated & $7 \times10^{-13}$ & $8 \times10^{-13}$ & $2 \times10^{-12}$\\
    \hline
    \end{tabular}
    \caption{Errors for the cubic equilibrium patch test on different types of meshes.}
    \label{tab:cubicpatch_equilibrium}
\end{table}

\subsection{Manufactured exact solutions}
\begin{figure}[H]
     \centering
     \begin{subfigure}{0.48\textwidth}
         \centering
         \includegraphics[width=\textwidth]{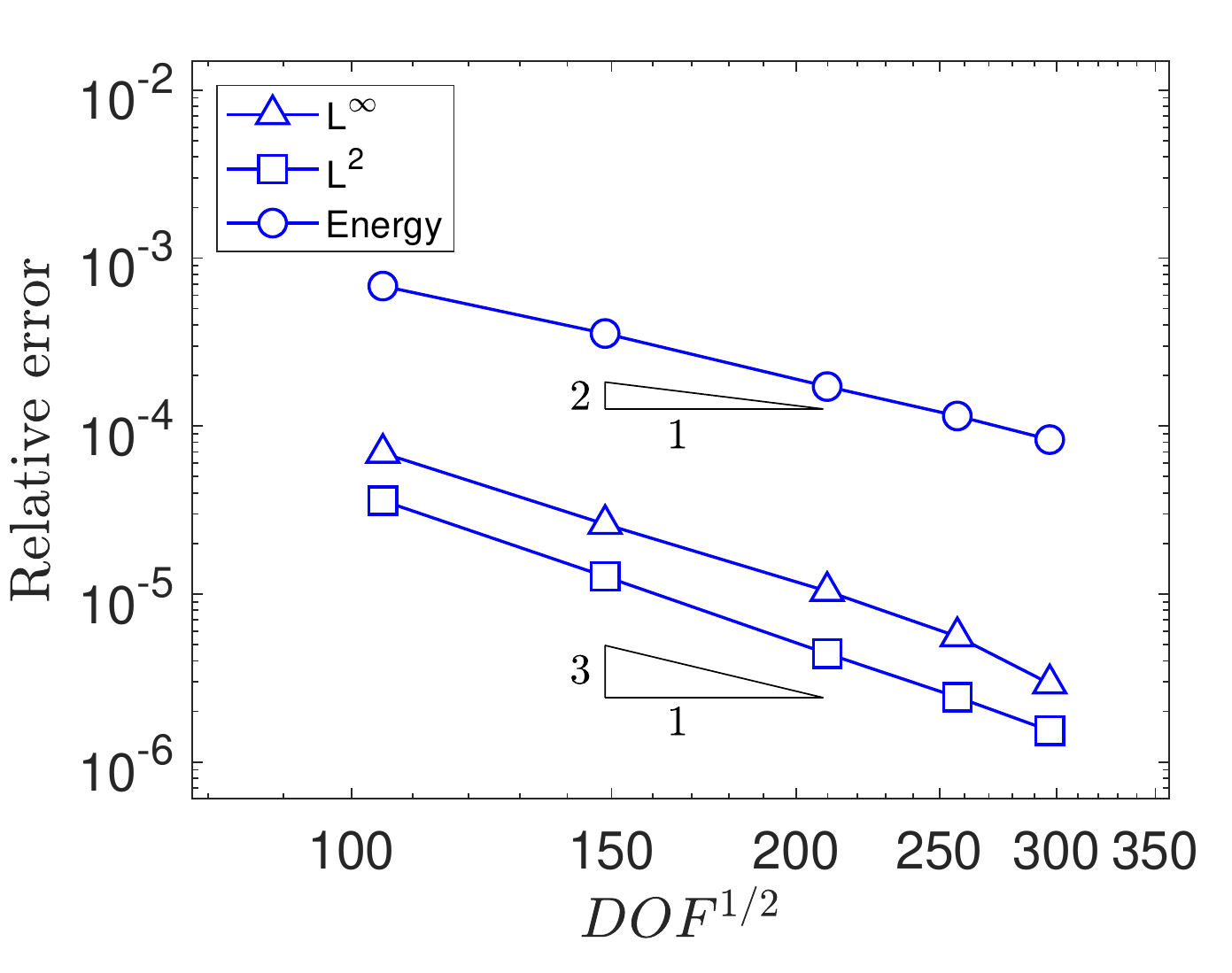}
         \caption{}
     \end{subfigure}
     \hfill
     \begin{subfigure}{0.48\textwidth}
         \centering
         \includegraphics[width=\textwidth]{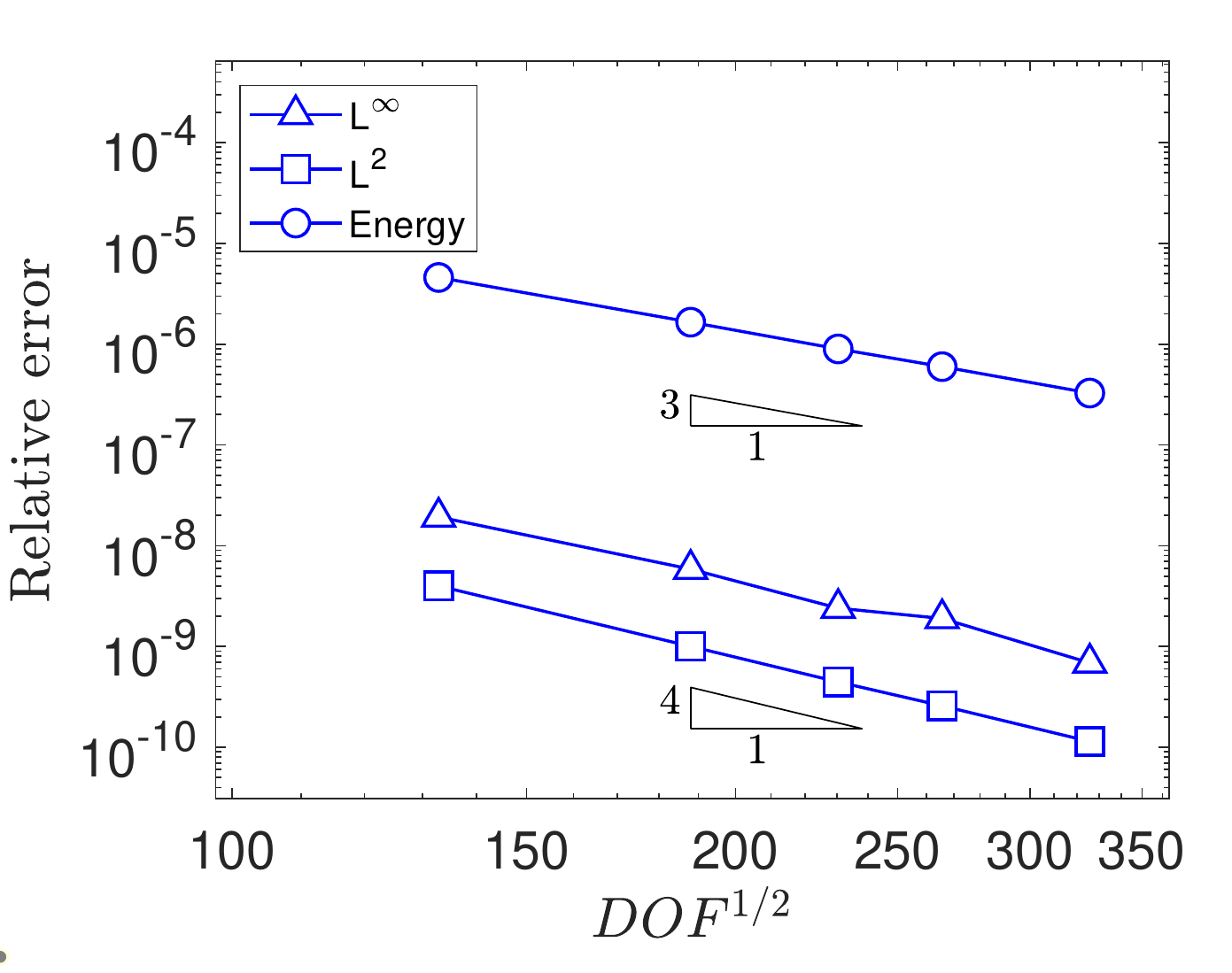}
         \caption{}
     \end{subfigure}
        \caption{Convergence curves for first manufactured solution on convex polygonal meshes with (a) $k=2$ and (b) $k=3$.}
        \label{fig:manufactured_1}
\end{figure}

\begin{figure}[H]
     \centering
     \begin{subfigure}{0.48\textwidth}
         \centering
         \includegraphics[width=\textwidth]{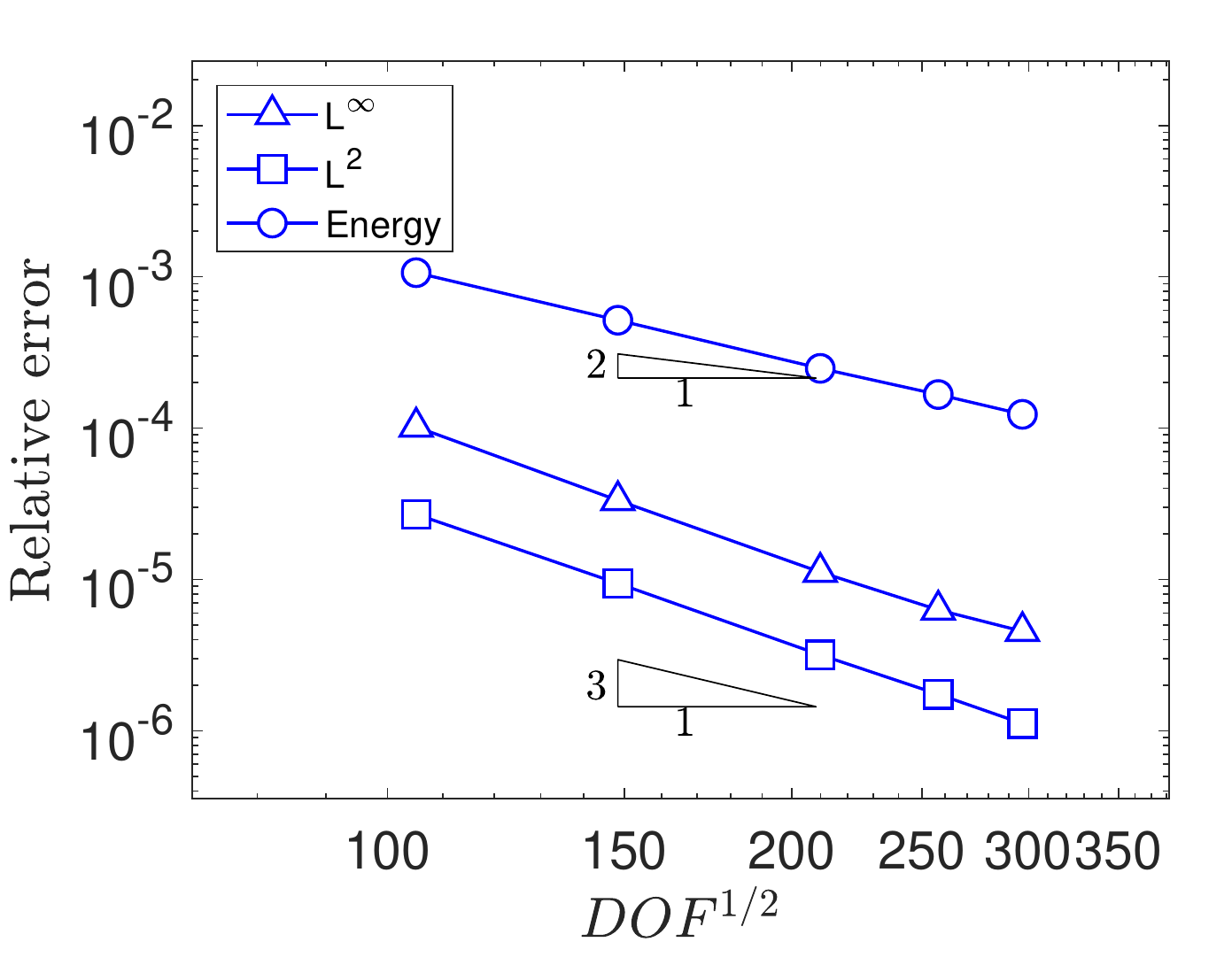}
         \caption{}
     \end{subfigure}
     \hfill
     \begin{subfigure}{0.48\textwidth}
         \centering
         \includegraphics[width=\textwidth]{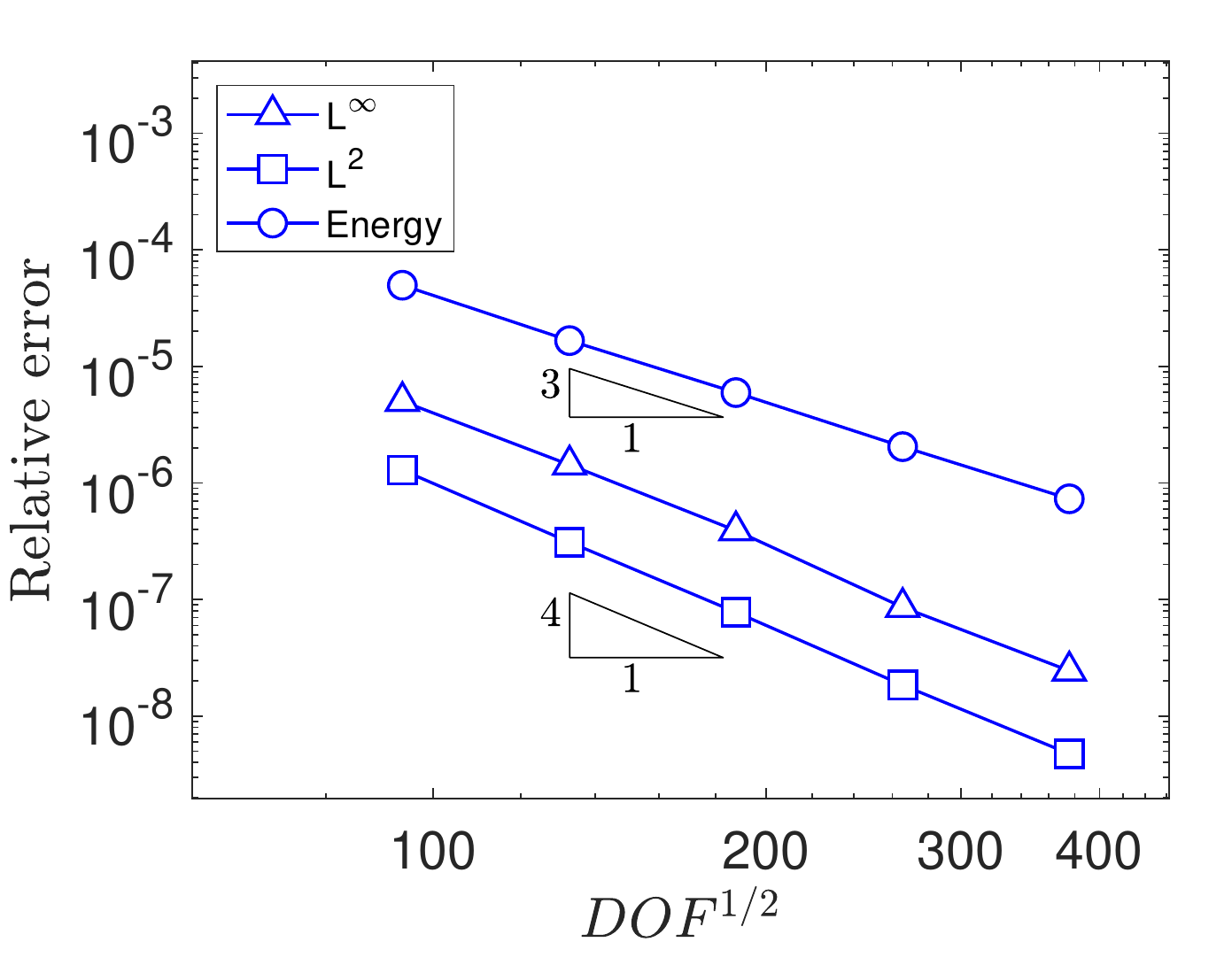}
         \caption{}
     \end{subfigure}
        \caption{Convergence curves for second manufactured solution on convex polygonal meshes with (a) $k=2$ and (b) $k=3$. }
        \label{fig:manufactured_2}
\end{figure}

We consider two manufactured problem as given in~\cite{enhanced:VEM} with known exact polynomial and nonpolynomial solutions over the unit square under plane stress conditions. The material properties are: $\ac{E_Y}=2.5$ psi and $\nu=0.25$. The exact solution and the associated loading for the first problem are: 
\begin{align*}
    &u(\vm{x}) = -\frac{x^6}{80} + \frac{x^4y^2}{2} - \frac{13}{16}x^2y^4 + \frac{3}{40}y^6 \ \ \textrm{and}  \ \ v(\vm{x}) = \frac{xy^5}{2}-\frac{5}{12}x^3y^3,\\
    &f(\vm{x})= \begin{Bmatrix}
         0 \\
         0
    \end{Bmatrix},
\end{align*}
and for the second problem are:
\begin{align*}
    &u(\vm{x}) = x\sin(\pi x) \sin(\pi y) \ \ \textrm{and}  \ \ v(\vm{x}) = y\sin(\pi x)\sin(\pi y),\\
    &f(\vm{x})= \begin{Bmatrix}
        \frac{11}{3}\pi^2 x \sin(\pi x)\sin(\pi y) - \frac{5}{3}\pi^2 y \cos(\pi x)\cos(\pi y) -7\pi\cos(\pi x)\sin(\pi y) \\
        \frac{11}{3}\pi^2 y \sin(\pi x)\sin(\pi y) - \frac{5}{3}\pi^2 x \cos(\pi x)\cos(\pi y) -7\pi\cos(\pi y)\sin(\pi x)
    \end{Bmatrix}.
\end{align*}
We include the results for both these tests
in Figures~\ref{fig:manufactured_1} and~\ref{fig:manufactured_2}. In both figures, we plot the discrete errors as a function of the square root of the number of degrees of freedom. From the plots, we observe that the convergence rates for $k=2,3$ in the $L^2$ and energy seminorm are in agreement with the theoretical rates. This shows that the stabilization-free virtual element method can reproduce the results from~\cite{enhanced:VEM}.

\subsection{Beam subjected to transverse sinusoidal loading}
\begin{figure}
     \centering
     \begin{subfigure}{\textwidth}
         \centering
         \includegraphics[width=\textwidth]{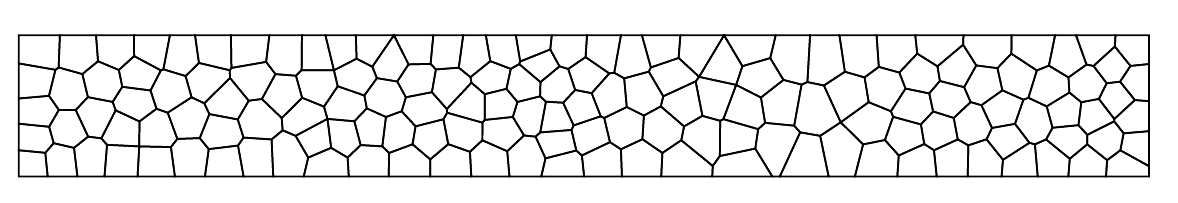}
         \caption{}
     \end{subfigure}
     \vfill
     \begin{subfigure}{\textwidth}
         \centering
         \includegraphics[width=\textwidth]{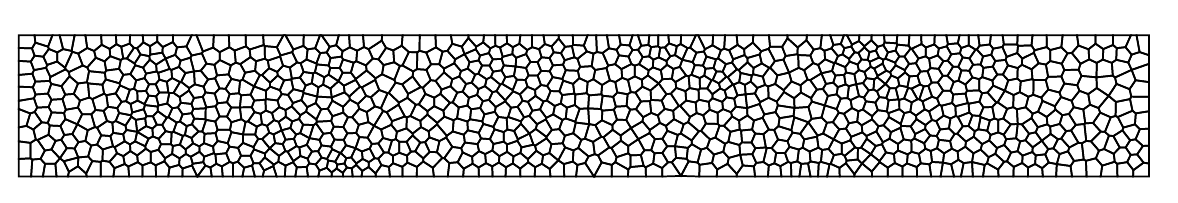}
         \caption{}
     \end{subfigure}
     \vfill
     \begin{subfigure}{\textwidth}
         \centering
         \includegraphics[width=\textwidth]{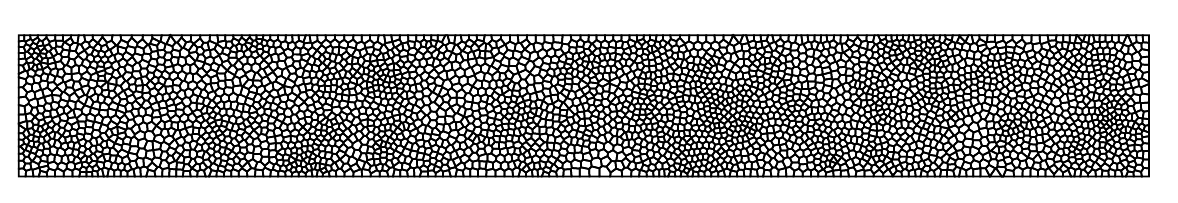}
         \caption{}
     \end{subfigure}
        \caption{Polygonal meshes for the loaded beam problem. (a) 150 elements, (b) 1000 elements and (c) 3500 elements.  }
        \label{fig:beammesh}
\end{figure}

\begin{figure}
     \centering
     \begin{subfigure}{0.48\textwidth}
         \centering
         \includegraphics[width=\textwidth]{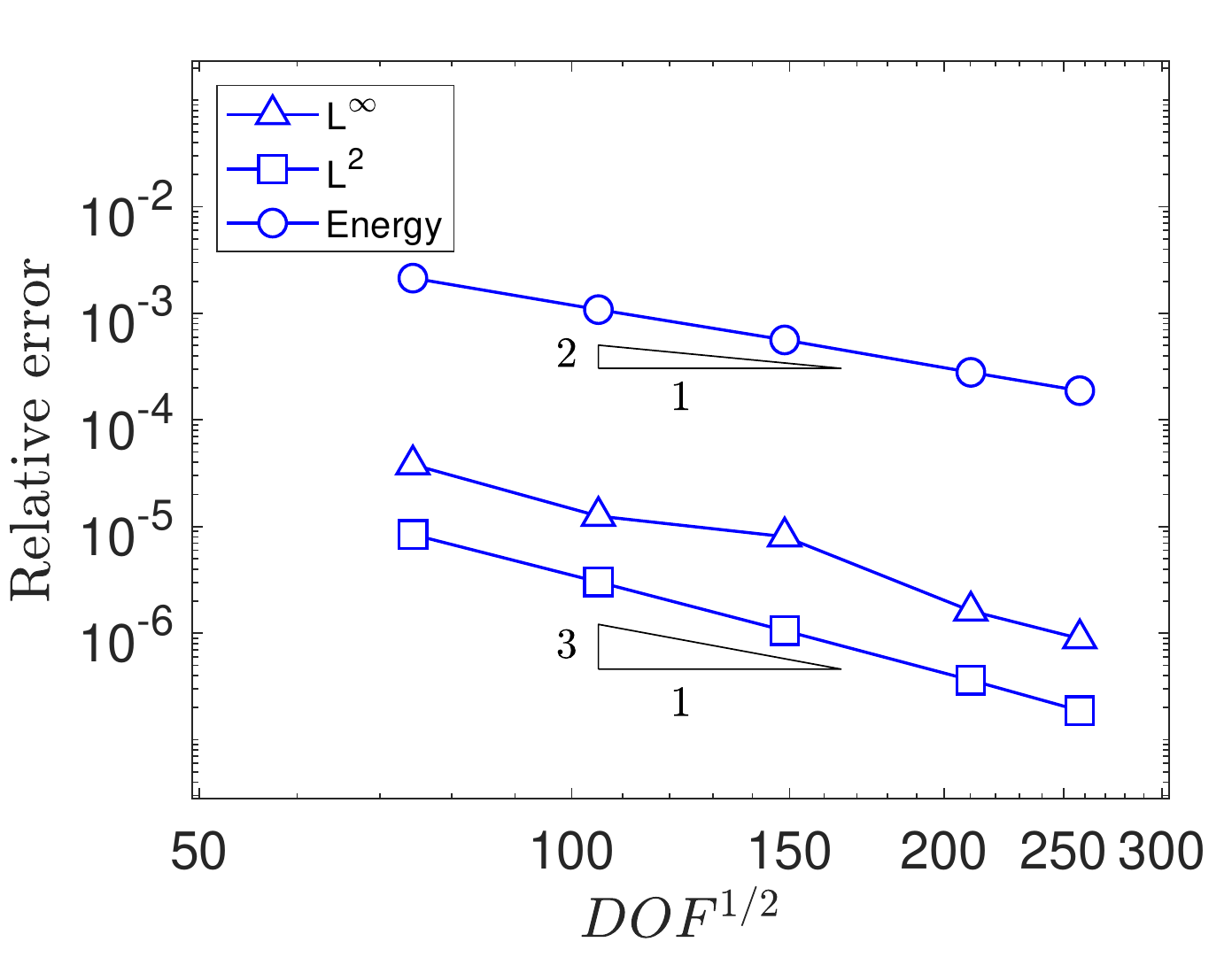}
         \caption{}
     \end{subfigure}
     \hfill
     \begin{subfigure}{0.48\textwidth}
         \centering
         \includegraphics[width=\textwidth]{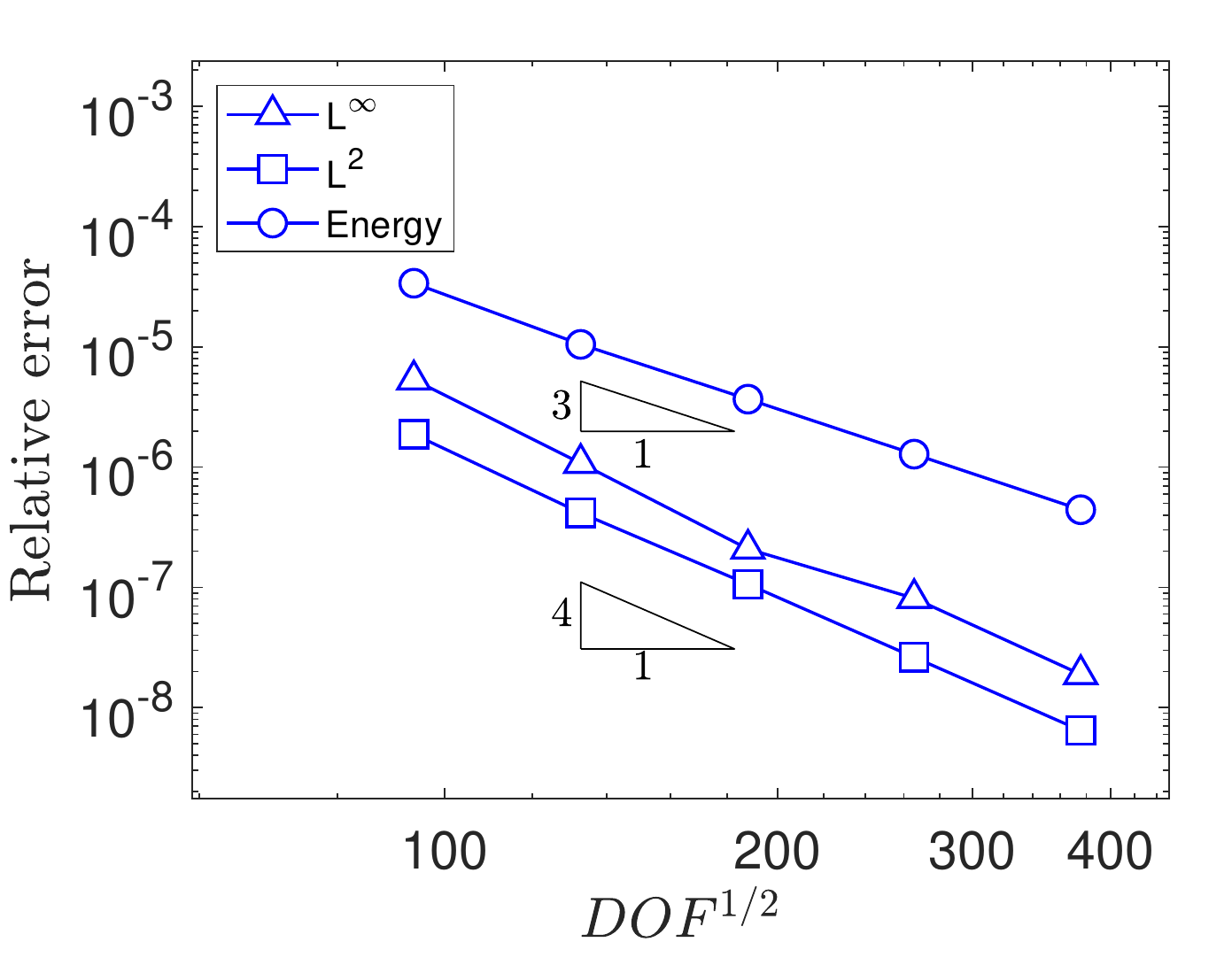}
         \caption{}
     \end{subfigure}
        \caption{Convergence curves of serenedipity VEM on convex meshes for sinusoidal loaded beam 
        problem. (a) $k=2$ and (b) $k=3$.}
        \label{fig:sinusoidal_load}
\end{figure}

We consider the problem of a simply-supported beam subjected to a transversely sinusoidal load~\cite{Sadd2009Elasticity}. The material properties are chosen as: $\ac{E_Y}=2\times 10^5$ psi and $\nu =0.3$, and plane stress conditions are assumed. The beam has length $L=8$ inch, height $D=1$ inch and unit thickness. We apply a sinusoidal load $P=-100\sin(\frac{\pi x}{L})$ lb along the top edge, and along the two side edges we prescribe shear stresses to keep the beam in equilibrium. This problem does not have a closed-form solution; however, it can be shown that a generalized solution (one that satisfies some of the boundary conditions in an average sense) can be found with a Fourier series Airy stress function. In~\cite{Sadd2009Elasticity}, the solution for this simply-supported beam is given as:
\begin{align*}
    u(\bm{x})= &-\frac{\beta}{E}\cos(\beta x)\left\{A(1+\nu)\sinh(\beta y)+B(1+\nu)\cosh(\beta y)\right. \\ &+C\left[(1+\nu)\beta y \sinh(\beta y)+2\cosh(\beta y)\right] \\ &+D\left[(1+\nu)\beta y\cosh(\beta y)+2\sinh(\beta y)\right]\left.\right\}+u_0,\\
    v(\bm{x}) =&-\frac{\beta}{E}\sin(\beta x)\left\{A(1+\nu)\cosh(\beta y)+B(1+\nu)\sinh(\beta y)\right. \\ &+C\left[(1+\nu)\beta y \cosh(\beta y)-(1-\nu)\sinh(\beta y)\right] \\ &+D\left[(1+\nu)\beta y\sinh(\beta y)-(1-\nu)\cosh(\beta y)\right]\left.\right\},
\end{align*}
where the constants $A, B, C, D, \beta, u_0$ are detailed in~\cite{Sadd2009Elasticity}. In Figure~\ref{fig:beammesh}, we show a few sample meshes for the beam, and in Figure~\ref{fig:sinusoidal_load} we show the convergence results. From these figures, we observe that  optimal convergence rates in Sobolev norms are achieved for both $k = 2$ and $k = 3$.

\begin{figure}
     \centering
     \begin{subfigure}{\textwidth}
         \centering
         \includegraphics[width=\textwidth]{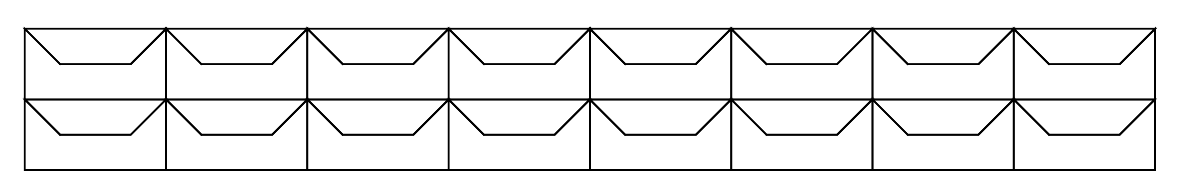}
         \caption{}
     \end{subfigure}
     \vfill
     \begin{subfigure}{\textwidth}
         \centering
         \includegraphics[width=\textwidth]{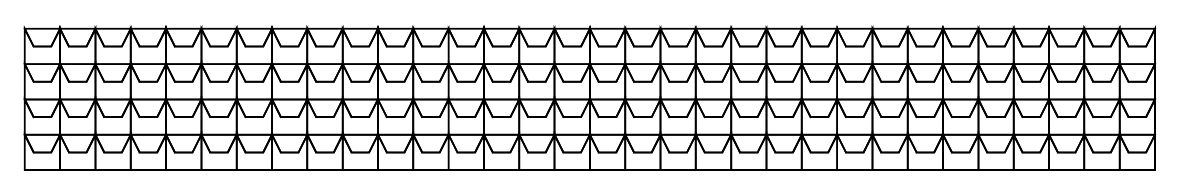}
         \caption{}
     \end{subfigure}
     \vfill
     \begin{subfigure}{\textwidth}
         \centering
         \includegraphics[width=\textwidth]{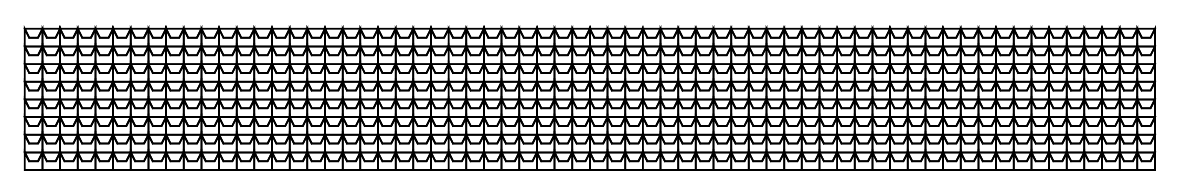}
         \caption{}
     \end{subfigure}
        \caption{Nonconvex polygonal meshes for the loaded beam problem. (a) 32 elements, (b) 256 elements and (c) 1024 elements.  }
        \label{fig:beammesh_non_convex}
\end{figure}

\begin{figure}
     \centering
     \begin{subfigure}{0.48\textwidth}
         \centering
         \includegraphics[width=\textwidth]{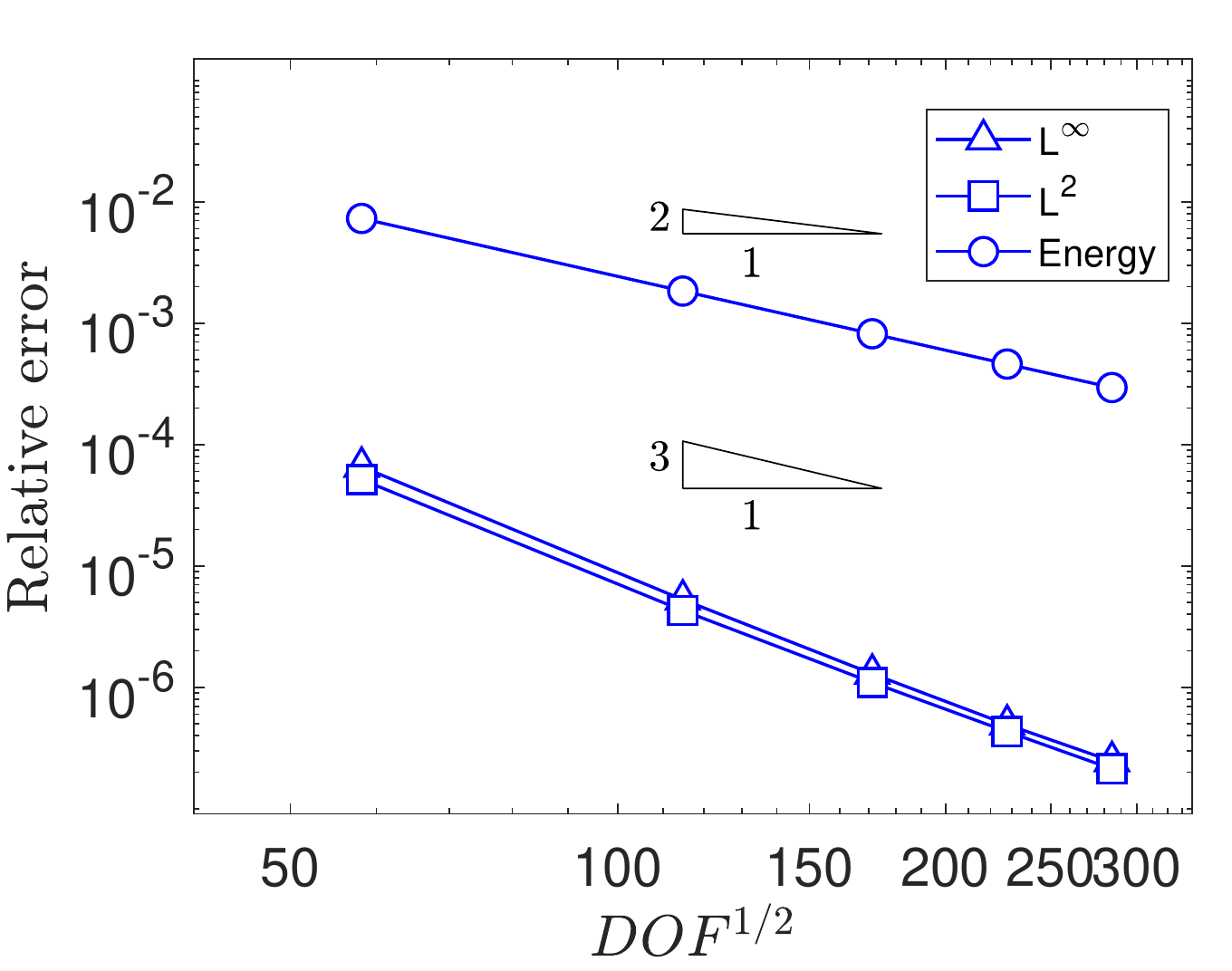}
         \caption{}
     \end{subfigure}
     \hfill
     \begin{subfigure}{0.48\textwidth}
         \centering
         \includegraphics[width=\textwidth]{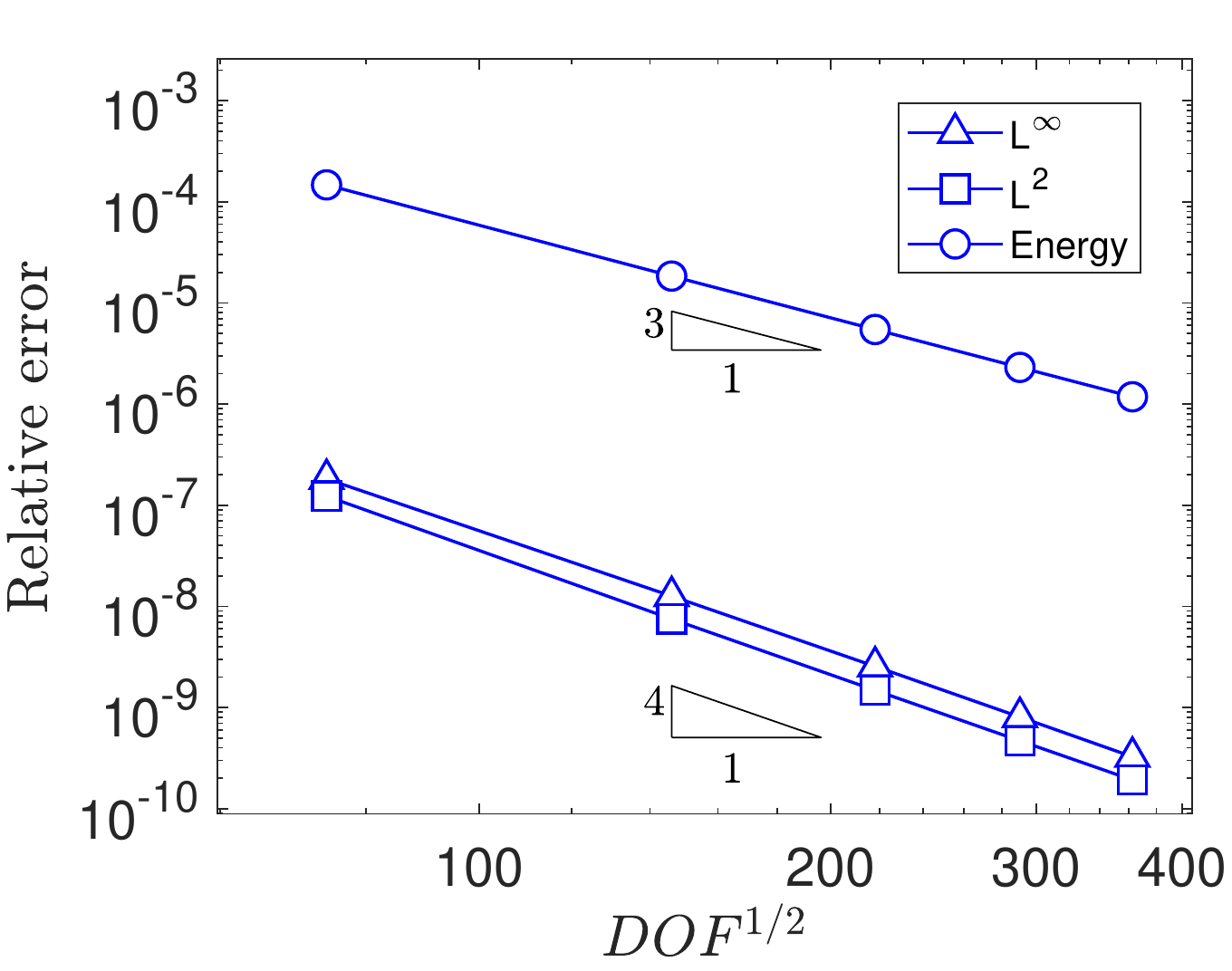}
         \caption{}
     \end{subfigure}
        \caption{Convergence curves of serendipity VEM on nonconvex meshes for sinusoidal loaded beam problem. (a) $k=2$ and (b) $k=3$.}
        \label{fig:sinusoidal_load_non_convex}
\end{figure}

We also test this problem with nonconvex meshes. We start with a uniform rectangular mesh, then we split each element into a convex quadrilateral and a non-convex hexagonal element. We show a few sample meshes in Figure~\ref{fig:beammesh_non_convex}. In Figure~\ref{fig:sinusoidal_load_non_convex}, the results show that the errors on nonconvex meshes still retains the optimal convergence rate.  

\subsection{Infinite plate with a circular hole under uniaxial tension}
\begin{figure}[H]
     \centering
     \begin{subfigure}{0.48\textwidth}
         \centering
         \includegraphics[width=\textwidth]{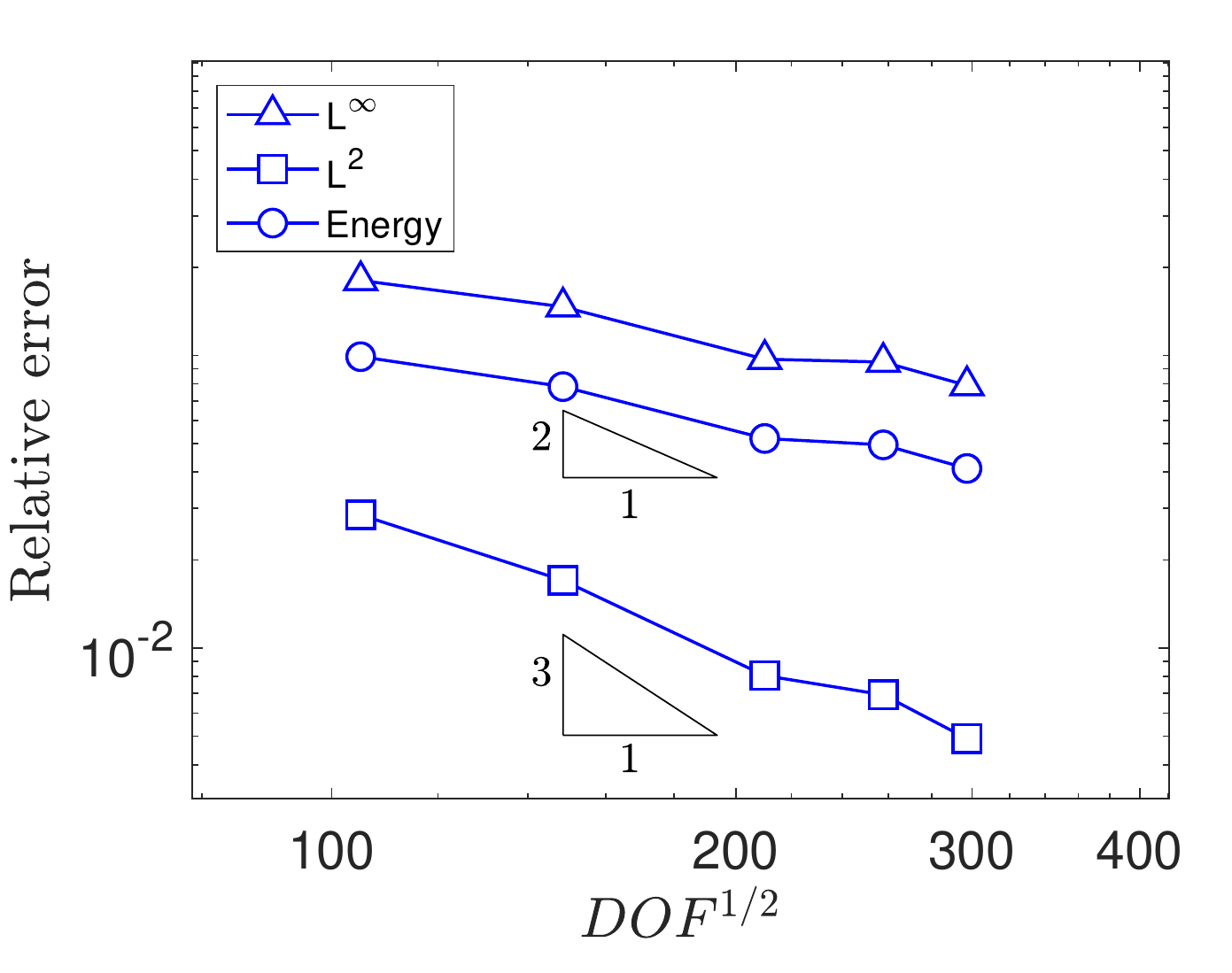}
         \caption{}
     \end{subfigure}
     \hfill
     \begin{subfigure}{0.48\textwidth}
         \centering
         \includegraphics[width=\textwidth]{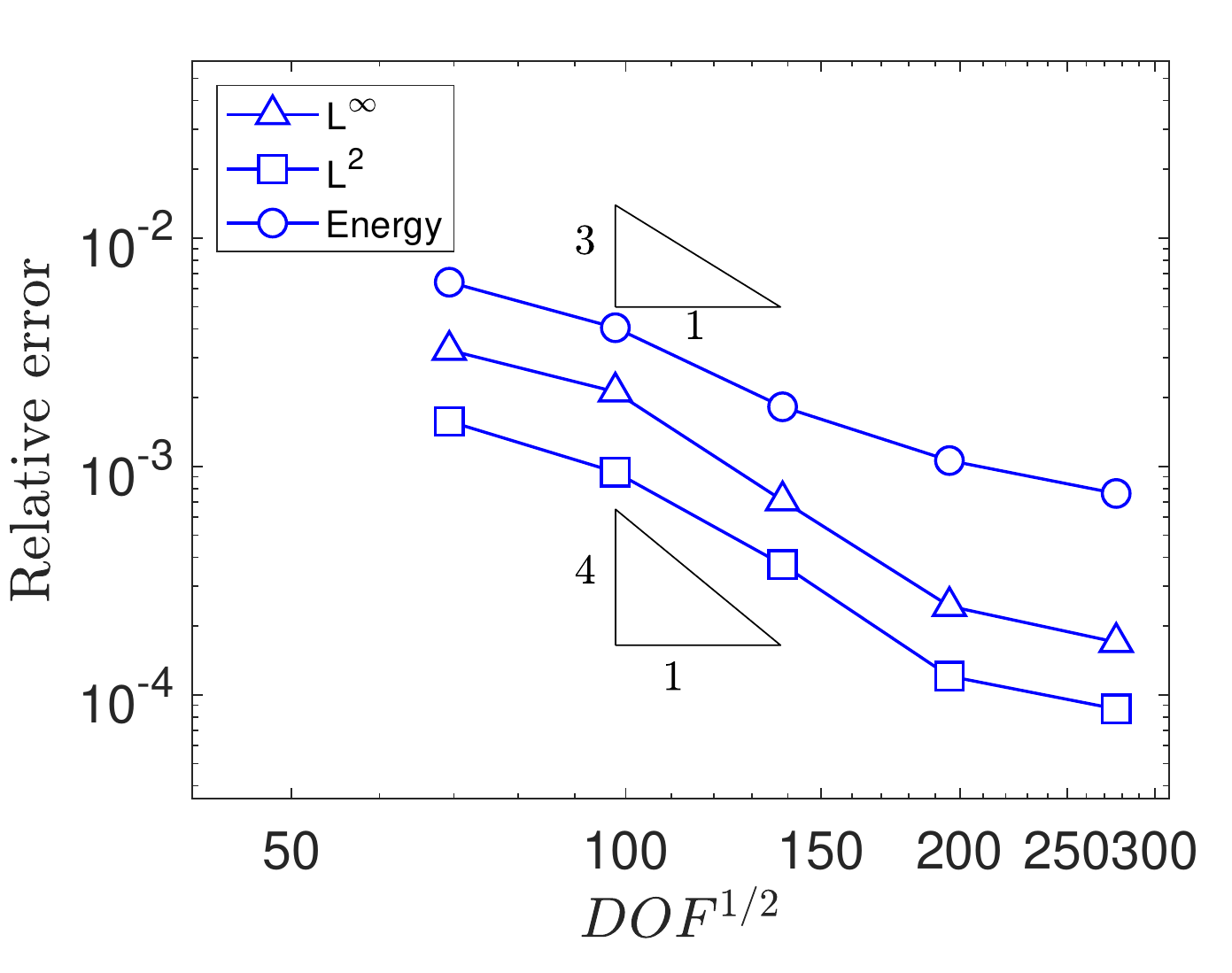}
         \caption{}
     \end{subfigure}
        \caption{Convergence curves of serendipity VEM for plate with a circular hole problem.
        (a) $k=2$ and (b) $k=3$.}
        \label{fig:plate_hole}
\end{figure}

Finally, we consider the problem of an infinite plate with a circular hole under uniaxial tension. The hole is subject to traction-free condition, while a far field uniaxial tension $ \sigma_0=1$ psi, is applied to the plate in the $x$-direction. We use the material properties $\ac{E_Y}=2\times 10^5$ psi and $\nu = 0.3$, with a hole radius $a=1$ inch. Due to symmetry, we model a
quarter of the finite plate ($L = 5$ inch), with exact boundary tractions prescribed as data. Plane strain conditions are assumed. It is known from~\cite{Artioli:2020:cmame,daveiga:2019:esiam}, that standard VEM methods with order $k\geq 2$ will suffer from loss of convergence rates when approximating domains with curved edges. We see this result in Figure~\ref{fig:plate_hole}, where both the second- and third-order methods failed to attain the optimal convergence rates. With this result, it is natural to look into the extension of stabilization free methods onto elements with curved edges. 
\begin{figure}[!h]
     \centering
     \begin{subfigure}{0.32\textwidth}
         \centering
         \includegraphics[width=\textwidth]{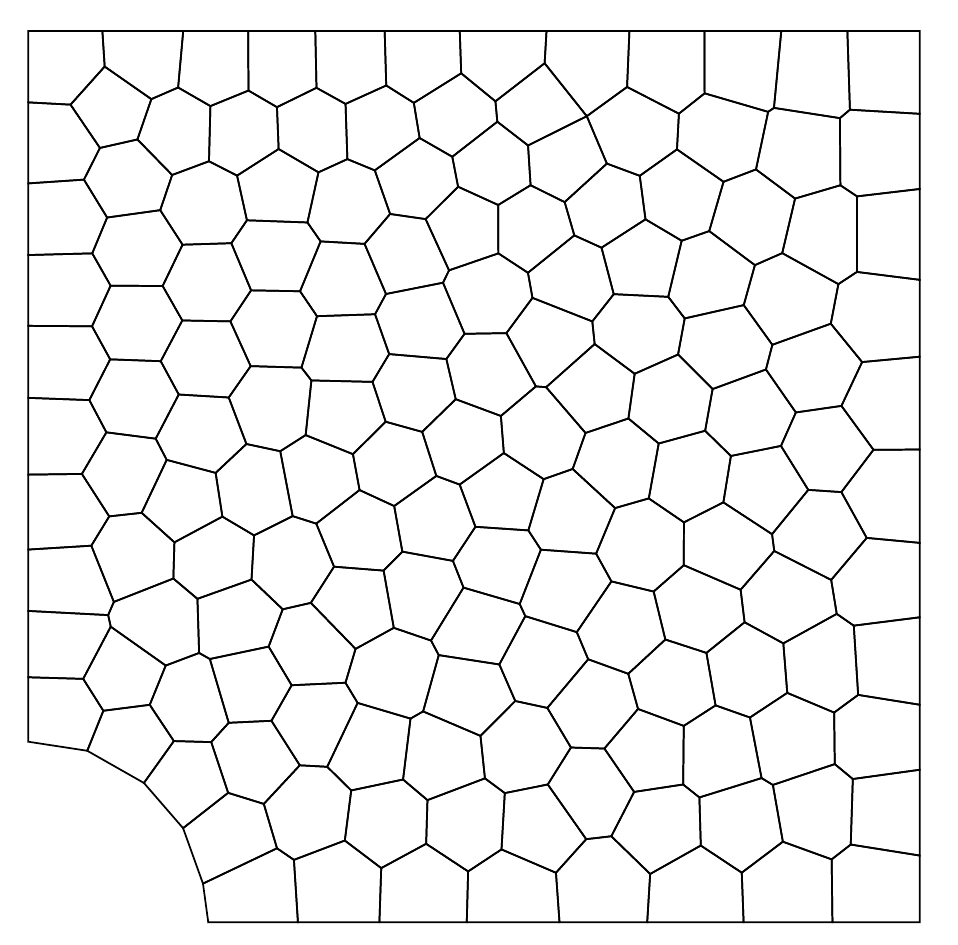}
         \caption{}
     \end{subfigure}
     \hfill
     \begin{subfigure}{0.32\textwidth}
         \centering
         \includegraphics[width=\textwidth]{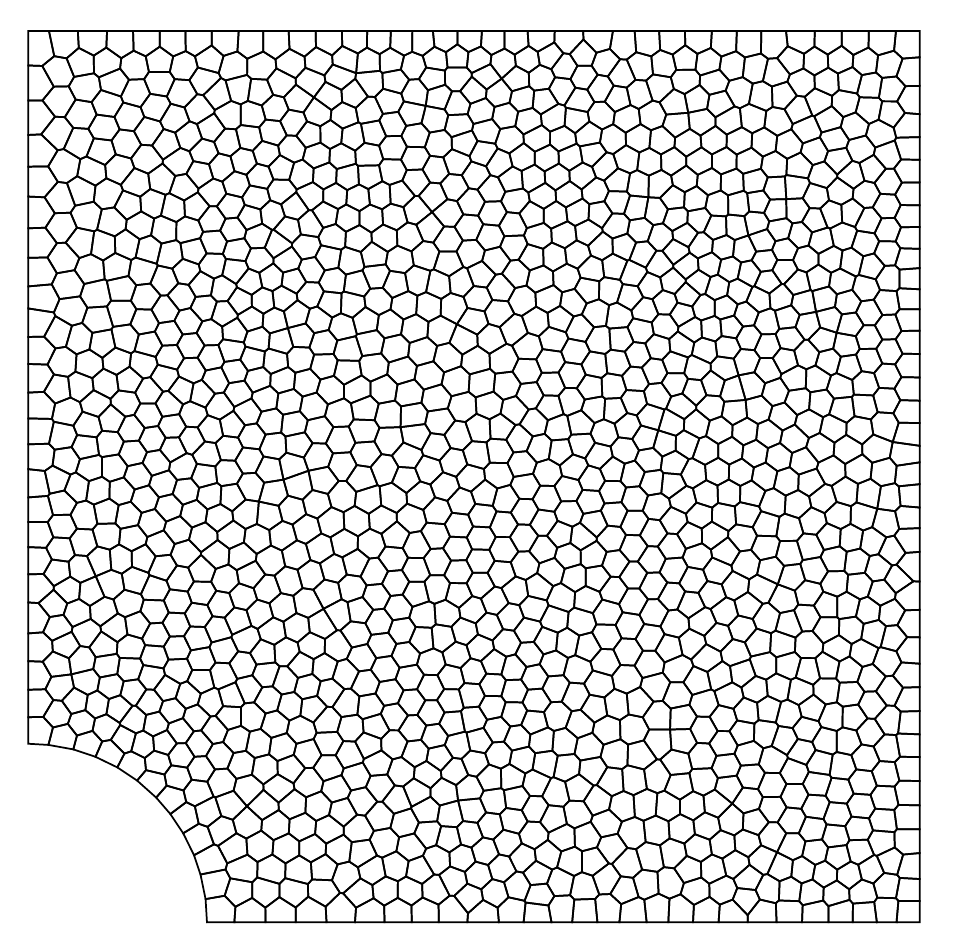}
         \caption{}
     \end{subfigure}
     \hfill
     \begin{subfigure}{0.32\textwidth}
         \centering
         \includegraphics[width=\textwidth]{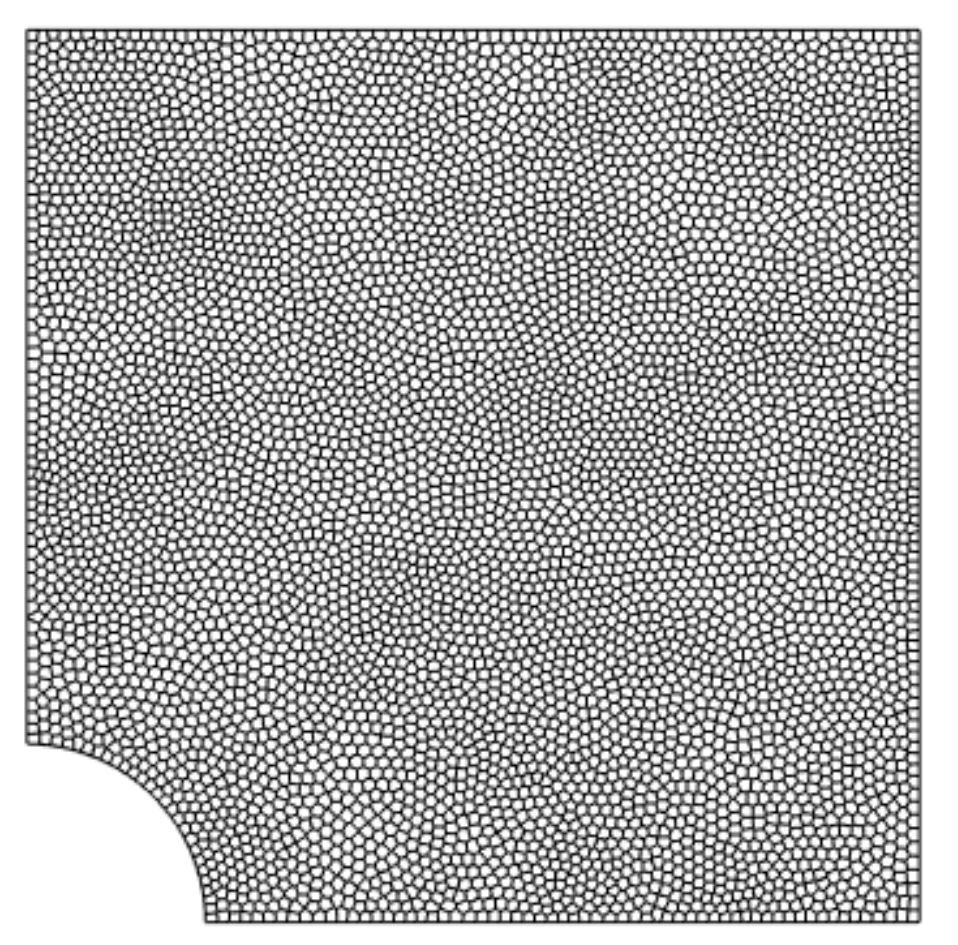}
         \caption{}
     \end{subfigure}
        \caption{Polygonal meshes for the plate with a circular hole problem. (a) 250 elements, (b) 1500 elements, and (c) 6000 elements.}
        \label{fig:holemesh}
\end{figure}

\section{Conclusions}\label{sec:conclusions}
In this paper, we studied a higher order (serendipity) extension of the 
stabilization-free virtual element method~\cite{berrone2021lowest,chen:2022:arxiv} 
for plane elasticity. To establish a stabilization-free method for solid continua, 
we constructed an enlarged VEM space 
that included higher order polynomial approximations of the strain field. To 
eliminate additional degrees of freedom we incorporated the serendipity approach 
into the VEM space~\cite{Veiga:2016:caf}. On each polygonal element we chose the 
degree $\ell$ of vector polynomials such that the element stiffness \ac{has the} correct 
rank. We set up the construction of the necessary projections and stiffness 
matrices, and then solved several
problems from plane elasticity using a second- and third-order method. 
For the patch test, we recovered the displacement and stress fields to near 
machine-precision. 
From an element-eigenvalue analysis, we numerically
examined a suitable choice of $\ell$ that was sufficient to 
ensure that the element stiffness matrix had no spurious zero-energy modes,
and hence the element was stable.
For a few manufactured problems and the cantilever beam problem 
under sinusoidal top load,
we found that the convergence rates of the second- and third-order
stabilization-free VEM in the $L^2$ norm and energy seminorm were in agreement with
standard VEM theoretical results. However, consistent with expectations, we have 
verified that the serendipity virtual element method on affine edges has reduced 
convergence rates for domains with curved edges~\cite{daveiga:2019:esiam}. As part 
of future work, extensions of stabilization-free virtual element to general curved 
domains and nonlinear plane elasticity are of interest.  

\section*{Acknowledgements}
The authors are grateful to Alessandro Russo for many helpful discussions.

\end{document}